\documentclass{amsart}
\usepackage{graphicx} 
\usepackage[colorlinks=true]{hyperref}
\usepackage{xcolor}
\usepackage{amsmath}
\usepackage{amsfonts}
\usepackage{amsthm}
\usepackage{placeins}
\usepackage{enumitem}
\usepackage{amssymb,mathrsfs}
\usepackage{float}
\usepackage{comment}
\usepackage{subcaption} 
\usepackage{}

\newtheorem{theorem}{Theorem}[section]
\newtheorem{lemma}[theorem]{Lemma}
\newtheorem{corollary}[theorem]{Corollary}
\newtheorem{proposition}[theorem]{Proposition}

\theoremstyle{definition}
\newtheorem{definition}[theorem]{Definition}
\newtheorem{example}[theorem]{Example}

\newtheorem{notation}[theorem]{Notation}

\theoremstyle{remark}
\newtheorem{remark}[theorem]{Remark}

\numberwithin{equation}{section}

\newcommand{\cH}{\mathcal{H}}
\newcommand{\calS}{\mathcal{S}}

\newcommand{\CC}{\mathbb{C}}

\newcommand{\calC}{\mathcal{C}}
\newcommand{\calN}{\mathcal{N}}
\newcommand{\DD}{\mathbb{D}}
\newcommand{\NN}{\mathbb{N}}
\newcommand{\Dbar}{\overline{\mathbb{D}}}
\newcommand{\eps}{\varepsilon}

\newcommand{\Fndca}{F_{n,d,c,a}}
\newcommand{\Fnca}{F_{n,c,a}}

\newcommand{\KFndca}{K(F_{n,d,c,a})}
 \newcommand{\JFndca}{J(F_{n,d,c,a})}
 \newcommand{\Mna}{M_{n,a}}
 \newcommand{\Mnda}{M_{n,d,a}}

\newcommand{\Vand}{\mathcal{V}_{n,d,a}}

\newcommand{\Cnda}{\calC_{n,d,a}}

\newcommand{\Ca}{\calC_{a}}
\newcommand{\ttbig}{\texttt{big}}
\newcommand{\tttiny}{\texttt{tiny}}
\newcommand{\ttmaxiter}{\texttt{maxiter}}

\newcommand{\abs}[1]{\lvert#1\rvert}

\newcommand{\field}[1]{\mathbb{#1}}
\newcommand{\An}[2]{\ensuremath{\field{A}}\left(#1, #2\right)} 

\newcommand{\Arg}[1]{\text{Arg}(#1)}

\newcommand{\moveimg}[2]{%
  \makebox[0pt][l]{\hspace{#1}#2}%
}

\begin{document}


\title[Maximally Generalized McMullen Maps]{On dynamical and parameter spaces for the Maximally Generalized McMullen family of rational maps}

\author[S. Boyd]{Suzanne Boyd}
\address{Department of Mathematical Sciences\\
University of Wisconsin Milwaukee\\
PO Box 413\\
Milwaukee, WI 53201, 
USA}
\email{sboyd@uwm.edu, ORCID: 0000-0002-9480-4848}
\thanks{The authors completed the research and writing of this article during an NSF-sponsored Research Experience for Undergraduates program at the University of Wisconsin-Milwaukee during June and July of 2026. We thank the NSF for the financial support of award \# 2447752, and we thank UWM for their hospitality.}

\author[J. Cervantes]{Jillian Cervantes}
\address{Department of Mathematical Sciences\\
University of Wisconsin Milwaukee\\
PO Box 413\\
Milwaukee, WI 53201, 
USA}
\email{cervan32@uwm.edu, ORCID: 0009-0001-3953-3221}
\thanks{J. Cervantes was supported by the National Science Foundation Graduate
Research Fellowship Program under Grant No 2025370184.}

\author[R. Foster]{Rachel Foster}
\address{Spelman College \\
350 Spelman Lane \\
Atlanta, GA 30314,
USA}
\email{rgfoster555@gmail.com, ORCID: 0009-0009-1768-0333} 
\thanks{Any opinions,
findings, and conclusions or recommendations expressed in this material are those of the
authors and do not necessarily reflect the views of the National Science Foundation.}

\author[O. Hwang]{Olivia Hwang}
\address{Department of Mathematics\\ Smith College \\
10 Elm Street \\
Northampton, MA 01063,
USA}
\email{ohwang@smith.edu, ORCID: 0009-0007-9874-9739}
\thanks{No generative AI was used in the preparation of this manuscript. }

\author[S. Marcano]{Sebastian Marcano}
\address{Department of Computational and Applied Mathematics\\
Rice University \\
PO Box 1892 \\
Houston, TX 77251,
USA}
\email{sm376@rice.edu}
\thanks{The authors have no relevant financial or non-financial competing interests to report.}

\author[M. Roughan]{Molly Roughan}
\address{Department of Mathematics and Statistics, Boston University\\
665 Commonwealth Ave\\ 
Boston, MA, 02215,
USA}
\email{mroughan@bu.edu, ORCID: 0009-0002-1318-6169}



\date{\today}

\begin{abstract}
We explore the maximally generalized McMullen family of rational maps $F_{n,d,c,a}(z) = z^n + \frac{a}{z^d} + c,$ for $n,d\in \mathbb{Z},$ with $n,d\geq 3$, $a\in \mathbb{C}^*$ and $c\in \mathbb{C}$. In the dynamical plane, we establish annular regions containing the Julia sets, and show that the geometric limits of both the Julia and filled Julia sets as $n,d\to \infty$ are contained in the unit circle. 
In parameter space, for any fixed $a\in \mathbb{C}^*$, we examine the boundedness locus in the $c$-plane. When $n=d$ we find that  as $n\to \infty$, the boundedness locus tends to the union of two unit circles centered at $\pm 2\sqrt{a}$.  For arbitrary $n,d$, we identify a spine consisting of unit circles determined by the critical values, and prove that the boundedness locus tends to the spine 
 as $n,d\to \infty$. 
\end{abstract}

\maketitle


\section{Introduction} \label{sec:intro}

The study of how dynamical and parameter spaces evolve under growth of the degree of a family of rational maps has produced a number of interesting geometric limit phenomena. For the unicritical polynomial family
\[
P_{n,c}(z)=z^n+c,
\]
Boyd and Schulz \cite{BoydSchulz} showed that as $n\to\infty$, the Julia sets tend to the unit circle while the corresponding Mandelbrot sets tend to the closed unit disk, except in the $|c|=1$ case which was then described by Kaschner and colleagues \cite{Kaschner2020,Kaschner2015}. 
For the family
\[
z\mapsto z^n+c+\frac{a}{z^n},
\]
Boyd and Schulz \cite{BoydSchulz} established analogous geometric
limit results for both Julia sets and boundedness loci.
These results suggest that pushing the degree (the $n$ in the above examples) higher often reveals simple geometric structures within parameter spaces.

A second line of investigation concerns the geometry of boundedness
loci for generalized McMullen maps. In
\cite{BoydHoeppner1,BoydMitchell},
boundedness loci, baby Mandelbrot sets, and spine structures were
studied in several one-dimensional slices of parameter space.

The present paper combines these two perspectives, studying both geometric limits under degree growth and the
parameter space geometry determined by the critical values
for the family
\[
F_{n,d,c,a}(z)
=
z^n+\frac{a}{z^d}+c,
\]
where $n,d\ge 3$ are fixed integers, and $a,c$ are complex with $a\neq 0$. Following Stoertz and collaborators, we refer to this as the {\em maximally generalized McMullen} family, as it simultaneously generalizes the classical McMullen maps obtained when $c=0$ and the generalized McMullen family obtained when $n=d$.

The study of McMullen-type families, their Julia sets, and their parameter spaces was a recurring theme throughout the work of Robert Devaney and many of his collaborators, and the questions considered here arise naturally within that line of research.
In the classical McMullen family \[ z^n+\frac{a}{z^d}, \] and its various generalizations, Robert Devaney and coauthors studied 
the dynamics and parameter spaces of McMullen maps and generalized McMullen maps, including the geometry of Julia sets, the structure of parameter spaces, and the occurrence of Mandelbrot-like phenomena
\cite{DevaneyRussell,DevaneyHalos,DevaneySurvey2013,devgar}. 

For a fixed choice of parameters, we define the {\em filled Julia set} in the usual manner:
\[
K(F_{n,d,c,a})
=
\left\{
z :
\{F_{n,d,c,a}^{\,m}(z)\}_{m=1}^{\infty}
\text{ is bounded}
\right\},
\]
and define the {\em Julia set} by
\[
J(F_{n,d,c,a})=\partial K(F_{n,d,c,a}).
\]

\begin{figure}[!ht]
    \centering

    \includegraphics[width=0.43\linewidth]{\detokenize{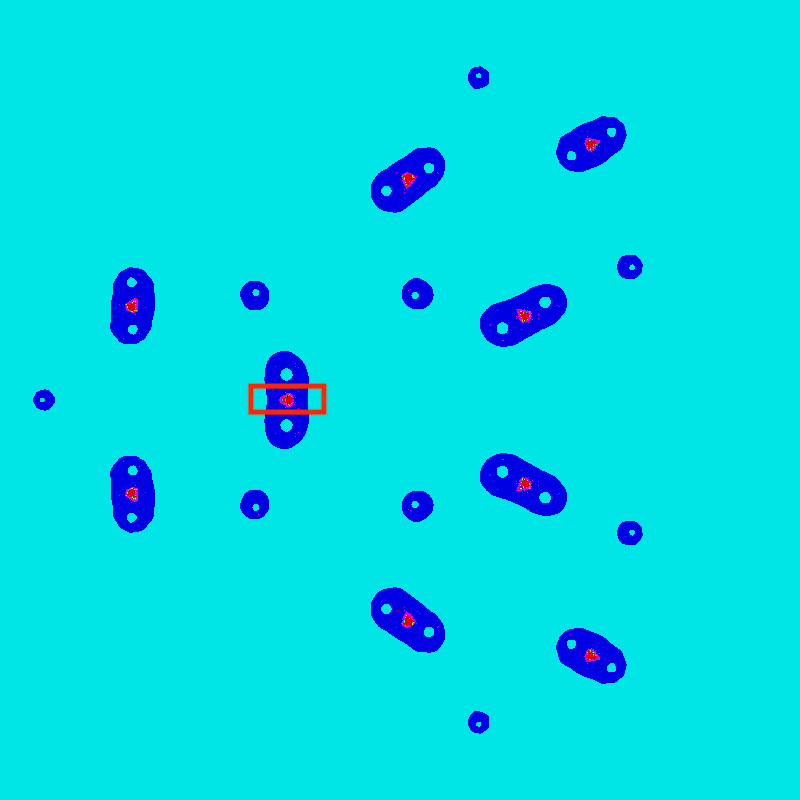}}
    \hspace{0.04\linewidth}
    \includegraphics[width=0.428\linewidth]{\detokenize{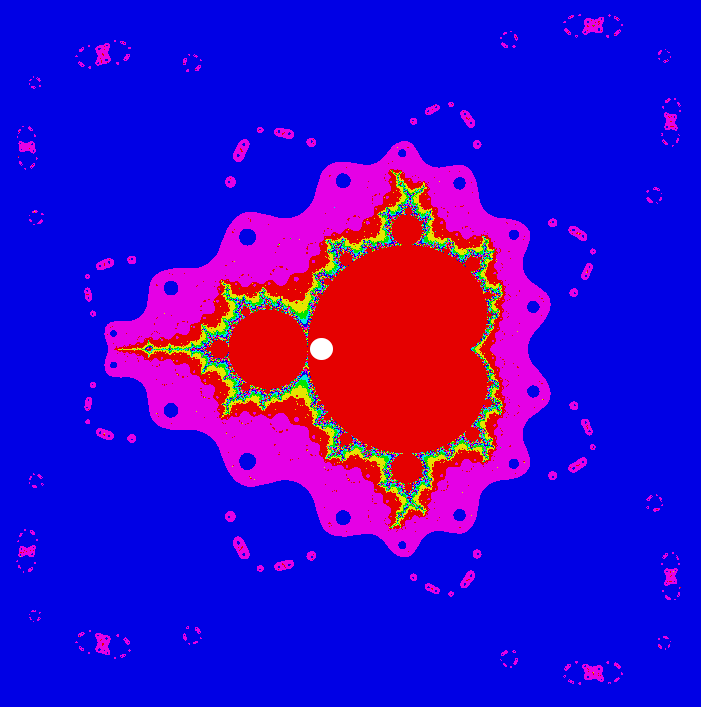}}
    \vspace{0.5em}
    
    \moveimg{4.7mm}{%
    \includegraphics[width=0.43\linewidth]{\detokenize{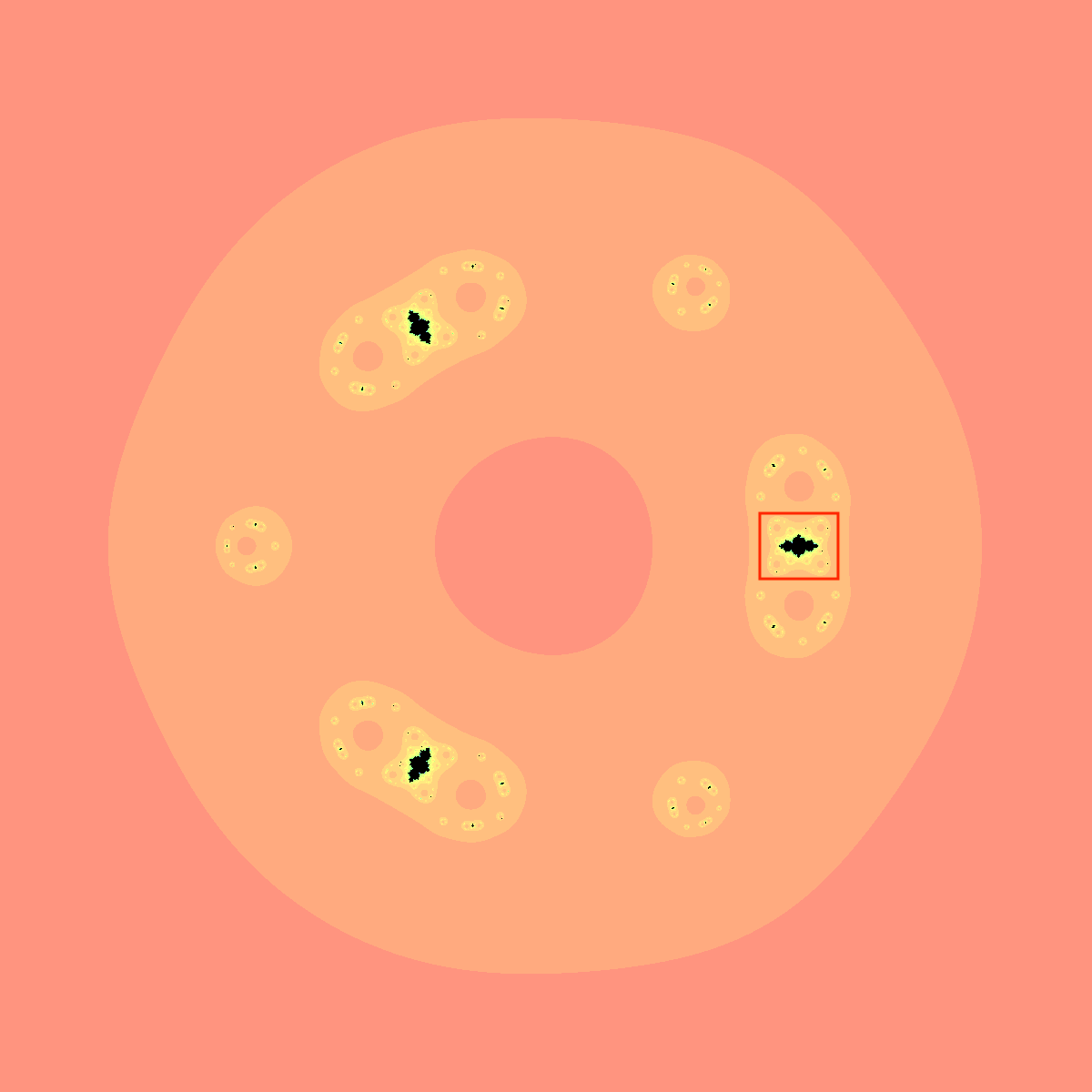}}
}
\hfill
\moveimg{-60mm}{%
    \includegraphics[width=0.429\linewidth]{\detokenize{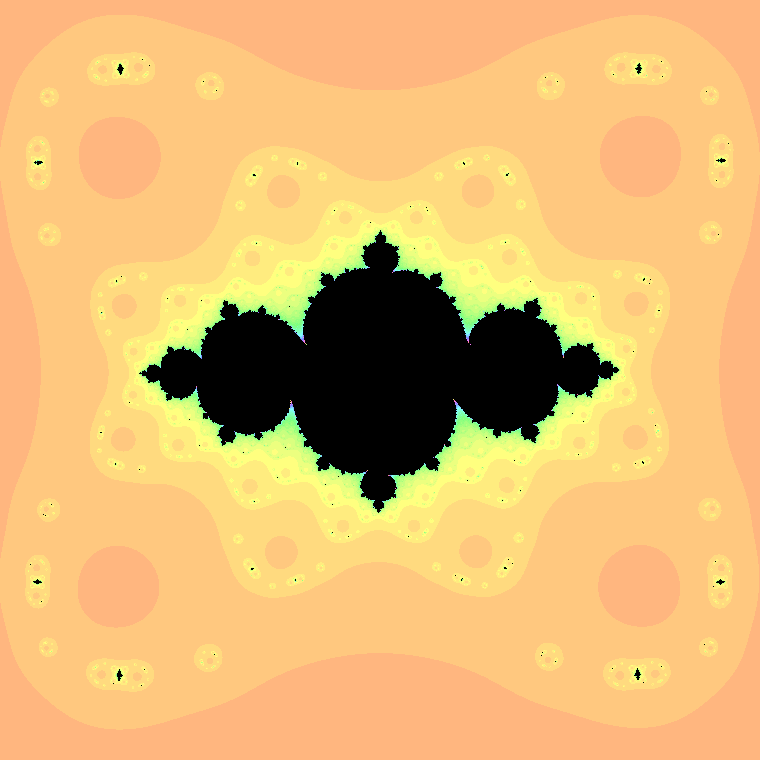}}
}  

    \caption{Let $n=6, d=3, a=1.$ Top-left:  Parameter window $[-3.5,3.5]\times [-3.5,3.5]$ in $c$-plane. The red-boxed baby Mandelbrot is magnified in the top-right ($c$-window: $[-1.1,-0.88]\times [-0.11,0.11]$). The white dot in the top right is at $c = (-1, -0.001)$; the bottom-left shows its corresponding Julia set ($z$-window: $[-2,2]\times [-2,2]$). There, the red box is magnified in the bottom-right (reveals an almost-``basilica'' baby Julia shape, in the $z$-window: $[0.809,1.043]\times [-0.115,115]$).}
    \label{fig:p_to_d}
\end{figure}

For fixed degrees and parameter $a\in\mathbb C^*$, we define the {\em boundedness locus} in the $c$-plane to be the set of parameters for which at least one critical orbit is bounded.
See Figure~\ref{fig:p_to_d} for both Julia set and boundedness locus examples. 

Our first objective is to understand the geometry of these sets as the degrees $n$ and $d$ become large. We identify annular neighborhoods of the unit circle which contain all Julia sets for all $n,d\geq 3$ (see Proposition~\ref{prop:K-in-annulus})
and then establish conditions under which the Julia sets converge, in the Hausdorff metric, into or onto $S^1$.

\begin{theorem}\label{thm:in+outJ1}
For any $c \in \CC$ and any $a \in \CC^*$, we have:
$$\lim_{n,d\to\infty} J(\Fndca) \subset S^1 \text{ and } \lim_{n,d\to\infty}  K(\Fndca) \subset S^1$$ if as $n,d\to\infty$ the limit(s) exist. 
Moreover, under the additional hypothesis that $ \ell = {\gcd(n,d)} \to \infty$ as $n,d\to \infty$, we have equality of the limits:
$$ \lim_{n,d \to \infty} \JFndca = \lim_{n,d \to \infty} \KFndca = S^1.$$
\end{theorem}

The key mechanism behind the equality statement is a rotational symmetry determined by $\ell = \gcd(n,d)$, which forces the critical points and critical values into highly structured configurations around the unit circle. As the quantity 
$\ell={\gcd(n,d)}$
grows, this symmetry forces increasingly many points of the Julia set to distribute around the unit circle.

Our second objective concerns parameter space. For fixed $a\in \CC^*$, we  study the boundedness locus in the $c$-plane:
\[
M_{n,d,a}
=
\left\{
c\in\CC :
\text{ at least one critical orbit of }
F_{n,d,c,a}
\text{ is bounded}
\right\}.
\]
 We begin with the special case $n=d$. In this setting we show that the boundedness locus is contained in two annular regions (Proposition~\ref{prop:GMparam-fat-annuli}) and converges to a pair of circles centered at $\pm2\sqrt a$. 
See Theorem~\ref{thm:GM-geometric-limit}.

The geometric limit  in the $a$-plane  for the $n=d$ case was studied in \cite{BoydSchulz}; here we instead investigate the geometry of the corresponding $c$-plane slices.  

We then pass to the general case $n\neq d$. 
The central objects in our parameter space analysis are the $\nu$ critical values.  We begin by establishing $\nu = \frac{n+d}{\gcd(n,d)}$ annular regions containing $\Mnda$   (see Proposition~\ref{prop:MGMparam-fat-annuli}). These determine a collection of circles in the parameter plane, which we call the \emph{spine} and denote by $\Cnda$ (see Lemma~\ref{lem:mgm-c-plane-spine}). The boundedness locus in Figure~\ref{fig:p_to_d} contains three spine circles.

We prove that, as the degrees grow, the boundedness locus is forced into arbitrarily small neighborhoods of this spine. Again we find as long as $\ell = \gcd(n,d)\to\infty$, the geometric limit of the boundedness locus coincides with the limiting spine.  Just as the unit circle governs the geometric limits in dynamical space, the spine governs the corresponding geometric limits in parameter space.

\begin{theorem} \label{thm:parameter_n_d_to_inf}
    For any $a\in \CC^*$, if $n,d\to \infty$ in such a way that $\ell=\gcd(n,d)\to \infty$,  then 
$$\lim_{n,d\to \infty} M_{n,d,a} = \lim_{n,d \to \infty} \calC_{n,d,a},$$
if and only if the limit on the right hand side exists.
    \end{theorem}

Throughout, we provide examples and figures to illustrate our results. Taken together, these results provide a description of the geometric limits arising in both dynamical space and parameter space for the maximally generalized McMullen family. On the dynamical side, Julia sets converge toward the unit circle. On the parameter side, boundedness loci are organized by explicit circle configurations determined by the critical values. These results extend previous work on generalized McMullen maps and demonstrate that many features of classical degree-growth phenomena persist in this substantially broader setting.

The paper is organized as follows. Section~\ref{sec:prelim} contains background and preliminary results. Section~\ref{sec:dynamicalplaneresults} studies Julia sets and proves Theorem~\ref{thm:in+outJ1}. Section~\ref{sec:gen-mcmullen-param-results} investigates parameter space in the special case $n=d$ and proves Theorem~\ref{thm:GM-geometric-limit}. Section~\ref{sec:MGM-param-results} introduces the spine construction and proves Theorem~\ref{thm:parameter_n_d_to_inf}. 
Additional figures and observations appear in Section~\ref{sec:pics}, and we conclude with directions for future research in Section~\ref{sec:future}.

\subsection*{Acknowledgments} We thank Brian Boyd for the program ``Dynamics Explorer'' \cite{detool} which we used to generate parameter and dynamical space images.  We also thank Kelsey Brouwer for the knowledge and inspiration she shared with us during her visits to our REU program. 

\subsection*{Dedication}
This article is dedicated to the memory of Robert L.\
Devaney. During the REU in which this work was completed, the student
authors began by learning the foundations of dynamics and complex dynamics from one of his introductory textbooks.
For the first author, this paper traces its origins to a conversation
more than fifteen years ago when Bob suggested studying maps of the form
\[
z^n+c+\frac{a}{z^n}.
\]
That suggestion led to many subsequent projects and publications,
culminating in the present work. We gratefully acknowledge Bob's profound influence as a mathematician, teacher, mentor, and ambassador for dynamical systems.

\section{Background and Preliminaries} \label{sec:prelim}

{\bf Notation.}  
First, some notation. We use $\CC^* = \CC\setminus\{0\}$, $\DD$ for the unit disk, $\Dbar$ for the closed unit disk, and $\DD_{r}$ for the disk of radius $r$ about the origin. 
For any $0 \leq r < R$, we denote the open annulus centered at the origin of inner radius $r$ and outer radius $R$ by:
$$ \An{r}{R} = \{ z \colon r < |z| < R \}.$$
If $z \in \CC$, we write $z + \An{r}{R}$ for the annulus centered at $z$ of inner radius $r$ and outer radius $R$.

We denote the $\eps$-neighborhood of a set $S$ by $\mathcal{N}_\eps(S)$, and the ball about a point $z$ of radius $\eps$ by $B_\eps(z)$.

\medskip

\subsection{Maximally Generalized McMullen maps}

We present here some basic facts about $F(z)$, just enough to get us started (also see \cite{Murali-PME}, and we go into more detail about this map in later sections when needed).  Set: 
\begin{equation} \label{eqn:basic-variables}
m:=n+d, \  \psi := \text{Arg}(a), \  \ell:=\gcd(n,d) \ \text{ and } \nu:=\frac{n+d}{\gcd(n,d)} = \frac{m}{\ell}.    
\end{equation}
Then $F$ has critical points at $0$ and at $\infty$, and there are 
$m$ additional critical points
given by:
\begin{equation*} \label{eqn:critpts}
    \left(\frac{da}{n}\right)^{1/m}.
\end{equation*}

Since these critical points are defined by an $m$th root, there are exactly
$m$ of them evenly spaced along a circle centered at the origin--we present the more detailed formula soon, but first, we calculate the critical values of $F$.

\begin{lemma} \label{lem:MGMcritvalues}
    The critical values of $F$ are:
\begin{equation*} \label{eqn:critvalues}
\frac{m}{d}
\left(\frac{da}{n}\right)^{n/m}
+ c.    
\end{equation*}
\end{lemma}

\begin{proof}
    Let $w$ be a critical point for the maximally generalized McMullen family given by $w=\left(\frac{da}{n}\right)^{1/m}$. We obtain the critical values by calculating
    $\Fndca(w)$. Equivalently, 
    $$
    \left(\frac{da}{n}\right)^{n/m}+\frac{a}{\left(\frac{da}{n}\right)^{d/m}}+c.
    $$ 
    We combine terms and simplify to 
    $$
    \frac{a(n+d)}{n^\frac{n}{m}(da)^\frac{d}{m}}+c.
    $$ Recall that $n+d=m$, which we use to further simplify to 
    $$
    \left(\frac{m}{d}\right)\frac{a^\frac{m-d}{m}d^\frac{m-d}{m}}{n^\frac{n}{m}}+c.
    $$ 
    Thus, we obtain
    $$
    \frac{m}{d}\left(\frac{da}{n}\right)^{n/m}+c.
    $$
    
\end{proof}

Next, we look more closely at the critical points and critical values.

\begin{lemma}\label{lem:MGM-k-maps}
Let $c\in\ \CC, a\in\CC^*$.     For $k \in \{0,\ldots,m-1\}$, we have the critical points for the maximally generalized McMullen family (excluding the critical points at $0$ and $\infty$) given by
$$
w_k := \left(\frac{d|a|}{n}\right)^{1/m}
\phi_k ^{1/m},
$$ 
where $\phi_k = e^{i(\psi+ 2\pi k}),$ with $\psi=\Arg{a}{}$, and 
$k=0,\ldots,m-1$. Let 
$$v_k := \frac{m}{d} \left(\frac{d|a|\phi_k }{n}\right)^{n/m} \qquad \text{(may also be denoted as } v).$$
Then 
$$
\Fndca(w_k) = v_k+ c,
$$
i.e., $v_k+c$ is the critical value corresponding to the critical point $w_k$.
\end{lemma}

\begin{proof}
First, for the critical points, observe that since $n,d >0$ are real numbers (in fact, integers) we have $\Arg{\frac{da}{n}} = \Arg{a}$, which we denote as $\psi$. Hence, 
$$\left(\frac{da}{n}\right)^{1/m}
= \left( \frac{d|a|}{n}\right)^{1/m} e^{i (\Arg{\frac{da}{n}} +2\pi k)\frac{1}{m}} 
=  \left( \frac{d|a|}{n}\right)^{1/m} e^{i (\psi +2\pi k)\frac{1}{m}}
= \left( \frac{d|a|}{n}\right)^{1/m} \phi_k^\frac{1}{m}.$$
    
    Now we find 
    $\Fndca(w_k)  $
    \begin{eqnarray*}
& = &  (w_k)^n + \frac{a}{w_k^d} + c 
=
\left(\left(\frac{d|a|}{n}\right)^{1/m}
\phi_k^{1/m}\right)^n + \frac{a}{\left({\left(\frac{d|a|}{n}\right)^{1/m}
\phi_k^{1/m}}\right)^d} + c\\
& = & 
\frac{d|a|\phi_k^{n/m}+|a|n\phi_k^{1-d/m}}{n^{n/m}(d|a|)^{{d/m}}}
= \frac{|a|^{n/m}\phi_k^{n/m}(d+n)}{n^{n/m}d^{{d/m}}}+c
= \frac{m}{d} \left(\frac{d|a|\phi_k }{n}\right)^{n/m} + c
\end{eqnarray*}
\end{proof}

Next, we describe how, due to symmetries of the map $F$, there are actually only
\[
\nu:= \frac{n+d}{\gcd(n,d)}= \frac{m}{\ell}
\]
 distinct critical values. 
Meyer, Wang, and Stoertz prove this symmetry in the dynamical plane as follows.
\begin{footnote}
{unpublished manuscript, 2024.}    
\end{footnote}

\begin{lemma} \label{lem:rotsym1}
For any point $z \in \mathbb{C}$,
\[
F\!\left(
z e^{\frac{2\pi i}{m}
\left(k+\frac{m}{\gcd(n,d)}j\right)}
\right)
=
F\!\left(
z e^{\frac{2\pi i}{m}k}
\right)
\]
for
\[
k = 0,1,\ldots,\frac{m}{\gcd(n,d)}-1,
\qquad
j = 0,1,\ldots,\gcd(n,d)-1.
\]
\end{lemma}

\begin{proof}
We have
\begin{align*}
F\!\left(
z e^{\frac{2\pi i}{m}
\left(k+\frac{m}{\gcd(n,d)}j\right)}
\right)
&=
F\!\left(
z e^{\frac{2\pi i}{m}k}
\left(e^{\frac{2\pi i}{\gcd(n,d)}}\right)^j
\right) \\
&=
z^n e^{\frac{2\pi i kn}{m}}
\left(e^{\frac{2\pi i n}{\gcd(n,d)}}\right)^j
+
a z^{-d} e^{-\frac{2\pi i kd}{m}}
\left(e^{\frac{2\pi i d}{\gcd(n,d)}}\right)^{-j}
+ c \\
&=
z^n e^{\frac{2\pi i kn}{m}}
+
a z^{-d} e^{-\frac{2\pi i kd}{m}}
+ c 
=
F\!\left(
z e^{\frac{2\pi i}{m}k}
\right),
\end{align*}
since $\gcd(n,d)$ divides both $n$ and $d$, implying
\[
\left(e^{\frac{2\pi i n}{\gcd(n,d)}}\right)^j
=
\left(e^{\frac{2\pi i d}{\gcd(n,d)}}\right)^j
= 1.
\]
\end{proof}

If we apply the above lemma to the case $k=0$, we obtain the following. 
\begin{corollary} \label{cor:rotsym}
    Rotation of $z$-values by integer multiples of 
\[
\frac{2\pi}{\gcd(n,d)}
\]
has no effect on their orbits under a map $F_{n,d,c,a}$. 
\end{corollary}

This corollary holds because under this rotation,  the original point and the rotated point map to the exact same
value. 

Moreover, Meyer, Stoertz, and Wang show that this is the only rotation which preserves the mapping.  

Given this, we can see why there are exactly $\nu = \frac{n+d}{\gcd(n,d)} = \frac{m}{\ell}$ critical values.

\begin{lemma} \label{lem:crit-val-count}
$\Fndca$ has exactly $\nu= \frac{m}{\ell}$ unique critical values. 
\end{lemma}

\begin{proof}
By Lemma~\ref{lem:rotsym1}, rotating a (critical or otherwise) point by an angle $ \frac{2\pi}{\ell}$ does not change its image.   Since the $m$ critical points are evenly spaced around a circle, they divide into equivalence classes of size $\ell = \gcd(n,d)$ of points having the same image (as there are $\ell$ rotations of angle $\frac{2\pi}{\ell}$ in a circle). So, each equivalence class maps to a different critical value,  so the number of distinct critical values is $ \frac{m}{\ell}.$ 
Note that increasing the index by $ \nu=\frac{m}{\ell} $ corresponds exactly to such a rotation. Hence $ F(w_{k+\nu})=F(w_k). $
\end{proof}

\subsection{Geometric Limits}

Here we state what we mean when we say the limit of a sequence of sets exists.

Recall the \textbf{Hausdorff metric} $d_\cH(A,B)$ for two sets $A$ and $B$ of a metric space $(X,d)$ is 
\begin{align*}
d_\cH(A,B) & =  \max\{\sup_{a\in A} \inf_{b\in B} d(a,b), \  \sup_{b\in B} \inf_{a\in A}d(a,b)\} \\
 & =  \max\{ \sup_{a\in A} d(a,B),  \ \sup_{b\in B} d(b,A). \}
\end{align*}

The Hausdorff distance between two sets $A, B$ can be thought of as the maximum of the maximum distance from $A$ to $B$ and the maximum distance from $B$ to $A$.    Consider the illustrative example $A=S^1$ and $B=\Dbar$.  We have $d_\cH(S^1, \Dbar) = 1$, since $0 \in \Dbar$ and $d(0, S^1) = 1$, even though $d(e^{it}, \Dbar) = 0$ for every $e^{it} \in S^1$.  Thus in each of our limit theorems, we need to supply two distance calculations.

\begin{definition}
\label{defn:set-limit}    
    If  $S_n$ is a sequence of compact subsets, then $S_n$ \textbf{converges} to compact subset $S$, and we write $\displaystyle \lim_{n \to \infty} S_n = S$,
if
for all $\eps > 0$ there is an $N\in \NN$ such that for all $n \geq N$ we have $d_\cH(S_n, S) < \eps$.
\end{definition}

A compact metric space (like $\hat{\CC}$) endowed with the Hausdorff metric is complete (i.e., Cauchy sequences converge). 

Mandelbrot and Julia sets for polynomials and rational functions are compact (\cite{Beardon}).

\section{Dynamical Plane Results}
\label{sec:dynamicalplaneresults}

In this section, we provide results on the location of and behavior of Julia sets for the maximally generalized McMullen maps. One tool is the symmetry results presented in Section~\ref{sec:prelim}.  We note that though the generalized McMullen is symmetric under the involution $h_a(z) = \frac{|a|^{1/n}}{z}$  (see Lemma 2.4 of \cite{BoydMitchell}), there is no such symmetry for the maximally generalized McMullen. Still, because $0$ is a pre-image of infinity, we get the same qualitative behavior that a neighborhood of $0$ maps outside an escape radius. In this section, we begin with this examination; that is, studying the location of the Julia sets from various perspectives.
We then use the symmetry to study the geometric limit of the Julia sets as $n,d$ both tend to $\infty$. 

\subsection{Location of Julia sets}
\label{subsec:dynamicalplaneresultspart1}

 First, for any $n,d$, we provide some bounds on the location of the Julia sets. (Compare the next result to Lemma 3.5 from \cite{BoydHoeppner1}.)

\begin{lemma} \label{lem:BoydHoeppner3.5}
Consider $c \in \CC$, $a \in \CC ^*$, and integers $n,d \ge 3$. Let $s=\max \{ 3,|c|,|a| \}.$ For any $|z| \ge s$, the orbit of z under $\Fndca$ escapes to infinity; 
i.e., $$\KFndca \subset \DD_s.$$
\end{lemma}

\begin{proof}
Note we have $s \ge 3$, hence for any positive integer $N\ge 3$, as $N-1 \ge 2$, we can conclude $s^{N-1} \ge 3^2 > 3$, hence $s^N > 3s$. 

We require $n,d \geq N$. Thus, by definition, $s^{n-2}\ge 3$ and $s^{d-2} \ge 3$. Now let $|z| \ge s$. We show by induction $|\Fndca^m (z)| \ge s^m$ for all $m \ge 1$. 

First, 
\[ 
|\Fndca (z)| = \left | z^n+c+\frac{a}{z^d} \right | 
\ge |z|^n-|c|-\frac{|a|}{|z|^d}
\ge 3s-s- \frac{|a|}{3s} \] \[
\ge 3s-s- \frac{|a|}{s}
=2s- \frac{|a|}{s}
= s+ \left ( s- \frac{|a|}{s} \right ) \ge s. 
\]

Now, suppose for some $m \ge 1$, $|\Fndca^m (z)| \ge s^m$. Then, 
\begin{align*}
&| \Fndca ^{m+1} (z)|  = \left | (\Fndca^m (z))^n + c + \frac{a}{(\Fndca^m (z))^d} \right | \\
& \ge |\Fndca^m (z)|^n - |c| - \frac{|a|}{|\Fndca^m (z)|^d}
\ge s^{mn}-s-\frac{|a|}{s^{md}}\\
& \ge s^{mn-m+1}-s-\frac{|a|}{s^{md}}
=s^{m(n-1)+1}-s-\frac{|a|}{s^{md}}\\
&=s^{m+1} \left ( s^{m(n-2)}-s^{-m}-\frac{|a|}{s^{md+m+1}} \right ) 
\ge s^{m+1}(3^m-1-1) \ge s^{m+1}.
\end{align*}

Thus by induction, $|\Fndca ^m (z)| \ge s^m$ for all $m \ge 1$. Because $s \ge 3 > 1$, the orbit of $z$ under $\Fndca$ escapes to infinity for any $|z|\geq s$.
\end{proof}

Examining the proof of the above Lemma~\ref{lem:BoydHoeppner3.5} we see it immediately yields a uniform bound for all $c,a$ in a bounded region. 

\begin{corollary}
    \label{cor:uniformbound1}
    If $Q,R>0$ are real numbers, then for all $c,a$ such that $|c| \leq R, 0 < |a| \leq Q$, and for any $n,d\geq 3$, we have 
$$K(\Fndca) \subset \DD_{\max\{ 3, R, Q \}}.$$
\end{corollary}

Next, we consider limiting cases, for large $n$ and/or $d$.
First, we observe with a straightforward adaptation that if we require only $n\to \infty$, no matter the $d\geq 3$, the Julia sets ``from the outside'' shrink down to the unit disk. (Compare to Lemma 5.1 of \cite{BoydSchulz} for the $n=d$ case.)

\begin{lemma} \label{lem:outsideJ1}
For any $c\in \CC$ and any $a\in \CC^*$, given any $\eps>0$, there is an $N\ge 3$ such that if $n \geq N$ and $d\ge 3$, then  $\KFndca \subset {\DD_{1+\eps}}$.
\end{lemma}

\begin{proof}
Analogous to the proof of Lemma 5.1 of \cite{BoydSchulz}, let $\sigma > \max\{1, |a|,|c|\}$. Choose $N\geq 3$ so that $(1+\eps)^N > 3\sigma$ and $\sigma^{N-2}>3$, and suppose $n\geq N$, $d\geq 3$. Consider any $z$ with $|z| \ge 1+\eps$. 

We claim that $|F_{n,d,c,a}^m(z)| \geq \sigma^m$ for all $m \geq 1$.

First, 
\begin{align*}
|F_{n,d,c,a}(z)| &= \left | z^n + c + \frac{a}{z^d}\right |    \geq |z|^n - |c| - \frac{|a|}{|z|^d} \ge |z|^n - |c| - \frac{|a|}{(1+\eps)^d} \\
& >(3\sigma-\sigma)- \frac{|a|}{(1+\eps)^d} =2\sigma-\frac{|a|}{(1+\eps)^d}=\sigma+\left ( \sigma-\frac{|a|}{(1+\eps)^d} \right ) \\
& > \sigma + \left (\sigma -{|a|} \right )
> \sigma.
\end{align*}

Now, suppose for some $m \geq 1$, $|F_{n,d,c,a}^m(z)| \geq \sigma^m$.  Then
\begin{align*}
& |F_{n,d,c,a}^{m+1}(z)|  = \left | (F_{n,d,c,a}^m(z))^n + c + \frac{a}{(F_{n,d,c,a}^m(z))^d} \right | \\
&  \geq   |F_{n,d,c,a}^m(z)|^n - |c| - \frac{|a|}{ |F_{n,d,c,a}^m(z)|^d}
 > \sigma^{mn} - \sigma - \frac{|a|}{\sigma^{md}}\\
&  > \sigma^{mn-m+1} - \sigma - \frac{|a|}{\sigma^{md}}
 = \sigma^{m+1} \left ( \sigma^{m(n-2)} - \frac{|a|}{\sigma^{md+m+1}} - \sigma^{-m} \right )\\
& > \sigma^{m+1} (3^{m}-1-1)  
 > \sigma^{m+1}.
\end{align*}

Thus by induction, $|F_{n,d,c,a}^m(z)| \geq \sigma^m$ for all $m \geq 1$.  Since $\sigma>1$, 
the orbit of $z$ under $F_{n,d,c,a}$ escapes to infinity, hence $z \notin K (F_{n,d,c,a})$.
\end{proof}

\begin{figure}[!ht] 
    \centering
    \includegraphics[width=0.3\textwidth]{\detokenize{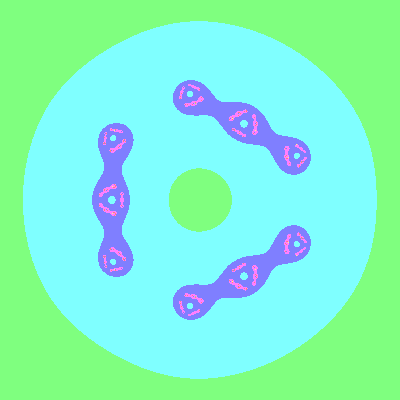}}\hfill
    \includegraphics[width=0.3\textwidth]{\detokenize{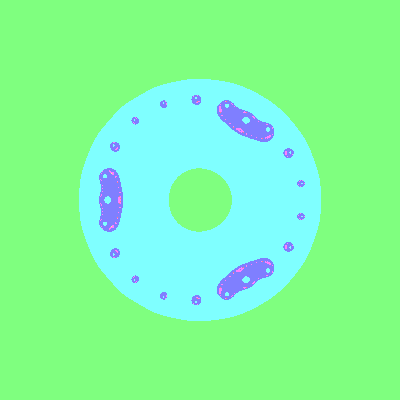}}\hfill
    \includegraphics[width=0.3\textwidth]{\detokenize{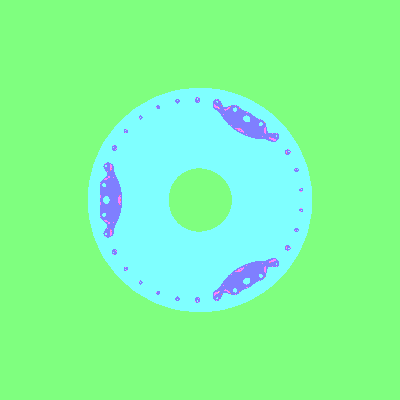}}\hfill
    \caption{$d=3, c=1, a=1$, and $n = 6, 18, 30$, and we draw the window $[-2,2]\times[-2,2]$. Observe the outer annulus bounding the Julia set approaching the unit circle as $n \xrightarrow{} \infty$. Also, $\ell= \gcd(n,d)=3$, observe the $3$-fold symmetry in the images.}
    \label{fig:n->inf}
\end{figure}

For an example of this phenomenon, see Figure~\ref{fig:n->inf}.

Note if we bound $|c|$ and $|a|$, we  obtain a common $N$ for the above Lemma~\ref{lem:outsideJ1}.

\begin{corollary}
    \label{cor:outside-J-common-N}
    If $Q,R>0$ are real numbers, then given any $\eps>0$ there is an $N\geq 3$ such that for all $n\geq N$ with $d\geq 3$, we have 
$$K(\Fndca) \subset \DD_{1+\eps}$$ for all $c,a$ such that $|c| \leq R, 0 < |a| \leq Q.$
\end{corollary} 

\begin{proof}
    Simply replace $\sigma$ in the proof of Lemma ~\ref{lem:outsideJ1} with $\sigma > \max\{1, Q, R\}$. 
\end{proof}

Next, we show that there is a neighborhood of $z=0$, inside which orbits also tend to infinity, creating an annulus where the Julia sets lie. 

\begin{lemma} \label{lem:insideescape}
    Let $c\in \CC$, $a\in \CC^*$, and $n,d\ge 3$. Set $s=\max \{ 3,|c|,|a| \}$, as in Lemma~\ref{lem:BoydHoeppner3.5}. If $p < \left(\frac{|a|}{s+1+|c|}\right)^{\frac{1}{3}}$, then $\DD_{p} \subset \CC \setminus K(F_{n,d,c,a})$.
\end{lemma}

\begin{proof}
    Recall from Lemma~\ref{lem:BoydHoeppner3.5} that $K(\Fndca)\subset \DD_{(\max\{3,|c|,|a|\})}$ for all $n,d\geq~3$. Suppose $|z| < p$, and note that $p<1$, so we have $|z|^n<1$. 
    We also have that $|z|^d< p^d\le \frac{|a|}{s+1+|c|}$. Then note that
        $$
        |\Fndca(z)|
        = \left|z^n + c + \frac{a}{z^d}\right| 
        \geq \frac{|a|}{|z|^d} - |z|^n - |c| \\
        \geq \frac{|a|}{\frac{|a|}{s+1+|c|}} - 1 - |c| 
        = s.
   $$

    Thus, by Lemma \ref{lem:BoydHoeppner3.5}, the orbit of $z$ escapes to infinity.
  \end{proof}

This immediately yields a uniform bound for bounded $|c|,|a|$, similar to corollaries~\ref{cor:uniformbound1} and~\ref{cor:outside-J-common-N}.

\begin{corollary}
    \label{cor:inside-J-RQ}
If $R>0$ and $0<Q_1<Q_2$, and $c,a\in \CC$ satisfy $|c|\le R$ and $Q_1\le |a|\le Q_2$, then for $s=\max\{3,R, Q_2\}$, and for 
\begin{equation*} \label{eqn:defn-p}
p := \left(\frac{Q_1}{1+R+\max\{3,R, Q_2\}}\right)^{\frac{1}{3}},     
\end{equation*}
and all $n,d\geq 3$, we have $$\DD_{p} \subset \CC \setminus K(F_{n,d,c,a}).$$
\end{corollary} 

\begin{proof}
By Lemma~\ref{lem:insideescape}, observe 
$$
p = \left(\frac{Q_1}{\max\{3,R, Q_2\}+1+R}\right)^{\frac{1}{3}}
\le 
\left(\frac{|a|}{\max\{3,|c|,|a|\}+1+|c|}\right)^\frac{1}{3}.$$
\end{proof}

Combining Lemmas~\ref{lem:BoydHoeppner3.5} and~\ref{lem:insideescape} 
immediately yields the following. 

\begin{proposition} \label{prop:K-in-annulus}
Let $R>0$ and $0< Q_1\leq Q_2$. For all $c\in\CC$ and $a\in\CC^*$ such that $|c|\le R,\, Q_1 \le |a| \le Q_2$, we have 
$$\KFndca\subset\An{p}{s},
$$
for $s=\max\{3,Q_2,R\}$ and 
$p = \left(\frac{Q_1}{1+R+\max\{3,R, Q_2\}}\right)^{\frac{1}{3}}.$   
\end{proposition}

\medskip

\begin{figure}[!ht]
    \centering
    \includegraphics[width=0.3\textwidth]{\detokenize{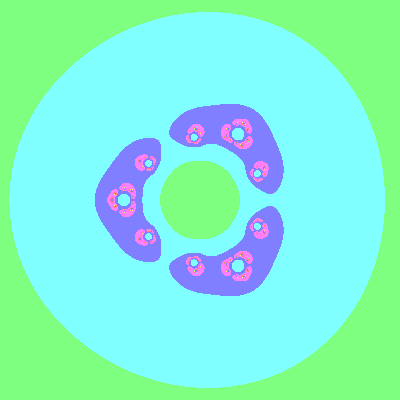}}\hfill
    \includegraphics[width=0.3\textwidth]{\detokenize{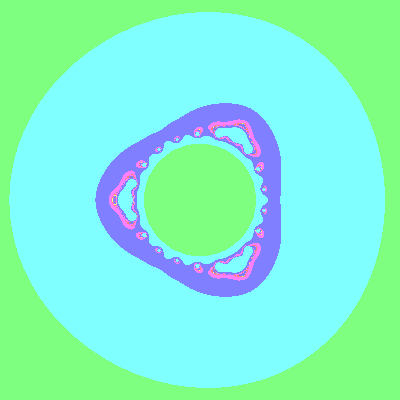}}\hfill
    \includegraphics[width=0.3\textwidth]{\detokenize{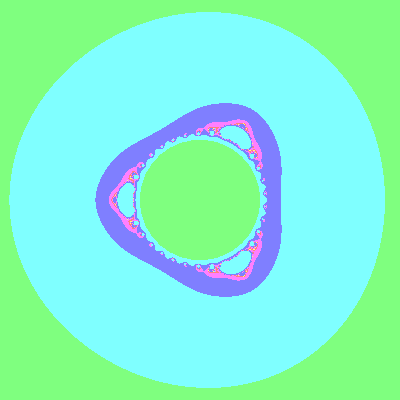}}\hfill
    \caption{\label{fig:dtoinfty}$n=3, c=1, a=1$, and $d = 6, 18, 30$, so $\ell = \gcd(n,d)=3$ (observe the $3$-fold symmetries), and we draw the window $[-3,3]\times[-3,3]$. Note the inner annulus bounding the Julia set approaches $S^1$ as $d \xrightarrow{} \infty$.}
    \label{fig:d->inf}
\end{figure}

\begin{example} \label{ex:J-annuli}
For example, for the $a,c$ values of Figures~\ref{fig:n->inf}, \ref{fig:d->inf}, \ref{fig:n-d-rel-prime}, and~\ref{fig:gcd-4-n-d-grow}, which are $c=1,a=1$, applying Proposition~\ref{prop:K-in-annulus} yields  $p = 0.58$ and $s = 3$; hence, we know these Julia sets all lie in the annulus $\An{0.58}{3}$. 
\end{example}

Next, we show that if $d\to\infty$, no matter the $n\geq 3$, the Julia sets push away from the origin as $d$ grows, to limit toward the Julia set escaping the open unit disk.

\begin{lemma} \label{lem:insideJ1}
    Consider $c\in \mathbb{C},$ $a\in \CC^*$, and $n \ge 3.$ Given $0 < 1 - \eps < 1$, there exists $ D\ge 3$ such that for all $ d \ge D$,  $\DD_{1 - \eps} \subset \mathbb{C} \setminus K(F_{n,c,d,a})$.

\end{lemma}

\begin{proof}
Let $z$ satisfy $|z| \le 1-\eps$.
   By the triangle inequality, 
        $$|\Fndca(z)|=\left | \frac{a}{z^d}+c+z^n \right |\ge\frac{|a|}{|z|^d}-|c|-|z|^n.$$ By assumption, $\frac{|a|}{|z|^d}\ge \frac{|a|}{(1-\eps)^d}$, and in order to guarantee escape using Lemma~\ref{lem:BoydHoeppner3.5}, we want
        $\frac{|a|}{(1-\eps)^d} \ge s+|c|+1>s+|c|+|z|^n$. Observe, if $D$ is so large that  $(1-\eps)^d \le \frac{|a|}{s+|c|+1}$ is satisfied for all $d\geq D$, then 
        $$
        \frac{|a|}{|z|^d}
        \ge \frac{|a|}{(1-\eps)^d} 
        \ge \frac{|a|}{\frac{|a|}{s+|c|+1}}
        =s+|c|+1,
        $$ 
        and so $\frac{|a|}{|z|^d} - |c| - 1 > s$.

        Thus, we have shown that $\Fndca(z)\geq s$ for all $|z|\le 1+\eps$ and the orbit of $z$ must escape to infinity by Lemma~\ref{lem:BoydHoeppner3.5}. Hence $z\notin K(\Fndca)$. 
\end{proof}

Figure~\ref{fig:d->inf} illustrates the behavior of the annulus in this case.

 Similarly to Corollary~\ref{cor:outside-J-common-N}, note that we can modify the proof of the above lemma to find a $D$ that works for all $|c|<R, Q_1<|a|<Q_2$, for any given $R>0$ and $0<Q_1\leq Q_2$. 

\begin{corollary}
\label{cor:inside-common-J}
Let $R>0$, $0<Q_1\leq Q_2$, and $0 < 1 - \eps < 1$. For all $c\in \CC$ and $a\in\CC^*$ such that $|c|\le R,\, Q_1 \le |a| \le Q_2$, we have 
    $$\DD_{1 - \eps} \subset \mathbb{C} \setminus K(F_{n,c,d,a}).$$
\end{corollary}

Combining all of the results in this section, we obtain that if {\em both} $n,d \to\infty$ then the filled Julia sets tend into the unit circle. Moreover, since $J = \partial K \subset K$, the Julia sets also get squeezed down into the unit circle. 
 
\begin{proposition} \label{prop:in+outJ1}
For any $c \in \CC$ and any $a \in \CC^*$, given any $\eps > 0$, there is an $N \geq 3$ such that for all $n,d \geq N$, we have $$K(F_{n,d,c,a}) \subset \mathbb{A}({1-\eps,1+\eps }).$$ 

Moreover, a common $N$ will suffice for all $c,a$ such that for some $R>0$ and $Q_1\le Q_2$, we have $|c|\le R, Q_1\le |a|\le Q_2$.
\end{proposition}

Hence, applying  Definition~\ref{defn:set-limit}, we immediately obtain the first statement of Theorem \ref{thm:in+outJ1}, that the Julia set and filled Julia set limits are contained in $S^1$, if they exist.

\subsection{Julia set limiting behavior}
\label{subsec:dynamicalplaneresultspart2}

In this subsection, we turn to considering when the limit of the prior theorem is equality. 

We first establish a further symmetry result for our family (compare this to Lemma 5.4 from \cite{BoydSchulz}).

\begin{lemma} \label{lem:MGMRcSinJ}
Recall $m=n+d$, $\ell=\gcd(n,d)$, and $\nu = \frac{m}{\ell}$. For any $c \in \CC$ and any $a \in \CC^*$,  if $\alpha$ is a primitive $\ell^{th}$ root of unity (i.e., $\alpha = e^{ i \frac{2\pi k}{\ell}},$ 
for any $k=0,\ldots, \ell-1$ which is co-prime to $\ell$; $k=1$ would suffice), then
 for any point $z \in \JFndca$ we have $\{ \alpha^j z \colon j \in \{ 0, 1, \ldots, \ell -1\} \}$ are 
 $\ell$
 distinct points in $\JFndca$, spread evenly around the circle of radius $\abs{z}$.
\end{lemma}

\begin{proof}
By Lemma~\ref{lem:rotsym1}, we have \[
F\!\left(
z e^{\frac{2\pi i}{m}
\left(k+\nu j\right)}
\right)
=
F\!\left(
z e^{\frac{2\pi i}{m}k}
\right)
\]
for
\[
k = 0,1,\ldots,\nu-1,
\qquad
j = 0,1,\ldots,\ell-1.
\]
Since $\frac{\nu}{m} = \frac{1}{\ell},$ given any $z\in \JFndca$, $\alpha^jz$ maps to $\Fndca(z)$. Since $\alpha$ is a primitive $\ell^{th}$ root of unity, $\alpha^\ell = 1$ but $\alpha^j \neq 1$ for any $j < \ell$, so since $0 \notin \JFndca$ (because $0$ maps to super attracting fixed point at $\infty$), but our $z\in \JFndca$, we know our $z \neq 0$, hence we have the set of all $\alpha^\ell z$, which has $\ell$ distinct points, and note for each $j$, $(\alpha^j)^\ell = 1.$  Since $\ell=\gcd(n,d)$, we also get $(\alpha^j)^n = (\alpha^j)^d =1.$
 Then, observe that:
 $$\Fndca( \alpha^j z) = (\alpha^j z)^n + c + \frac{a}{(\alpha^j z)^d} = (\alpha^j)^n z^n + c+ \frac{a}{(\alpha^j)^d z^d} = \Fndca(z).$$  
 Thus, each $\alpha^j z$ maps to the image of $z$.  Since $J$ is totally invariant (\cite{Beardon}), if $z$ is in $\JFndca$, then so is each $\alpha^j$, and of course each has norm $\abs{z}$.
\end{proof}

Now we provide criteria for when the Julia sets ``fill up'' the unit circle under degree growth.

\begin{lemma} \label{lem:MGMS1inJ}
Let $c \in \CC, a \in \CC^*$.  If $\{n_k, d_k\}$ are two sequences of integers such that $n_k,d_k\geq 3$ for all $k$, and $n_k, d_k, $and $\ell_k = {\gcd(n_k,d_k)} \to \infty$ as $k\to\infty$, then for each $\eps>0$, there exists $Q \geq 3$ s.t.\ for all $k \geq Q$, we have $B_\eps(e^{i\theta})\,\cap \, J(F_{n_k,d_k,c,a}) \neq \emptyset$.
\end{lemma}

(Comparing this to Lemma 5.5 of \cite{BoydSchulz}, we see the $n=d$ case is much simpler as $\ell=\gcd(n,n)=n$, so the symmetry is $n$-fold, hence the Julia sets always fill up $S^1$ as $n\to\infty$). 

\begin{proof}

Let $\eps >0$ be given.
By Proposition~\ref{prop:in+outJ1}, there is an $N \geq 3$ such that for all  $n,d\geq N$, $K(F_{n,d,c,a}) \subset \mathbb{A}({1-\frac{\eps}{2},1+\frac{\eps}{2} })$. So first, take $Q_1$ sufficiently large that $n_k,d_k > N$ for all $k\geq Q_1$.

Let $\gamma>0$ be the angle so that the polar rectangle defined by 
$U=\{re^{i\tau}: \  |r-1 |\le \frac{\eps}{2},  \ \ \ |\tau-\theta|\le \gamma\}$ 
is contained in $B_\eps(e^{i\theta})$, for any $e^{i\theta}\in S^1$. (Observe that $\gamma$ is independent of $\theta$; one could calculate $\gamma$ for the case of $\theta=0$ then apply for any $\theta$). 

We want to guarantee a point of $\JFndca$ in $U$. Notice first that for any point  $z_k\in\JFndca$ we have $|z_k-1|<\frac{\eps}{2}$ (so the radius is contained in the ball). Then, since $\ell_k \to \infty$, there is a $Q_2 \geq Q_1$ such that for all $k\geq Q_2$,  $\ell_k > \frac{\pi}{\gamma}$. Hence $2\gamma \ge\frac{2\pi}{\ell_k}$, which guarantees, due to the symmetry lemma, there is a point of $J(F_{n_k,d_k,c,a})$ in $U$ for all $k\ge Q_2$. 

We know every $J(F)$ is nonempty, since $F$ is a rational map (\cite{Beardon}).  Hence for any $k$, there exists a point $z_{k} \in J(F_{n_k,d_k,c,a})$.  
Lemma~\ref{lem:MGMRcSinJ} gives us that $z_{k} \in J(F_{n_k,d_k,c,a})$ implies the points $\alpha^j z_{k}$ are also in $J(F_{n_k,d_k,c,a})$, 
where $\alpha$ is a primitive $\ell_k$th root of unity, and these points are evenly distributed along the circle of radius $\abs{z_{k}}$ and are also in $J(F_{n_k,d_k,c,a})$. 
Thus given $\eps$, there is a $Q_3\geq Q_2$ such that for all $k\geq Q_3$, we know for some $j\in \{0,1, \ldots \ell-1\}$, there is a point of $\alpha^j z_k$ in $U\subset B_\eps(e^{i\theta})$. 

Thus, we have $B_\eps(e^{i\theta})\,\cap \, J(F_{n_k,d_k,c,a}) \neq \emptyset$, for all $k\geq Q_3$. Let the $Q$ of the lemma's conclusion be $Q_3$.

\end{proof}

Thus,  we find that, if we push both $n$ and $d$ towards $\infty$ in such a way that $\ell$ also tends to infinity, the Julia sets line up along and fill out the unit circle. (Compare to Theorem 1.4 of \cite{BoydSchulz}.)

 \begin{proof}[Proof of Theorem \ref{thm:in+outJ1}]
 The subset statement follows immediately from Proposition~\ref{prop:in+outJ1}.
 
For the equality, statement, first recall for any $n,d, c,a$, $\JFndca \subseteq \KFndca$.  

Proposition~\ref{prop:in+outJ1} implies that both $\displaystyle \inf_{z \in \JFndca} d(z, S^1)$ and $\displaystyle \inf_{z \in \KFndca} d(z, S^1)$ tend $0$ as $n,d\to\infty$.

Lemma~\ref{lem:MGMS1inJ} implies (under our hypotheses)
that for any $s \in S^1$,  $d(s, \JFndca)\to 0$ as $n,d \to \infty$, and since $\JFndca \subset \KFndca$, this implies $d(s, \KFndca) \to 0$ as $n,d \to \infty$.

Thus, $d_\cH (S^1, \JFndca)$ and $d_\cH(S^1, \KFndca)$ both tend to zero as $n,d, \ell\to \infty$, hence $\lim_{n,d\to\infty} \JFndca = \lim_{n,d\to\infty} \KFndca = S^1$ (by Definition~\ref{defn:set-limit}).
\end{proof}

In the following Example~\ref{ex:nucnd}, we have sequences with $\ell_k \to \infty$, hence they satisfy the additional theorem hypothesis.

\begin{example} \label{ex:nucnd}
    Let  $d_k=\{4,12,30,\ldots\}$ be a sequence tending to $\infty$ as $k\to \infty$, and set $n_k=2d_k=\{8,24,60,\ldots\}$. Notice that since $n_k$ is a multiple of $d_k$ that $\ell_k= d_k$ for any $k$. Then $\ell_k$ grows as $k$ grows. 

    We can generalize this to any case where $n_k$ and $d_k$ is a multiple of the other. Let $n_k$ be a sequence and $d_k=cn_k$ for some integer $c$ such that both $n_k$ and $d_k$ approach $\infty$. 
    Since $d_k=cn_k$, we have $\gcd(n_k,d_k)=n_k$. 
    Notice the same is true for the case $n_k=cd_k$.

Hence in these cases, the hypotheses of Lemma~\ref{lem:MGMS1inJ} are  satisfied. See Figure~\ref{fig:gcd-n-d-grow}.

\begin{figure}[!ht] 
    \centering
    \includegraphics[width=0.3\textwidth]{\detokenize{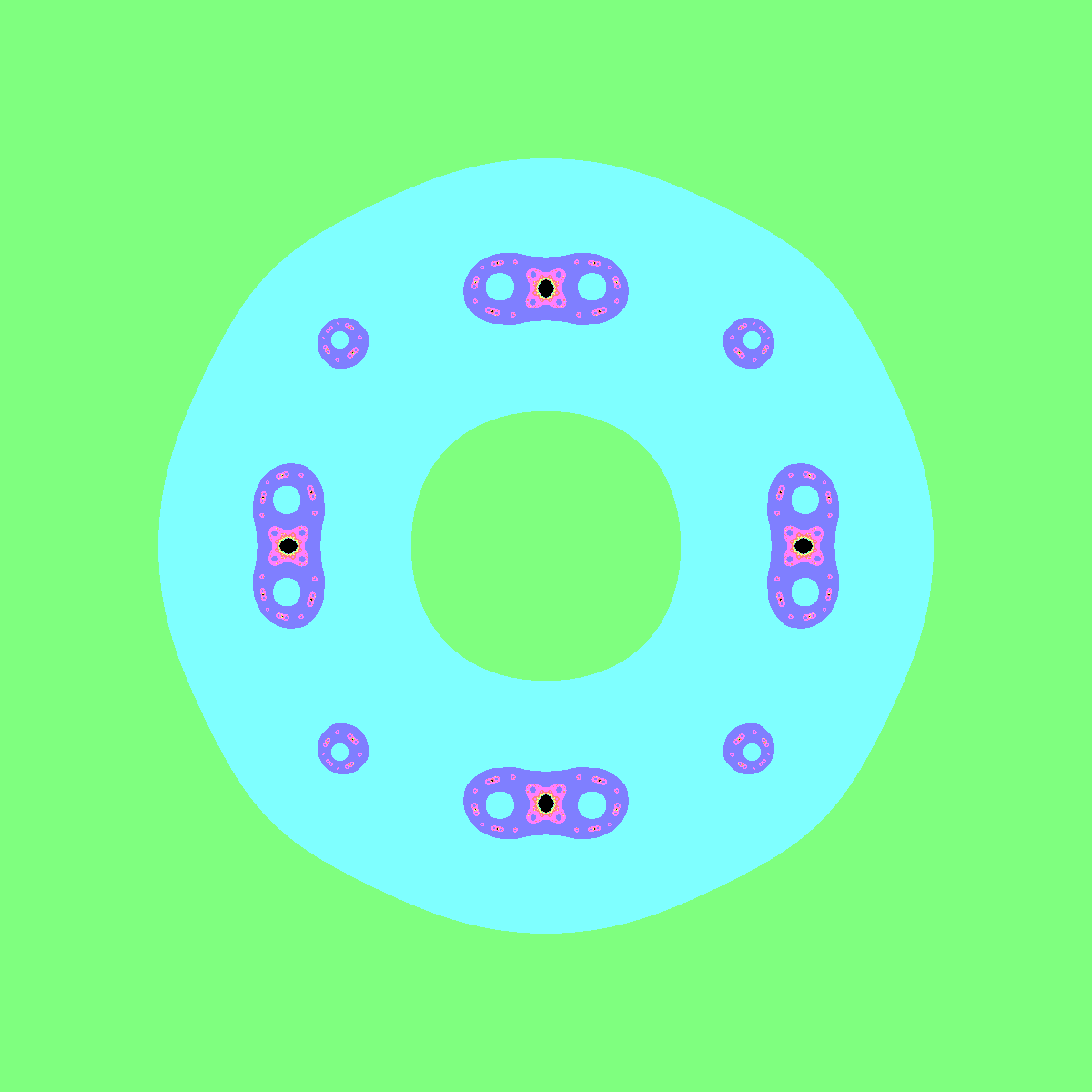}}\hfill
    \includegraphics[width=0.3\textwidth]{\detokenize{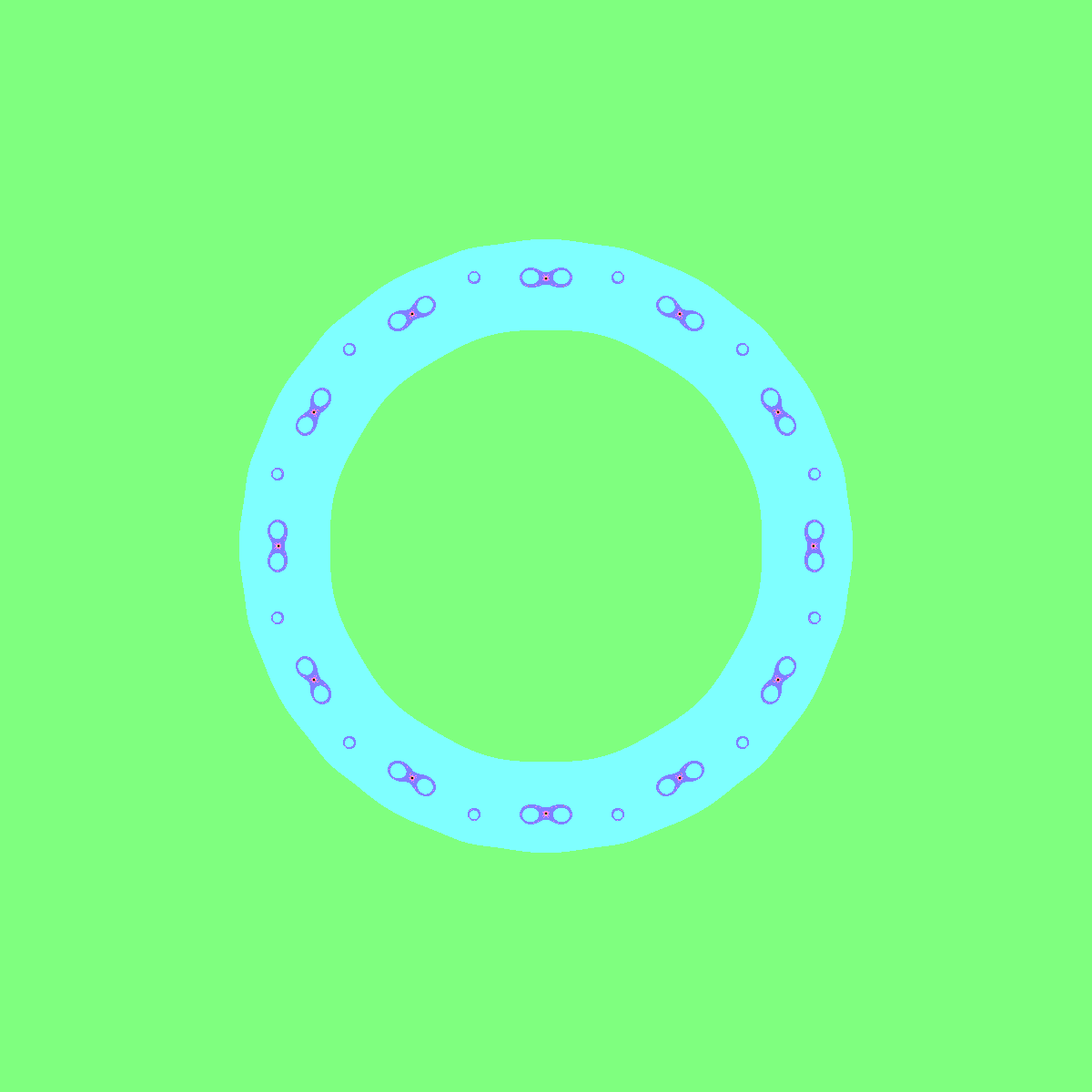}}\hfill
    \includegraphics[width=0.3\textwidth]{\detokenize{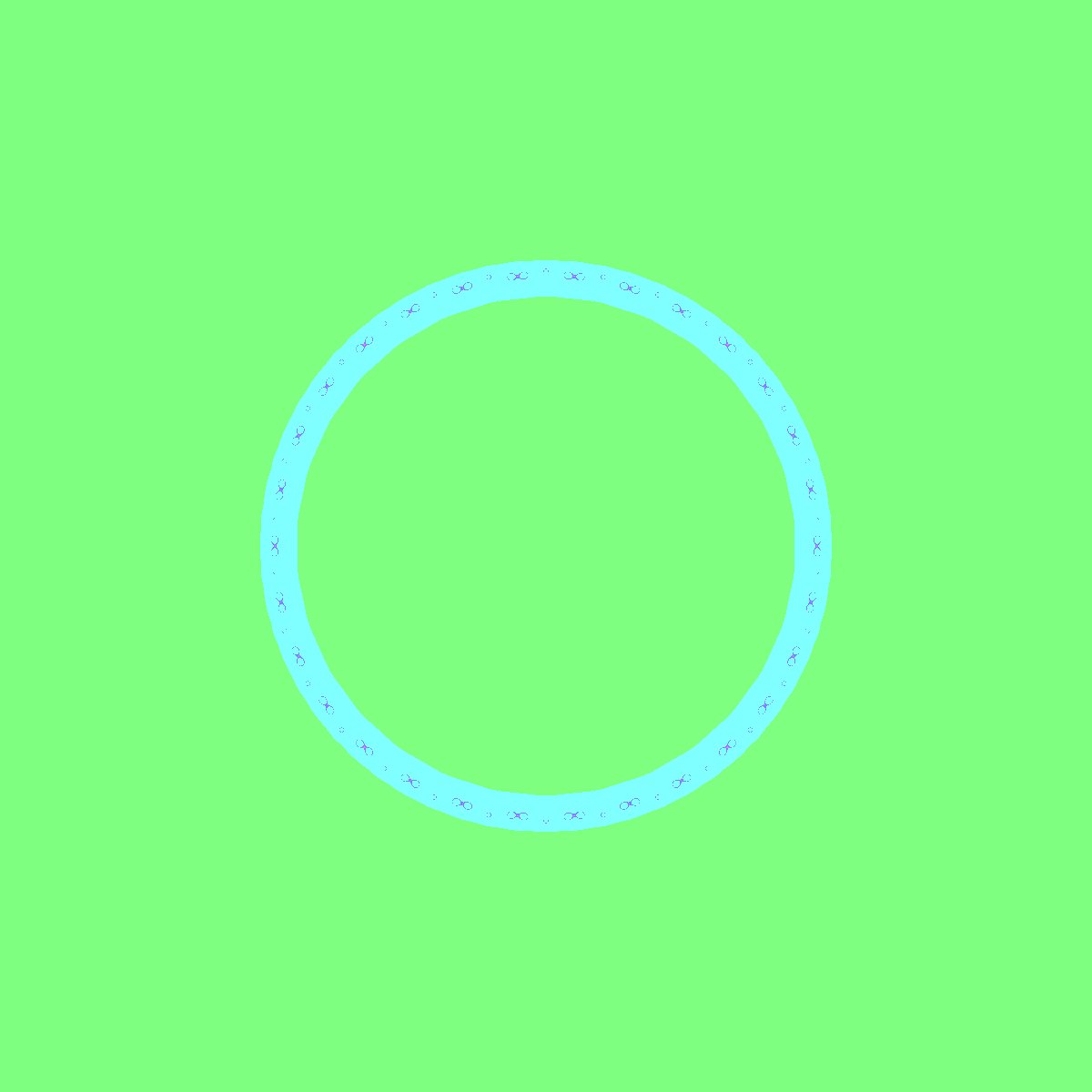}}\hfill
    \caption{$a=1,(n,d)=(8,4),(24,12),(60,30)$, and corresponding $c=(-0.95,0),(-0.91,0),(-0.8975,0)$ in the negative real axis (left to right). We draw the window $[-2,2]\times[-2,2]$. Here $\ell = \gcd(n,d)\to\infty$ as both $n$ and $d$ grow and the Julia sets fill the unit circle. 
    }
    \label{fig:gcd-n-d-grow}
\end{figure}

\end{example}

\begin{remark}
    Note if $n$ (or $d$) remains bounded, then $\gcd(n,d)$ cannot go to infinity. So, in order to meet the conditions of Proposition~\ref{prop:MGM-geom-lim-param},  we must have both $n$ and $d$ tending to infinity.
  \end{remark}

\section{Results on Parameter Space: the $n=d$ case}
\label{sec:gen-mcmullen-param-results}

In this section, we move to the parameter space for the simpler $n=d$ case to study the location of the boundedness locus in the $c$-plane. It was previously established in \cite{BoydSchulz} that the boundedness locus of $\Mnda$ in the $a$-plane is the polar lima\c{c}on $r=\frac{1}{2} (c + \cos\theta)$ shifted by $\frac{1}{4}(c^2-1)$. We use similar techniques to find the location of the boundedness locus in the $a$-plane, for the $c$-plane. 

In the $n=d$ case, there are simply two critical values: $c\pm 2\sqrt{a}$. We first discuss the approach for finding the boundedness locus.
\begin{remark}
    In general, consider the locus in $\CC^2$:
$$
\mathscr{L} = \{ (c,a) : \ c\pm 2\sqrt{a} = 1 \}.
$$ 
In the case of holding $c$ fixed, 
the intersection of this locus with any straight slice of $\CC^2$ obtained by fixing $c$ and allowing $a$ to vary is the polar lima\c{c}on $L_c$ examined in \cite{BoydSchulz}. 
On the other hand, we show below that holding $a$ fixed and examining a $c$-plane yields a curve defined by two circles. In general, one might be interested in a different 1-dimensional manifold in the parameter space (such as diagonal slices with $c=ta$ for various $t$ as in \cite{BoydHoeppner1}), or more dynamically significantly defined subfamilies).
In any case, if one is interested in the location of the boundedness locus, our approach to Theorem~\ref{thm:in+outJ1} yields the strategy employed here; that is, one must examine where the critical values, for the restrictions on $c,a$ of interest, are of modulus $1$. (In the next subsection, we also apply this strategy to examine where the modulus of the critical values of $\Fndca$ is $1$). 
\end{remark}

\medskip

Now, in this subsection, we use the simpler notation $\Fnca$ for the map 
$$
\Fnca(z) = z^n + c + \frac{a}{z^n}.
$$ 

\begin{notation}    
For a fixed $a\in \CC^*$, we let $\Mna$ denote the subset of the $c$-plane consisting of all $c$-values for which at least one critical orbit of $\Fnca$ is bounded.
\end{notation}

By Proposition~\ref{prop:in+outJ1} applied to the $n=d$ case (equivalently, Corollary 5.3 of \cite{BoydSchulz}), we know that as $n$ grows, the Julia sets lie in annuli shrinking down to $S^1$. 

We leverage the result here to show that the boundedness locus in the $c$-plane tends to the locus of $c$-values for which the critical values have modulus one, as follows.

\begin{proposition}
    \label{prop:c-plane-spine}
    Let $a\in \CC^*$. Let
    $$\Ca  = \Ca^+ \cup  \Ca^- :=  \{ e^{i\theta}\pm 2\sqrt{a} : \theta \in [0,2\pi)] \} = S^1 \pm 2\sqrt{a}.$$ 
    This is the locus of $c$-values for which the critical values of the Generalized McMullen map $\Fnca$ have modulus exactly $1$.

Moreover, for each $a\in \CC^*$, $\Ca$ is a union of two circles ($\Ca^\pm$), with centers at $\pm 2\sqrt{a}$ and each of radius $1$. The circles intersect as follows:
\begin{enumerate}
    \item The circles are disjoint if $|a| > 1/4;$
    \item The circles intersect in one point if $|a| = 1/4;$
    \item The circles intersect in two points if $|a| < 1/4$.
\end{enumerate}
Note the circles are only identical if $a=0$, which case we have excluded. 
\end{proposition}

\begin{proof}
    Recall the formula for critical points is $c\pm 2\sqrt{a}$, hence we want $|c\pm 2 \sqrt{a}|=1$, which means $c + 2 \sqrt{a}=e^{i\theta}$ and/or $c - 2 \sqrt{a}=e^{i\theta}$. We solve for $c$ and get 
    \[
    \{c=e^{i\theta}+2\sqrt{a}:\theta\in[0,2\pi)\}\cup\{c=e^{i\theta}-2\sqrt{a}:\theta\in[0,2\pi)\}. 
    \] 

    The description of $\Ca$ as circles is straightforward. The circles intersect precisely once when $2\sqrt{a}=1$, or $a=1/4$.
\end{proof}

\begin{figure}[!ht]
    \centering
    \includegraphics[width=0.3\textwidth]{\detokenize{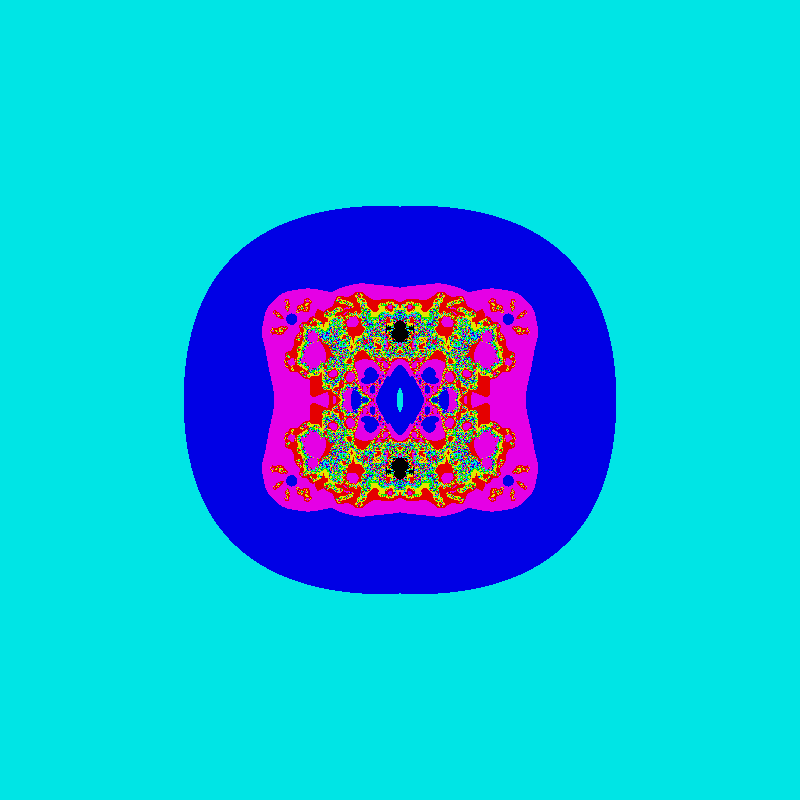}}\hfill
    \includegraphics[width=0.3\textwidth]{\detokenize{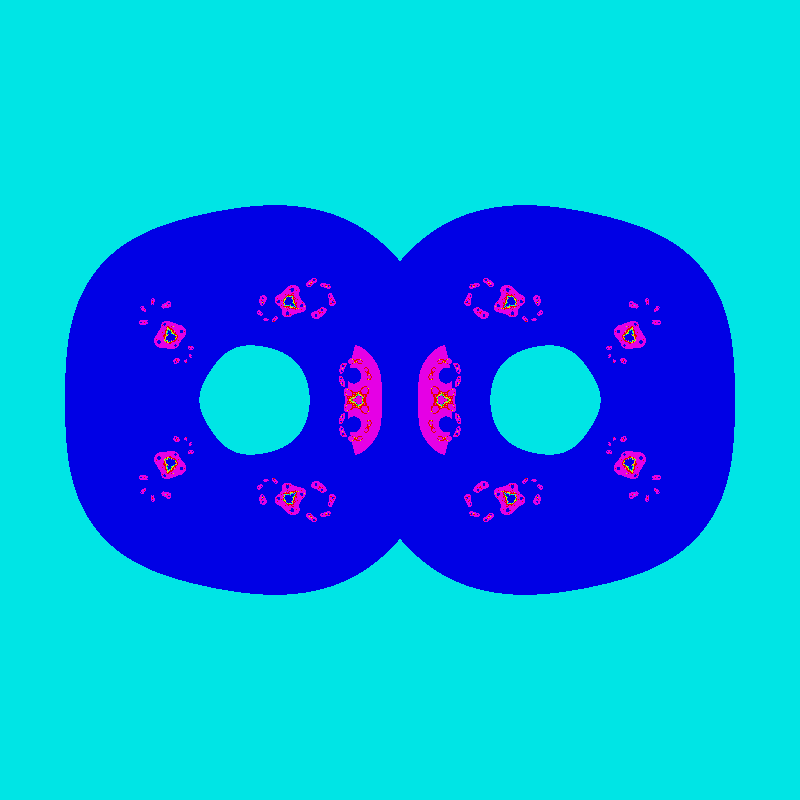}}\hfill
    \includegraphics[width=0.3\textwidth]{\detokenize{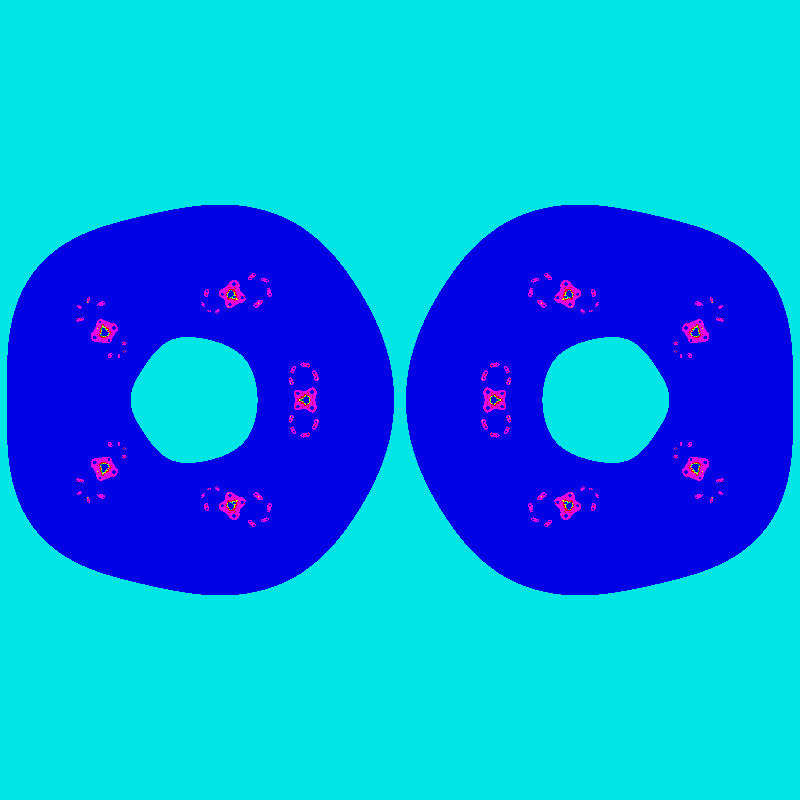}}\hfill
    \caption{The $c$-plane window $[-3.5,3.5]\times[-3.5,3.5]$ for $n=d=5$, and $a=0.01,0.4,0.8,$ left to right. The two clusters which seem to contain $\Mna$ separate slightly before $a=0.8$, though we know the theoretical spines separate at $a=0.25$.}
    \label{fig:two_clusters1}
\end{figure}

We call $\mathcal{C}_a$ the spine of the boundedness locus $\Mna$, as we will show that the boundedness locus is contained in annuli around these circles. The thickness of the annuli depends on $a$ and $n=d$. 
Note that when $n=d$, the two clusters separate at a threshold value of $|a|$ that depends on the common value of $n$ and $d$. For instance, Figure~\ref{fig:two_clusters1} illustrates the case $n=d=5$. For further insight into how the threshold value of $|a|$ required to separate the two clusters varies with $n=d$, refer to Section~\ref{sec:pics}.

Next, recall from Lemma~\ref{lem:BoydHoeppner3.5} that points $z$ with large modulus escape. Now, we establish something similar in parameter space.  

\begin{lemma} \label{lem:disk-for-Mna}
     Given $a\neq 0$, there exists a $T>0$ (which depends only on $a$) so that if $|c| \ge T$, then $c \notin \Mna$ for all $n \geq 3$; that is, both critical orbits escape. 
     Hence,  for all $n\geq 3$:
     $$\Mna \subset \DD_{T}.$$
\end{lemma}

\begin{proof}
In order to show that the orbits of the critical values escape, it is sufficient to show that the norm of their images is at least the escape radius $s = \max(3,|c|,|a|)$ established in Lemma~\ref{lem:BoydHoeppner3.5}. Recall that the critical values are $c\pm2\sqrt{a}$, so for this to be true, we will show:
\begin{align*}
   | \Fndca(c\pm2\sqrt{a})|&=\left |(c\pm2\sqrt{a})^n+\frac{a}{(c\pm2\sqrt{a})^n}+c \right |\\&\ge|c\pm2\sqrt{a}|^n- \frac{|a|}{|c\pm2\sqrt{a}|^n}-|c| \ge s.
\end{align*} 

To show this, we use $|c \pm 2\sqrt{a}|^n \ge (|c|-|2\sqrt{a}|)^n,$ and 
first we establish that $(|c|-|2\sqrt{a}|)^n>|a|$ in order to guarantee $-\frac{|a|}{|c\pm2\sqrt{a}|^n}>-1$. For simplicity, denote $|c|=X$ and $|2\sqrt{a}|=r_a$. \\

Using the binomial theorem, 
\begin{align*}
    (|c|-|2\sqrt{a}|)^n
    &=
    (X-r_a)^n
    =
    \sum_{k=0}^{n}((-1)^k{n\choose k}X^{n-k}r_a^k)\\
    &=
    X^n-nX^{n-1}r_a+\ldots+ (-1)^{n-1}X r_a^{n-1}  +(-1)^n r_a^n.
\end{align*}

Because $r_a$ is simply a constant, this is a positive polynomial based on the first $X^n$ term, meaning $(|c|-|2\sqrt{a}|)^n\to\infty$ as $|c|\to\infty$. Thus there must be a $T_1>0$ such that $(X-r_a)^n=(|c|-|2\sqrt{a}|)^n>|a|$ for some $X=|c| \ge T_1$ and all $n\geq 3$ (as larger $n$ corresponds to larger $X^n$, we can defined $T_1$ using $n=3$). Hence, for $|c| \ge T_1$, we have $-\frac{|a|}{|c\pm2\sqrt{a}|^n}>-1$, so:
\begin{align*}
    |\Fndca(v^\pm)|&
    \ge|c\pm2\sqrt{a}|^n-\frac{|a|}{|c\pm2\sqrt{a}|^n}-|c|
     >|c\pm2\sqrt{a}|^n-1-|c|
     \\ &
     \ge (|c|-|2\sqrt{a}|)^n-1-|c|
      =(X-r_a)^n-1-X,
    \end{align*}
so we want $(X-r_a)^n-1-X > s,$ for $X$ sufficiently large. 

Because $a$ is fixed, the worst case to show is when $s=|c|=X$. Again because $X^n$ is the leading term and is positive, 
we know that $(X-r_a)^n-1-X>X$ is true for sufficiently large $X$. Thus there exists a $T_2>0$ such that for all $X \ge T_2$, we have $(X-r_a)^n-1-X-\max(3,|a|,X)>0$ (for all $n\geq 3$ again). Then, we can define $T = \max(T_1,T_2)$, and then for all $|c| \ge T$, we know both critical orbits map outside the disk of radius $s$, which means both escape, hence $c \notin \Mna$. Thus, $\Mna \subset \DD_T$. 
\end{proof}

We use the above result to establish that all boundedness loci (for any $n,d\geq 3$ and any $a\in \CC^*$) lie in a union of annuli about the spine circles. 

\begin{proposition}
        \label{prop:GMparam-fat-annuli}
Let $a\in \CC^*$. Using Lemma~\ref{lem:disk-for-Mna},
let $T>0$ be such that $\Mna \subset \DD_T$ for all $n\geq 3$.
Let $s = \max\{3, |a|, T \}$ 
and $p = \left( \frac{|a|}{1+T+s} \right)^{1/3}$. Then,
for all $n \geq 3$, $\Mna$ lies in annular regions about the circles $\calC_a^{\pm}$; in particular:
\begin{eqnarray*}    
M_{n,a} & \subset & 
\An{p}{s} \pm 2\sqrt{a} 
\\
& = & \{ c \ : \ p < |c - 2\sqrt{a}|< s  \}
\cup \{ c \ : \ p < | c + 2\sqrt{a}| < s  \}.
\end{eqnarray*}
\end{proposition}

\begin{proof}
We shall prove that
$\DD_T \setminus (\An{p}{s} \pm 2\sqrt{a} )\subset \DD_T \setminus \Mna.$  To show this,  by Proposition~\ref{prop:K-in-annulus} (applied to $n=d$), because $a$ is fixed and $c\in \DD_T$ so $|c| \le T$, for all $n\geq 3$, we know $K(\Fnca) \subset \An{p}{s}$, for this $a$ and all $c\in \DD_T$.
 If $c\in \DD_T$, but not in either annulus, then 
the difference between $c$ to any annulus center is either $>s$ or $<p$, but that difference is exactly the formula for the critical values of the map for that $c$-value: $c\pm 2\sqrt{a}$. Equivalently, we say that $c\pm 2\sqrt{a}\notin\An{p}{s}.$ Hence, for such a $c$ value, the critical values are not in $\An{p}{s}$ and so 
neither can be in $K(\Fnca)$. Thus both critical values escape, so $c\notin \Mna$
and we have 
$$\DD_T \setminus (\An{p}{s} \pm 2\sqrt{a} )\subset \DD_T \setminus \Mna.$$

Now, suppose $c\in \Mna$; we want to show $c\in (\An{p}{s} \pm 2\sqrt{a} ). $ Since $c\in \Mna$ and $c\in \DD_T$, we know $c\notin \DD_t \setminus \Mna$. Hence, by the conclusion of the above paragraph, $c \notin \DD_T \setminus (\An{p}{s} \pm 2\sqrt{a} )$. Since $c\in \Mna \subset \DD_T, c\in \DD_T$. Thus it must be that $c\in (\An{p}{s} \pm 2\sqrt{a} )$, as desired.
\end{proof}

Similarly to the above, we use Lemma~\ref{lem:disk-for-Mna} next to yield
the following result on the $\eps$-neighborhood about $\Ca$.

\begin{proposition}
 \label{prop:GMparameps}
    Let $a\in \CC^*$, be fixed, and $0 < \eps < 1$. There is an $N\geq 3$ such that for all $n\geq N$, we have 
    \begin{eqnarray*}
    \Mna & \subset & 
    \calN_\eps(\Ca) 
 \   = \ \An{1-\eps}{1+\eps} \pm 2\sqrt{a}. 
    \end{eqnarray*}
    
    Thus, $\lim\limits_{n\to\infty} \Mna \subset \calC_a $  if this limit exists.
\end{proposition}
\begin{proof}

First, by Lemma~\ref{lem:disk-for-Mna} we know there is a $T>0$  so large that $\Mna \subset \DD_T$. 

Then (as in the proof of Proposition~\ref{prop:GMparam-fat-annuli}) it suffices to show that
$\DD_T\setminus \calN_\eps(\Ca) \subset \DD_T \setminus \Mna.$  Note that in $\DD_T$, $|c| < T$ and $a$ is fixed. So, by Proposition~\ref{prop:in+outJ1} (applied to $n=d$), 
there is an $N$ such that for all $n\geq N$, we have $K(\Fnca) \subset \An{1-\eps}{1+\eps}$ for this $a$ and all $c$ with $|c| < T$.
Let $n\geq N$. If $c\in \DD_T$, but not in either annulus, then (similar to the previous proof) 
the difference between $c$ to any annulus center is either $>1+\eps$ or $<1-\eps$, but that difference is exactly the formula for the critical values of the map for that $c$-value. Hence, for such a $c$ value, the critical values are not in $\An{1-\eps}{1+\eps}$ and so
neither can be in $K(\Fnca)$. Thus both critical values escape, so $c\notin \Mna.$

The final claim follows immediately, using  Definition~\ref{defn:set-limit}.

\end{proof}

We now show that the geometric limit of $\Mna$ fills up the spine (compare to Proposition 5.13 of \cite{BoydSchulz}). 

\begin{proposition}
    \label{prop:GM-geom-lim-param}
    Let $a\in \CC^*$,  $\eps>0$. If $\chi\in \Ca$, then there exists an $N$ such that for all $n\geq N$, \[
    B_\eps(\chi)\cap \Mna \neq \emptyset.
    \]
\end{proposition} 

\begin{proof}
    We want to find where $\Mna$ tends as $n\to \infty$, so we begin by finding some $c$ that are always contained in $\Mna$; the simplest such $c$ correspond to maps with at least one critical fixed point.  By Lemma 3.1 from \cite{BoydHoeppner1}, the formula for critical points is $a^{1/(2n)}=|a|^{1/(2n)}e^{i\left(\frac{\psi+2\pi k}{2n}\right)}$ for $k=0, \dots,2n-1$, where the $n$ critical points with odd $k$ map to $v^-=c-2\sqrt{a}$ and the $n$ critical points with even $k$ map to $v^+=c+2\sqrt{a}$. 
That is, $F(|a|^{1/2n} e^{i (\psi + 2\pi k)/(2n)} ) = c- 2\sqrt{a},$ if $k$ is odd, else $c+2\sqrt{a}$, if $k$ is even. 
   
    Now, we can explicitly write out the $c$-values for which the corresponding maps have these fixed critical points as
    \[
    C_{n,a}= C^+_{n,a} \cup C^-_{n,a}, 
    \] where  
    $$
    C^+_{n,a} = \{  |a|^{1/(2n)}e^{i(\psi+2\pi j)/2n}+2\sqrt{a} :j=0,\ldots,2n-1 \text{ and } j \text{ is odd}\}$$ 
    and 
    $$  C^-_{n,a}=\{  |a|^{1/(2n)}e^{i(\psi{a}+2\pi j)/2n}-2\sqrt{a}:j=0,\ldots,2n-1, \text{ and } j \text{ is even}\}.$$

Now $\chi \in \Ca =  \Ca^+ \cup \Ca^-$, where  $\Ca^\pm$ is the circle of radius 1 with center $\pm 2\sqrt{a}$. Without loss of generality, suppose $\chi \in \Ca^+$, so $\chi$ is a point a distance $1$ away from $2\sqrt{a}$. Thus $\chi=2\sqrt{a} + e^{i\theta}$ for some $\theta\in [0,2\pi)$. 

    Let $c^*$ be the point in $C^+_{n,a}$ closest to $\chi$, so $c^*=2\sqrt{a} + |a|^{1/(2n)} e^{i(\psi+2\pi j)/(2n)}$ for some (odd) $j$ in $\{0, \ldots, 2n-1\}$. Since $\chi=e^{i\theta}+2\sqrt{a} $, consider the triangle with vertices $2\sqrt{a}, \chi,$ and $c^*$. 
    The angle facing the edge from $\chi$ to $c^*$ is $\alpha_j = \min|\theta-(\psi+2\pi j)/2n|$. 
    Then by the law of cosines, the distance $d(\chi,c^*)$ between $\chi$ and $c^*$ is $$d(\chi, c^*)=\sqrt{|a|^{1/n}+1-2|a|^{1/(2n)}\cos(\alpha_j)}.$$ If we know $\eps>d(\chi, c^*)$, then $B_\eps(\chi)$ contains a point in $c$.

That is, we need that as $n\to \infty$, $d(\chi,c^*)$ tends to $0$. First, as $n\to \infty$, we know $|a|^{1/(2n)}\to 1$. 
Second, there are $n$ points in $C^+_{n,a}$ evenly distributed about the circle with center $2\sqrt{a}$ and radius $|a|^{1/(2n)}$. Thus as $n\to \infty$, we can make the angle $\alpha_j$ as small as we like, so its cosine tends to $1$. Hence the expression for $d(\chi,c^*)$ tends to $\sqrt{1+1-2} = 0$.

\end{proof}

Now combining Propositions~\ref{prop:GMparameps} and~\ref{prop:GM-geom-lim-param} (and Definition~\ref{defn:set-limit}), we immediately obtain: 

\begin{theorem} \label{thm:GM-geometric-limit}
For any $a\in \CC^*$, for the family $\Fnca(z) = z^n + c + \frac{a}{z^n},$ we have
$$
\displaystyle \lim_{n\to\infty} \Mna = S^1 \pm 2\sqrt{a}.
$$  
\end{theorem}

\medskip

    Finally, we observe  that $|a|$ determines the direction in which the ``centers'' in the boundedness loci, $c \in C_{n,a}$, approaches $\Ca$ as $n\to\infty$.

\begin{corollary}
   For all $n$, if $|a|<1$, then as $n$ increases,
$C_{n,a}^\pm$ approaches $\Ca^\pm$ from the inside of the circle $\Ca^\pm$, whereas if $|a|>1$, then the approach is from the outside. 
\end{corollary}

\begin{proof}
    Let $|a|<1$. Suppose $c^* \in C^+_{n,a}$. 
    From the equation defining $c^*$, we have the relationship $|c^*-2\sqrt{a}|=|a|^{\frac{1}{2n}}$. For any $n$, $|a|^{\frac{1}{2n}}<1$, so $|c^*-2\sqrt{a}|<1$. Thus the point $c^*$ is inside of the circle of radius $1$ about $2\sqrt{a}$, which is $\Ca^+$. 
    Now consider the case when $|a|>1$. Using the same equality, we get $|a|^{\frac{1}{2n}}>1$, and thus $|c^*-2\sqrt{a}|>1$.
\end{proof}

Next we proceed to the case in which $n$ may not be equal to $d$. All of the results in this section above can be derived from setting $n=d$ in the results in the following section, but we included them here for benefit of the reader's understanding. That is, we found the $n\neq d$ case more accessible after examining the $n=d$ case.

\section{Results on Parameter Space: allowing $n\neq d$}\label{sec:MGM-param-results}

In this section, we expand on the results from Section~$\ref{sec:gen-mcmullen-param-results}$ to the $n \ne d$ case, and examine the boundedness locus of the maximally generalized McMullen map 
$$\Fndca(z) = z^n  + \frac{a}{z^d} + c.$$  

In the $n = d$ case, there are two critical values $c\pm2\sqrt{a}$ which corresponds with two circles that make up the spine of the boundedness loci. In the $n\neq d$ case, our formula for critical points and values is dependent on $n, d, a$, thus we have a number of critical values which depends on $n,d$. Using symmetry, we established a formula for these critical points and values in Section~\ref{sec:prelim}; 
Lemma~\ref{lem:MGM-k-maps} gives the critical point formula $w_k = \left(\frac{d|a|}{n}\right)^{1/m}
(\phi_k )^{\frac{1}{m}},$ and its corresponding critical value formula $\Fndca(w_k) = v_k +c =\frac{m}{d} \left(\frac{d|a|\phi_k }{n}\right)^{n/m} + c$, where $\phi_k =e^{i(\psi+2\pi k)}$, with $\psi =\Arg{a}.$ 

\begin{notation} \label{not:Vand}
    For fixed $a\in \CC^*$, let $\Mnda$ denote the subset of the $c$-plane consisting of all $c$-values for which at least one critical orbit of $\Fndca$ is bounded.

Let $\Vand$ denote the  roots in the set $\frac{m}{d}\left(\frac{d a}{n}\right)^\frac{n}{m}$ (recall due to the symmetry of the map, there are only $\nu = \frac{n+d}{\gcd(n,d)} = \frac{m}{\ell}$  distinct values rather than $m$). Note that $\Vand$ is the set of all $v_k.$
That is, 
\begin{equation*} \label{eqn:Vnda}
\Vand = \left\{ v_k 
: k=0,\ldots,m-1
\right\}.
\end{equation*}

Finally, we denote the modulus of the $v_k$'s by:
    $$
r_{n,d,a} 
:= |v_k| = \frac{m}{d}|w_k|^n = \frac{m}{d}\left(\frac{d|a|}{n}\right)^\frac{n}{m}.
$$ 
  \end{notation}

Before presenting the details, we provide an overview of the main strategy of this section. We have $m=n+d$ critical points $(w_k$) arranged in a circle about $0$. Due to the symmetry of the map, they can be grouped into $\frac{m}{\ell}$ equivalence classes of $\ell=\gcd(n,d)$ points having the same image, that is, which map to the same critical value. The critical values are $F(w_k) = v_k +c$ (where $v_k$ is given by a formular in terms of $a,n,d,$ and $e^{i2\pi k})$. Essentially because $|z|^n$ or $|z|^{-d}$ is large or tiny for $|z|$ far from $1$, we get that the filled Julia sets tend to $S^1$ as $n,d$ grow (that is the key to our Section~\ref{sec:dynamicalplaneresults} results). Now, $\Mnda$ is the locus for which a critical orbit is bounded, meaning a critical orbit is in its filled Julia set, hence it a critical value $c+v_k$ must have norm close to $1$ for large $n,d$. So we ask for what $c$ is  $|c+v_k|=1$? Precisely on the circle of radius 1 about the center point $c=-v_k$, which is where the norm of $c+v_k$ would be $0$ if $c=-v_k$. Thus, the ``spine'' is the set of $c$-values which are $1$ away from the $-v_k$'s, and the set of all $v_k$'s is a set of $\nu = \frac{m}{\ell}$ unique values (because the set of $c+v_k$'s is $\nu$ values that's the number of critical values).  Finally, some $c$-values which we know are in $\Mnda$ are the ones with fixed critical point; that means for each of the $m$ critical points $w_k$, we get to solve the equation $w_k = v_k+c$; this produces $m$ $c$-values, our $c_k$'s, and since those group into sets of $\ell$ points which each map to the same critical value (by the symmetry of the map), we end up with $\ell$ points $c_k$ spread around each circle of a norm close to $1$ for $n,d$ large about the center point $-v_k$. Thus, $\ell\to\infty$ guarantees the convergence to the spine. 

Now, our first results of this section establish the location of spines, an analog of Proposition~\ref{prop:c-plane-spine} to the maximally generalized McMullen case.

\begin{lemma}\label{lem:mgm-c-plane-spine}
    Let $a\in \CC^*$ and $n,d\geq 3$. Then
\begin{eqnarray*}
    \Cnda & := & \cup_{k=0}^{m-1} \  \Cnda^k     
:= S^1 - \Vand 
= \bigcup_{k=0}^{m-1} \ \left\{ e^{i\theta}-
    v_k  :
    \theta \in [0,2\pi)] \right\}
\end{eqnarray*}    
    is the locus of $c$-values for which the critical values of the maximally generalized Mcmullen map $\Fndca$ have modulus exactly $1$.
    
Moreover, this set consists of $\nu = \frac{m}{\ell} = \frac{n+d}{\gcd(n,d)}$ circles of radius $1$, with centers the negatives of the points $v_k$ in $\Vand$, which 
are spread evenly around the circle (centered at the origin) of radius $r_{n,d,a},$ which satisfies
$$
|a|^\frac{n}{m} \leq r_{n,d,a} = |v_k|  \leq  2|a|^\frac{n}{m}. $$ \end{lemma}

\begin{proof}
    Recall the formula for critical values $(\frac md)(\frac{da}{n})^{\frac nm}+c = \Vand+c,$ where $m=n+d$. Thus we want $|(\frac md)(\frac{da}n)^{\frac nm}+c|=1$ or equivalently $(\frac md)(\frac{da}{n})^{\frac nm}+c=e^{i\theta}$, and then we solve for $c$, which results in the roots as claimed above. 
    
    As we observed in Lemma~\ref{lem:crit-val-count}, after taking the $m^{th}$ root and then raising to the $n^{th}$ power, some of the values may repeat, so we end up with $\nu = \frac m{\gcd(n,d)}$ unique critical values. So this set consists of $\nu$ circles of radius $1$ with centers all of the values $v_k$ in $\Vand$
which means the centers (the points $v_k$ in the set $\mathcal{V}_{a,n,d}$) are spread evenly around the circle (centered at the origin) of radius $r_{n,d,a}$.

   Finally, Meyer, Stoertz and Wang also show that we have $1\leq \frac{m}{d} \left( \frac{d}{n} \right)^{\frac{n}{m}}\leq 2$ for all $n,d\geq 3$, $\frac{m}{d} \left( \frac{d}{n} \right)^{\frac{n}{m}}$ reaches its global maximum of $2$ along the line $n=d$, and for any $n$ fixed, $\frac{m}{d} \left( \frac{d}{n} \right)^{\frac{n}{m}}\to 1$ as $d\to\infty$. 
    
Since $r_{n,d,a}  = \frac{m}{d} \left( \frac{d}{n} \right)^{\frac{n}{m}}|a|^\frac{n}{m}$, this yields the desired bound, in terms of $n,d,$, and $|a|$, on the norm of the $v_k$'s.

    \end{proof}

Now we provide a description of for which $a$-values the circles in the spines are separated.

\begin{proposition} \label{prop:cluster-separation}
     For any $n,d\geq 3$, where $m=n+d$, let 
    $$\rho({n,d}): =
    \dfrac{ 
    \left(\frac{d}{m}\right)^{\frac{m}{n}}\left(\frac{n}{d}\right)}
    {\left(\sin\left(\frac{\pi}{\nu}\right)^{\frac{m}{n}} \right)}.
    $$
    
 For each $a\in \CC^*$, the circles in the spine $\calC_{a,n,d}$ intersect as follows:
\begin{enumerate}
    \item if 
    $|a| > \rho({n,d})$, then the circles are disjoint; 
    \item if $|a| = \rho({n,d})$, then the circles are mutually tangent (each two adjacent circles intersect in one point); or
    \item if $|a| < \rho({n,d})$, then the circles intersect in more points (as $|a|\to 0$ the circles tend to coincide). 
\end{enumerate}
\end{proposition}

\begin{proof}
    We have $\nu$ unique critical values, each of which corresponds to a circle in $\calC_{a,n,d}$ (in the $c$-parameter plane), and the centers of the circles (the points $v_k$ in $\Vand$) are spaced evenly around a circle centered at the origin of radius 
$
r_{n,d,a}=\frac{m}{d}\left(\frac{d|a|}{n}\right)^\frac{n}{m}.
$ 
     We want to find the values of $a$ when the circles are mutually tangent. At this $a$-value, the centers of the $\nu$ circles form a regular $\nu$-gon. By basic trigonometry, we find that $\sin(\pi/\nu)=1/(r_{n,d,a})$ and solving for $|a|$ we obtain  
     $$
     |a| = \left(\frac{d}{m}\right)^{\frac{m}{n}}\left(\frac{n}{d}\right)/\left(\sin\left(\frac{\pi}{\nu}\right)^{\frac{m}{n}} \right).$$ The other two cases follow.
\end{proof}

\begin{figure} [!t]
    \centering
    \includegraphics[width=0.5\linewidth]{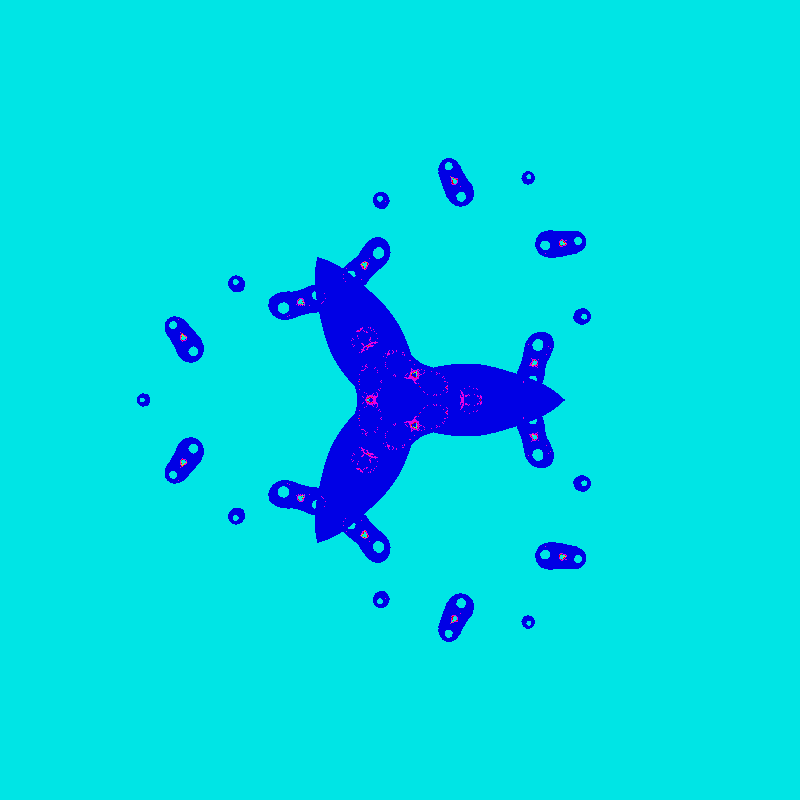}
    \caption{The $c$-plane window $[-3.5,3.5]\times[-3.5,3.5]$ for $d=5$, $n=10$, $a=0.48$.  }
    \label{fig:lemma5.3}
\end{figure}

Note, we can observe this in Figure~\ref{fig:lemma5.3} for the case in which $d=5$ and $n=10$, where the three clusters are barely separated when $a=0.48$, which agrees with the result above. We calculate the theoretical spine separation for more of our figures in the following example. For further insight into how the value of $|a|$ required to separate clusters differs across cases of lower $n,d$, refer to Section~\ref{sec:pics}.

\begin{figure}[!ht]
    \centering
    \includegraphics[width=0.49\textwidth]{\detokenize{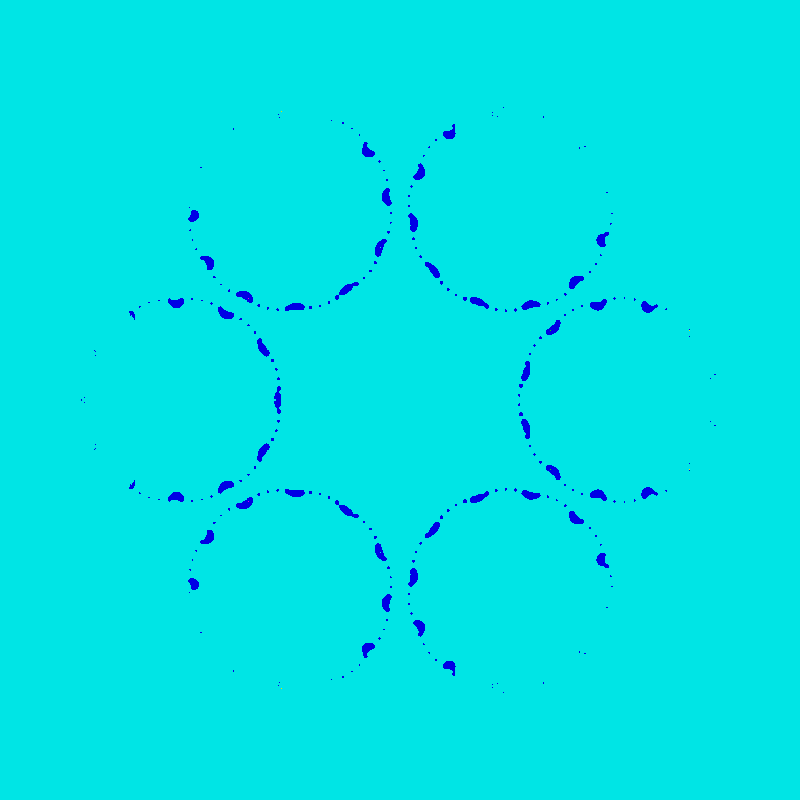}}
    \includegraphics[width=0.49\textwidth]{\detokenize{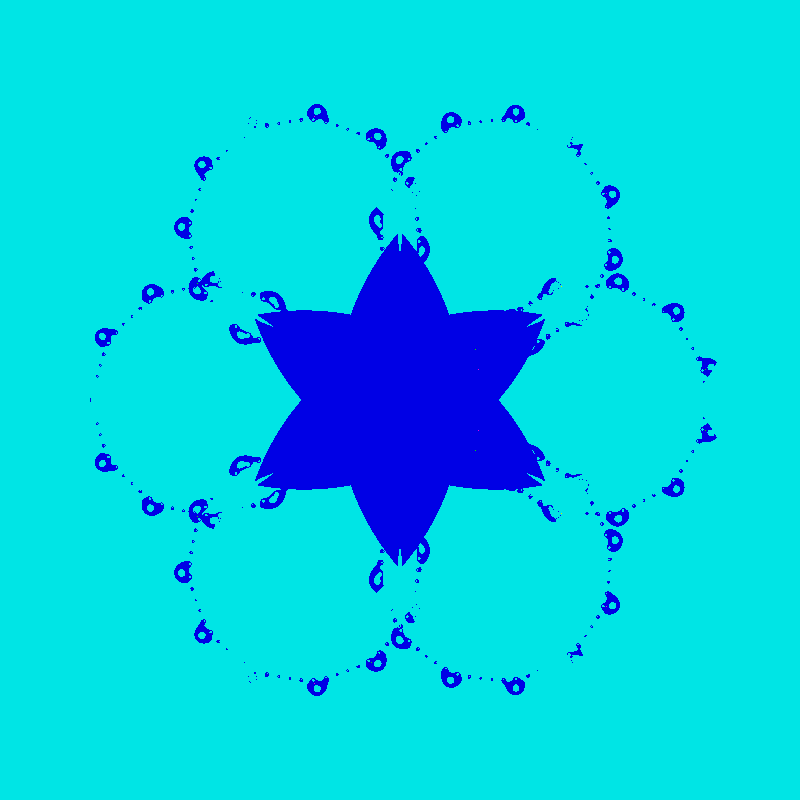}}
    \caption{$(n,d)=(60,12)$ (left), $(n,d)=(12,60)$ (right), and for both, $a=1.5$ and the $c$-plane window drawn is $[-4,4]\times[-4,4]$. Notice that when $d$ is small and $n$ is large, the baby Mandelbrot sets point inward, whereas when $n$ is small and $d$ is large, they point outward. Note that the star artifact at the center of the right picture is not part of the clusters.}
    \label{fig:in/out-param}
\end{figure}

\begin{example} \label{example:cluster-separation}
    Applying Proposition~\ref{prop:cluster-separation} to the examples in our figures yields the following $|a|=\rho({n,d})$-values for spine separation. Note these theoretical values for large $n,d$ are quite different from the actual $a$-values at which we observed cluster separation for various lower values of $n,d$, as described in the figures.
    \begin{itemize}
         \item Figure~\ref{fig:in/out-param}, left $(n,d)=(60,12)$: $|a| =$1.34;
        \item Figure~\ref{fig:in/out-param}, right $(n,d)=(12,60)$: $|a| =$4.29;
        \item Figure~\ref{fig:clusters1}, $(n,d)=(10,5)$: $|a| =$0.48,         \item Figure~\ref{fig:clusters2}, $(n,d)=(20,10)$: $|a| =$0.48;
       \item Figure~\ref{fig:clusters3}, $(n,d)=(20,5)$: $|a| =$1.04
         \item Figure~\ref{fig:diff_in_d/n}, left, $(n,d)=(40,20)$: $|a| =$0.48, and
     \item Figure~\ref{fig:diff_in_d/n}, right, $(n,d)=(20,40)$: $|a| =$0.23
    \end{itemize}
\end{example}

Next, we generalize Lemma~\ref{lem:disk-for-Mna} for the $n\neq d$ case, to find a large disk in the $c$-parameter plane which is guaranteed to contain $\Mnda$.

\begin{lemma} \label{lem:MGMdisk-for-Mdna}
Given $a\in \CC^*$, and $n,d\geq 3$, there exists a $T>0$ so large that if $|c| \ge T$, then $c \notin \Mnda$; that is, all critical orbits escape. Hence, $\Mnda \subset \DD_T$.
\end{lemma}

\begin{proof} This proof is completely analogous to Lemma~\ref{lem:disk-for-Mna}, replacing all $\pm2\sqrt{a}$ with $v$ where $v+c$ denotes any critical value.
\end{proof}

Next, we want to transfer knowledge of the location of the Julia sets to the locations of the boundedness locus. First, a helper lemma.

\begin{lemma} \label{lem:K-in-annulus-lemma}
    Fix $a\in \CC^*$. Suppose for a (compact) collection of $c$-values, call it $\mathcal{S}$, we know that $K(\Fndca)\subset \An{t}{u},$ for some $0<t<u$, and for any collection of $n,d\geq 3$ (infinite collections allowed). Then 
    $$
    \mathcal{S}\cap \Mnda = \{ c \in \mathcal{S} : \Fndca \text{ has a bounded critical orbit } \} \subset 
    \bigcup_{v \in \Vand} \An{t}{u} + v.
    $$
\end{lemma}

\begin{proof}
Let $R > \sup_{c \in \mathcal{S}} |c|$.  So, $\mathcal{S} \subset \DD_R$. 

Similar to the reasoning in the proof of Proposition~\ref{prop:GMparam-fat-annuli}, if $c\in \DD_R$, but not in either annulus, then 
the difference between $c$ to any annulus center is either $>u$ or $<t$, but that difference is exactly the formula for the critical values of the map for that $c$-value: $c\pm v$, for $v\in \Vand$. Hence, for such a $c$ value, the critical values are not in $\An{t}{u}$, hence 
neither can be in $K(\Fnca)$ by our assumption, so both escape, so $c\notin \Mnda \cap \mathcal{S}.$
\end{proof}

We now apply this lemma to a few situations in which we have an annular location of the Julia set. 

\begin{proposition}\label{prop:MGMparam-fat-annuli}
Let  $a\in \CC^*$. Using Lemma~\ref{lem:MGMdisk-for-Mdna},
let $T>0$ be such that $\Mnda \subset \DD_T$ for all $n,d\geq 3$.
Let $s = \max\{3, |a|, T \}$ 
and $p = \left( \frac{|a|}{1+T+s} \right)^{1/3}$. Then for all $n,d \geq 3$, $\Mnda$ is contained in annular regions about the circles $\Cnda$; in particular:
$$
M_{n,d,a} \subset 
\bigcup_{v\in \Vand} v + \An{p}{s}
$$
\end{proposition}

\begin{proof}
 By Proposition~\ref{prop:K-in-annulus} we have $$ K(\Fndca) \subset \An{p}{s},$$
 for the given $a$, for all $c\in \DD_T$, and for all $n,d\geq 3$. So, we combine this with Lemma~\ref{lem:K-in-annulus-lemma}, taking $\calS := \DD_T$, we know that (for any $n,d\geq3$): $$\calS \cap \Mnda \subset \bigcup_{v\in\Vand}v+\An{p}{s}.$$ 

 Since $\Mnda \subset \DD_T$, we get $\calS \cap \Mnda = \DD_T \cap \Mnda = \Mnda$. Thus we have 
 $\Mnda \subset \bigcup_{v\in\Vand}v+\An{p}{s}.$ 

 To verify this, refer to Lemma~\ref{lem:mgm-c-plane-spine}.
\end{proof}

\begin{figure}[!ht]
    \centering
    \includegraphics[width=0.3\textwidth,trim=0 0 100 0,clip]{\detokenize{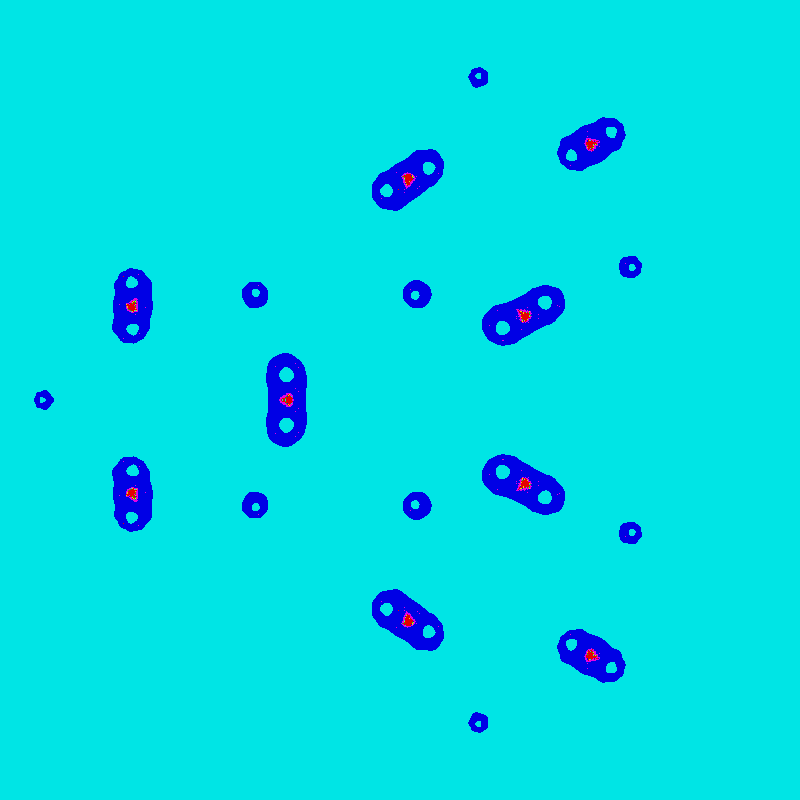}}\hfill
    \includegraphics[width=0.3\textwidth,trim=0 0 100 0,clip]{\detokenize{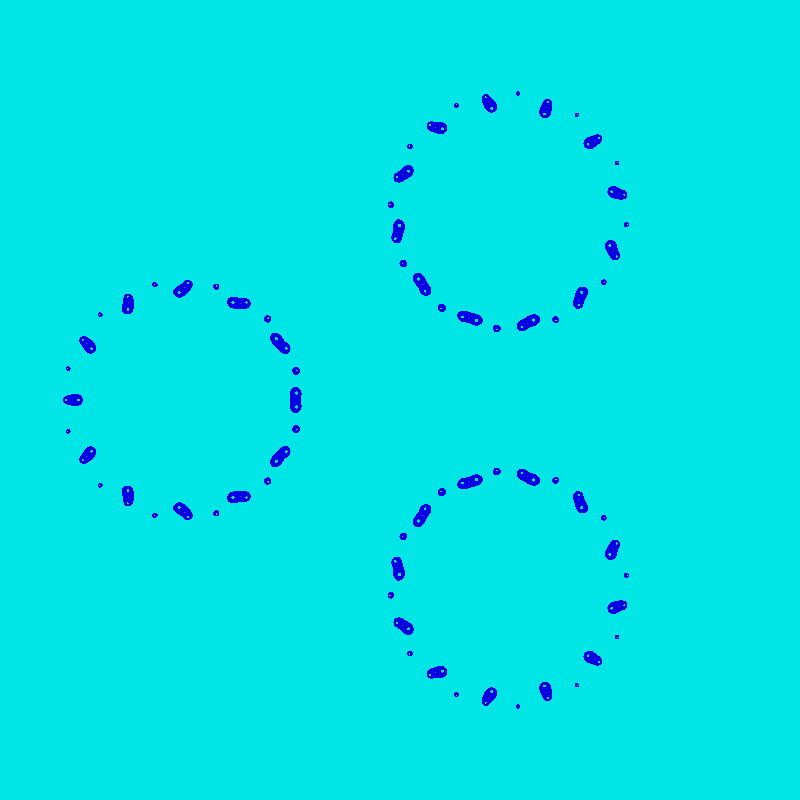}}\hfill
    \includegraphics[width=0.3\textwidth,trim=0 0 100 0,clip]{\detokenize{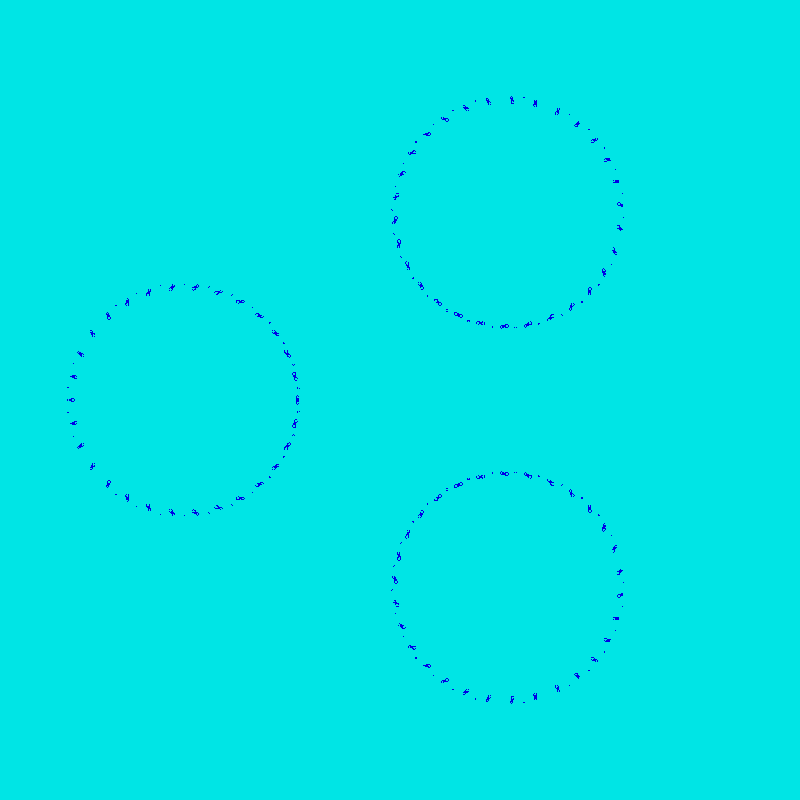}}
    \caption{$n>d$ with $a=1$. From left to right: $(n,d)=(8,4)$, $(24,12)$, and $(60,30)$, where the $c$-window shown is $[-3.5,3.5]\times[-3.5,3.5]$. Here, $\ell = \gcd(n,d)=d$ in each case and $\nu=3$ is constant. Since $\ell \to\infty$, the boundedness locus fills in the three circles in $\mathcal{C}$ in the limit.}
    \label{fig:circles}
\end{figure}

\begin{figure}[!ht]
    \centering
    \includegraphics[width=0.24\textwidth,trim=0 0 0 0,clip]{\detokenize{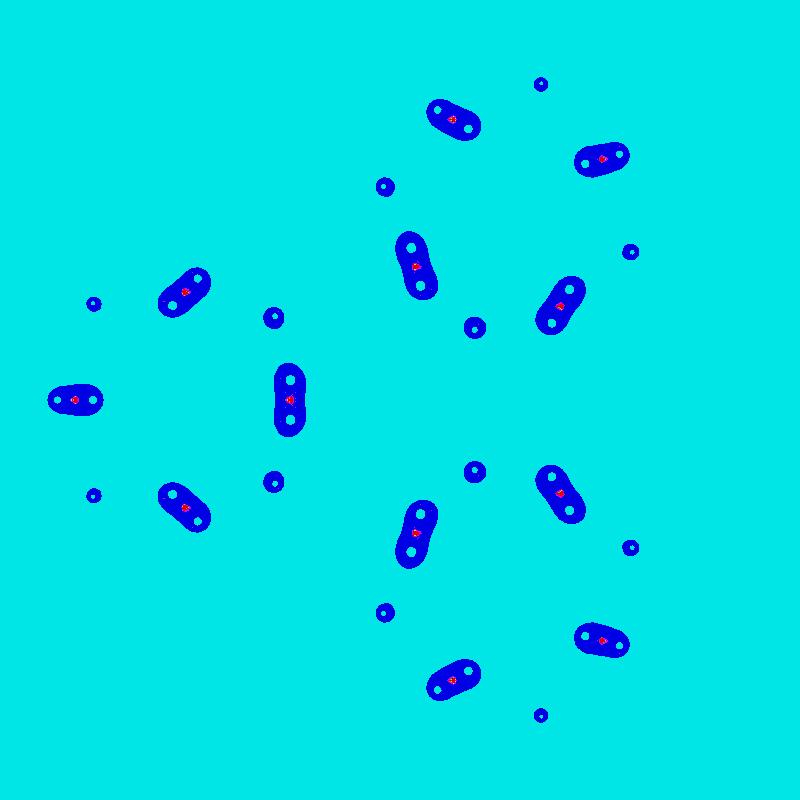}}\hspace{0.15cm}%
    \includegraphics[width=0.24\textwidth,trim=0 0 0 0,clip]{\detokenize{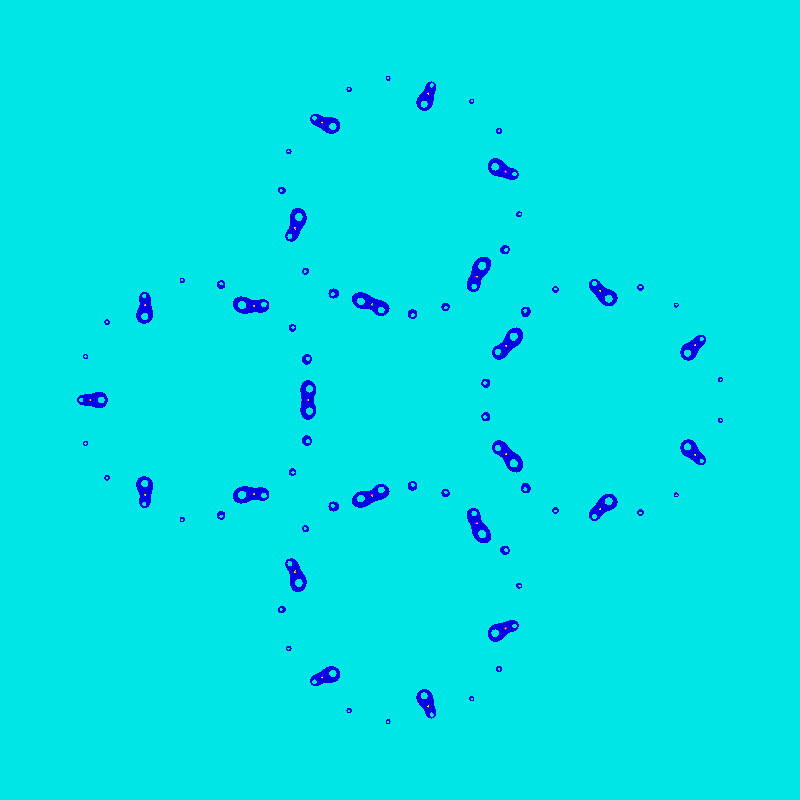}} \hspace{0.1cm}%
    \includegraphics[width=0.24\textwidth,trim=0 0 0 0,clip]{\detokenize{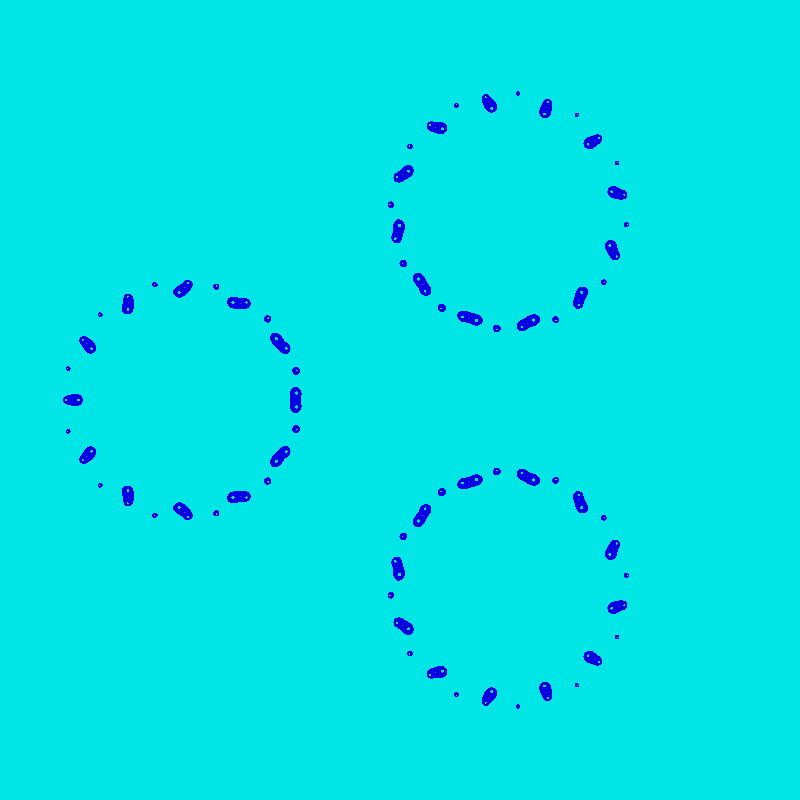}}\hspace{0.15cm}%
    \includegraphics[width=0.24\textwidth,trim=0 0 0 0,clip]{\detokenize{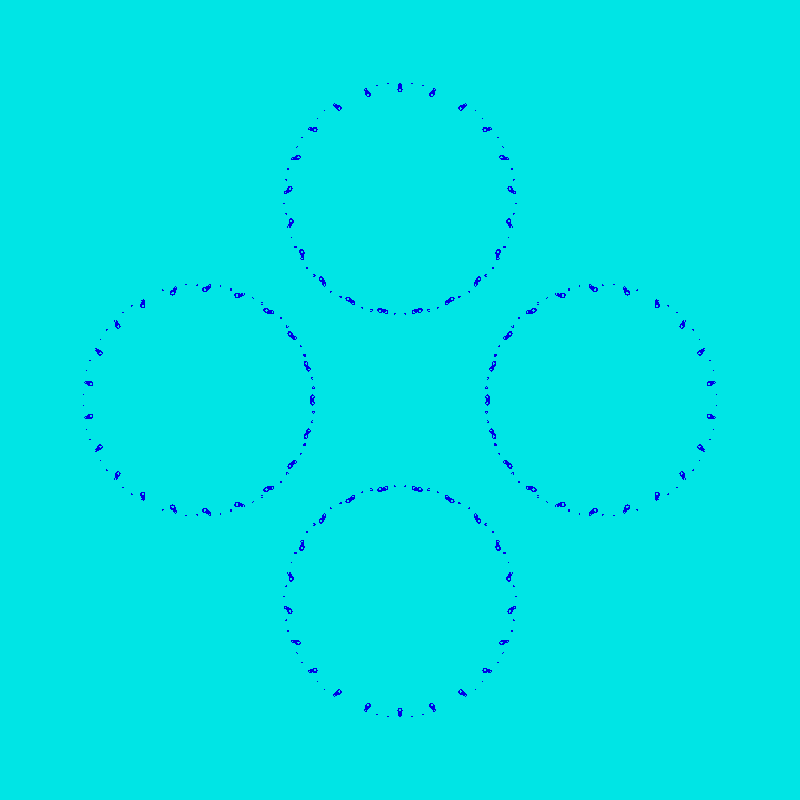}}
    
    \caption{$n>d$ with $a=1$. From left to right: $(n,d)=(6,3), (18,6), (24,12),$ and $(63,21)$; $c$-windows all $[-3.5,3.5]\times[-3.5,3.5]$. The sequence satisfies $n=2d$ for odd indicies, and $n=3d$ for even indices; in all cases $\gcd(n,d)=d$. As $\gcd(n,d)\to\infty$, the baby Mandelbrot sets fill each $\mathcal{C}$, while the number of circles in $\mathcal{C}$ alternating between three and four along the sequence.}
    \label{fig:zig_zag}
\end{figure}

\begin{example} \label{example:p-s-param}
    Applying Proposition~\ref{prop:MGMparam-fat-annuli} to the examples in our figures yields the following  annuli about the $v_k$'s which contain the boundedness loci: 

    \begin{itemize}
        \item Figure~\ref{fig:lemma5.3}:   $p = 0.391, s = 3.5,$ (for $a = 0.48$),
        \item Figure~\ref{fig:in/out-param}: $p = 0.55, s = 4,$ (for $ a = 1.5$), 
        \item Figure~\ref{fig:circles}: $p = 0.5, s = 3.5,$ (for $a = 1$), and
        \item Figure~\ref{fig:zig_zag}: $ p = 0.5, s = 3.5,$ (for $a = 1$). 

    \end{itemize}
\end{example}

Next, we use the information in Proposition~\ref{prop:MGMparam-fat-annuli} to squeeze $\Mnda$ toward $\Cnda$ as $n,d$ grow.

\begin{proposition}
    \label{prop:MGMparameps}
    Let $a\in \CC^*$ be fixed and $0 < \eps < 1$. There is an $N\geq 3$ such that for all $n,d \geq N$, we have 
    \[
    \Mnda \subset
    \calN_\eps(\Cnda).
    \]
    
Observe 
$$ \calN_\eps(\Cnda)
= \bigcup_{v\in \Vand} (v + \An{1-\eps}{1+\eps}).$$
\end{proposition}

\begin{proof} 
First, by Lemma~\ref{lem:MGMdisk-for-Mdna},  there is a $T>0$  so large that $\Mnda \subset \DD_T$.
 By Proposition~\ref{prop:in+outJ1}, for the given $a$, for all $c$ with $|c|\le T$, there exists an $N$ such that for all $n,d\geq N$, we have: 
 $$ K(\Fndca) \subset \An{1-\eps}{1+\eps}.$$
 Then, similarly to the proof of Proposition~\ref{prop:MGMparam-fat-annuli}, by Lemma~\ref{lem:K-in-annulus-lemma}, we  take $\calS := \DD_T$ and get 
$$\calS \cap \Mnda \subset \bigcup_{v\in \Vand} (v + \An{1-\eps}{1+\eps}) = \calN_\eps(\calC_{a,n,d}).$$ 
Notice we have $\calS \cap \Mnda =\DD_T \cap \Mnda = \Mnda$ (since $\Mnda \subset \DD_T$), and so
$ \Mnda \subset \calN_\eps(\calC_{a,n,d}).$ 
\end{proof}

The above immediately implies the following. 

\begin{corollary} \label{cor:GM-geom-limit}
    If $n_i,d_i$ are sequences that $\to\infty$ as $i \to\infty$, then
$$\lim_{i\to \infty} M_{n_i,d_i,a} \subset \lim_{i \to \infty} \calC_{n_i,d_i,a},$$
if and only if the limit on the right hand side exists.
\end{corollary}

Now, we can also ``mix and match'' the above to yield the following corollaries on annuli containing $\Mnda$. 

\begin{corollary}
    \label{cor:MGM-dbig}
Let $a\in\CC^*, \eps>0$. 
Let $T>0$ be such that $\Mnda \subset \DD_T$ for all $n,d\geq 3$, as guaranteed by Lemma~\ref{lem:MGMdisk-for-Mdna}.
Let $s = \max\{3, |a|, T \}$.

    Then there is a $D\geq 3$ such that for all $d\geq D$, and all $n\geq 3$ we have 
    \[
    \Mnda\subset \  \bigcup_{v\in \mathcal{V}_{a,n,d}} v+\An{1-\eps}{s}. 
    \]
\end{corollary}

\begin{proof}
Combine Corollaries~\ref{cor:uniformbound1} and~\ref{cor:inside-common-J}  
 as in the proof of Propositions~\ref{prop:MGMparam-fat-annuli} and~\ref{prop:MGMparameps}.
\end{proof}

\begin{corollary}
    \label{cor:MGM-nbig}
    Let $a\in\CC^*, \eps>0$. 
Let $T>0$ be such that $\Mnda \subset \DD_T$ for all $n,d\geq 3$, as guaranteed by Lemma~\ref{lem:MGMdisk-for-Mdna}.
Let $s = \max\{3, |a|, T \}$ 
and $p = \left( \frac{|a|}{1+T+s} \right)^{1/3}$.

    Then there is an $N\geq 3$ such that for all $n\geq N$, and all $d\geq 3$ we have 
    \[
    \Mnda \subset \ \bigcup_{v\in \mathcal{V}_{a,n,d}} (v+\An{p}{1+\eps}). 
    \]
\end{corollary}

\begin{proof}
Similar to the above proof, we combine 
Corollary~\ref{cor:inside-J-RQ} with Corollary~\ref{cor:outside-J-common-N}.  
\end{proof}

Observe Figure~\ref{fig:in/out-param} to visualize Corollaries ~\ref{cor:MGM-dbig} and~\ref{cor:MGM-nbig}.

\bigskip

Next, we establish an analog of Proposition~\ref{prop:GM-geom-lim-param} in the $n\neq d$ case. 
In this case, the spine defined in Lemma~\ref{lem:mgm-c-plane-spine} changes with $n$ and $d$; compare to the $n=d$ case where the spine is fixed. Still, our analog result shows that $\Mnda$ approaches this moving spine as $n,d\to \infty$. 

\begin{proposition}
    \label{prop:MGM-geom-lim-param}
    Let $ a\in \CC^*$ and  $\eps>0$. Suppose  $n_i, d_i$ are sequences such that   $\ell_i = \gcd(n_i,d_i)\to\infty$ as $i \to \infty$,
    and $\{ \chi_i \}$ is a complex-valued sequence such that $\chi_i \in \calC_{n_i,d_i,a}$. 
    Then there exists an $N$ such that for all $ i \geq N$ 
    \[
    B_\eps(\chi_i)\cap M_{n_i, d_i, a} \neq \emptyset.
    \]
\end{proposition}

\begin{proof}

    We must establish where $M_{n_i, d_i, a}$ tends as $i\to \infty$. For simplicity, we focus on one set of $n,d$ values at a time, so we do not include the $i$ index below. As in the proof of Proposition~\ref{prop:GM-geom-lim-param}, we  begin by finding some $c$ that are always contained in $\Mnda$; maps with at least one critical fixed point.  
    We refer to these values as the ``centers'' of $\Mnda$, as for any baby Mandelbrot in the boundedness locus, one such value is at its center.

Recall from Notation~\ref{not:Vand}:
$$
\Vand = \left\{ v_k : k=0,\ldots,m-1
\right\}.$$

Now, for each $k \in \{ 0, \ldots, m-1\}$, we want to find a $c$-value, $c_k$, such that for the map $F_{n,d,a,c_k}$, the critical point $w_k$ is fixed; that is, $F$ maps $w_k$ to its appropriate critical value, so
 $ F(w_k) = w_k.
$ 
Equivalently (recall Lemma~\ref{lem:MGMcritvalues}), 
$$
v_k +c_k = w_k, \text{ hence } 
c_k  =-v_k+w_k.
$$

Recall $\Cnda = \bigcup_{k=0}^{m-1} \Cnda^k$ is the union of $\nu$ circles of radius 1 with centers at the points $v_k$ in $\Vand.$

Notice the above equation 
 tells us that the distance from 
to the center of $\Cnda^k$ (the $v_k$'s) to the $c_k$'s is $|w_k|.$ 
By the Section~\ref{sec:dynamicalplaneresults} results, since the Julia sets approach the unit circle as $n,d \to \infty$, we know $|w_k| \to 1$ as $n,d\to \infty$ (remember, $a$ is fixed).

Next we describe the symmetry of $c_k$-values about $\Cnda^k$.
Suppose $\chi $ is in the $k^{th}$ circle, denoted by $\Cnda^k$. 

By Lemma~\ref{lem:rotsym1} we have that for any $z^*$, for each $j = 0, \ldots, \ell-1$: 
$$
F\left(z^*e^{i\frac{2\pi}{m}\left(k+\frac{m}{\ell}j\right)}\right)
=
F\left(z^*e^{i\frac{2\pi}{m}k}\right).
$$ 
Thus for $z^*= \left(\frac{d|a|}{n}e^{i\psi}\right)^{1/m} $, we have $w_k = z^* e^{i\frac{2\pi}{m}k},$ thus
$
F\left(z^*e^{i\frac{2\pi}{m}\left(k+\frac{m}{\ell}j\right)}\right)
= F(z^* e^{\frac{i 2\pi k}{m}} ) = F(w_k),  
$
so 
\begin{eqnarray*}
c_k + v_k  & =&  F(w_k)  
= F\left(z^*e^{i\frac{2\pi}{m}\left(k+\frac{m}{\ell}j\right)}\right)
= z^*e^{i\frac{2\pi}{m}\left(k+\frac{m}{\ell}j\right)}
\\ &= & 
 \left(\frac{d|a|e^{i\psi}}{n}\right)^{1/m} e^{i\frac{2\pi}{m}\left(k+\frac{m}{\ell}j\right)}
 = \left(\frac{d|a|\phi_k}{n}\right)^{1/m} e^{i\frac{2\pi}{m}\left(\frac{m}{\ell}j\right)}
= (w_k) \ e^{i\frac{2\pi}{\ell}j}.
\end{eqnarray*}

Thus the equation for all $c$-values with $w_k$ as the fixed critical point is:
\begin{eqnarray*}
    c_{k,j} 
    &=& -  \frac{m}{d}\left(w_k\right)^{n}  
+ \left(w_k\right) e^{i\frac{2\pi}{\ell}j}
\end{eqnarray*}
for each $j = 0, \ldots, \ell-1$.

Notice this tells us that there are $\nu= \frac{m}{\ell}$ clusters of $c_{k,j}$-values, where the $k^{th}$ cluster has $\ell$ number of $c_{k,j}$-values which lie in a circle of modulus $|w_k|$ about the center of $\Cnda^k$, and are evenly spaced by the angle $\frac{2\pi}{\ell}$. The centers $-v_k$ of these $\nu$ circles $\Cnda^k$ themselves lie in a circle about the origin of radius $r_{n,d,a} = |v_k|$. 

Since $\chi \in \Cnda^k$, let $c^*$ be the point $c_{k,j}$ which is closest to $\chi$, so $c^* = c_{k,j}$ for some $j\in \{0,\ldots, \ell-1\}$. 
Since $\chi=- \frac{m}{d}(w_k)^{n}+ e^{i\theta}$,
we have $\chi - c^* = e^{i\theta} - \left(w_k\right)e^{i\frac{2\pi}{\ell}j}. $ 
We already know that $|w_k| \to 1$ as $n,d\to \infty$; then we just need the $e^{i\frac{2\pi}{\ell}j}$ values for all $j\in \{0, \ldots, \ell-1\}$ to fill up the unit circle as $n,d\to\infty$, which occurs precisely when $\ell\to \infty$. 
\end{proof}

Combining the above results, we are ready to prove Theorem~\ref{thm:parameter_n_d_to_inf}.

    \begin{proof}[Proof of Theorem \ref{thm:parameter_n_d_to_inf}]
The equality of the limits follows from Propositions ~\ref{prop:MGMparameps} and ~\ref{prop:MGM-geom-lim-param}, if they exist.
  
Examining Propositions~\ref{prop:MGMparameps} and~\ref{prop:MGM-geom-lim-param}, we conclude that for a sequence of $n,d$ such that $\lim_{n,d \to \infty} \calC_{n,d,a}$ exists, the limit of $M_{n,d,a}$ must also exist (the sequence must be Cauchy). 

But if the right hand side does not converge, neither can the sequence of the boundedness loci. 
    \end{proof}

Thus, success of Theorem \ref{thm:parameter_n_d_to_inf} depends on the particular sequences $n_i, d_i$. We provide an example where the limits exist and one where they do not.

\begin{example}
\label{ex:Mfillspine}
Let $d_i$ be any sequence such that $d_i\to \infty$ as $i\to\infty$, and set $n_i:= 2 d_i$. Then $m_i = n_i+d_i = 3d_i$, $\ell_i=\gcd(n_i,d_i) = d_i$ (so $\ell_i \to \infty$ so these sequences satisfy the hypotheses of Theorem~\ref{thm:parameter_n_d_to_inf}), and
$\nu_i  = \frac{m_i}{\ell_i}= 3.$ 
Thus,
$
v_k = 
3\left(\frac{1|a|}{2}\right)^{2/3} e^{i (\psi + 2\pi k) \frac{2}{3}},$ with unique values for $k=0,1,2$.
Hence, $\calC_{n_i,d_i,a}$ is always three circles of radius $1$ with centers at 
$$\{ v_0, v_1, v_2 \} 
= \left\{ 
3\left(\frac{|a|}{2}\right)^{2/3}, \ 
3\left(\frac{|a|}{2}\right)^{2/3} e^{i (\psi + 2\pi) \frac{2}{3}}, \ 
3\left(\frac{|a|}{2}\right)^{2/3} e^{i (\psi + 4\pi)\frac{2}{3} } \right\}.$$
These do not change with $n,d$, so $\calC_{n_i,d_i,a}$ is the same for all $i$. Thus $\lim_{i\to\infty} \calC_{n_i,d_i,a}$ exists. 

Figure~\ref{fig:circles} illustrates such an example: for $a=1$ and $d_i = 2i$, the sets $M_{n_i,d_i,a}$ are drawn for $i=2, 6, 15$, with $n_i=2d_i$ (note we use $i$ simply as an index here, not as the imaginary number).

Of course, there is nothing special about the $2$ in $n_i = 2 d_i$; any integer $N$ would suffice and give a constant number of circles $N+1$, and the same result. 
\end{example}

\begin{example}
\label{ex:Mdoesntfillspine}
    To generate an example of sequences $n_i, d_i$ where the limits of $\{ \calC_{n_i,d_i,a} \}$ do not exist, one could perform a ``zipper merge'' on two different sequences defined like the prior example. That is, let $d_i$ be any sequence such that $d_i\to \infty$ as $i\to \infty$. For odd $i$ let $n_i: = 2 d_i$, and for even $i$ let $n_i := 3 d_i$. Then the sequence of sets $\calC_{n_i,d_i,a}$ alternates between three circles and four circles, thus the limit does not exist. See Figure~\ref{fig:zig_zag}. Here, $a=1$. The sets $M_{n_i,d_i,a}$ are shown with $n_i=2d_i$ for the first and third elements and $n_i=3d_i$ for the second and fourth elements.
\end{example}

Now we observe from what direction in parameter space the boundedness loci tend to the spine.

\begin{corollary}
    Let $n_i, d_i$ be sequences such that $n_i, d_i \to \infty$ as $i\to\infty$. Consider an $i\in \NN$.
    \begin{itemize}
        \item  If $|a| < \frac{n_i}{d_i},$ then the ``centers'' of $M_{n_i,d_i,a}$ 
    (the $c$-values with fixed critical point $w_k$) are interior to the spine circles $\Cnda$;
    \item  if  $|a| > \frac{n_i}{d_i},$ the centers of $M_{n_i,d_i,a}$ are exterior to the spines; 
    \item if $|a| =  \frac{n_i}{d_i},$ the centers of $M_{n_i,d_i,a}$ lie on the spines. 
    \end{itemize}
   Thus, 
    if there is an $I$ such that $|a|<\frac{n_i}{d_i}$ for all $i\geq I$, then $\Mnda$  approaches $\Cnda$ generally from the interior of these circles, 
    whereas if $|a|>\frac{n_i}{d_i}$ for all $i$ larger than some $I$, then they approach from the exterior. 
\end{corollary}

\begin{proof}
     Suppose $c^*=c_{k,j}$ for some fixed $k$ and arbitrary $j.$ Then we know $|w_k|=\left(\frac{d_i|a|}{n_i}\right)^{1/m_i}$ is the distance from the center of the $k^{th}$ circle to $c^*.$ So if we compare this to $1$, we need only know the value of $\frac{d_i|a|}{n_i}$, which results in our comparison of $|a|$ to $\frac{n_i}{d_i}.$
\end{proof}

For example, in Examples~\ref{ex:Mfillspine} and~\ref{ex:Mdoesntfillspine}, we have $|a|=1 < \frac{n_i}{d_i}=2$ or $3$ for all $i$, and we observe that the distance from the centers to the baby Mandelbrots  grows to $1$ as $i$ increases, i.e. the circles expand up into the spines.  
\medskip

\section{Additional Pictures}
\label{sec:pics}

In this section, we provide a few additional observations and visuals of $\Mnda$ in parameter space. We also explain our methodology for creating these visualizations, and provide some comments on the interplay between the theory and the computer-generated images.

\subsection{Cluster Separation Observations}
First, as mentioned earlier, Meyer, Stoertz, and Wang proved  if $|a|$ is sufficiently large, there are $\ell=\gcd(n,d)$ baby Mandelbrots per cluster.
Though it is difficult to discern in most of our images, in all of the $c$-planes where $a$ is large enough to achieve some good separation of the clusters, magnification appears to reveal baby Mandelbrot sets strewn throughout the clusters. In fact, if a baby Mandelbrot exists, its ``center'' is a parameter $c$ where a critical point is a fixed point; recall that finding such parameters was our strategy to prove Proposition~\ref{prop:MGM-geom-lim-param}. Moreover, there are, in some of our Julia set pictures, intricate baby Julia set shapes which appear upon magnification. See Figure~\ref{fig:p_to_d}.

Now, Propositions~\ref{prop:c-plane-spine} and~\ref{prop:cluster-separation} provide an $|a|$-value for which the spines of the boundedness locus, the collection of $\nu=(n+d)/\gcd(n,d)$ circles, separate in the limit under degree growth---for large $n$ and/or $d$. For finite $n,d$, the boundedness locus has a thickness around these circles, thicker for lower $n,d$. We provide below some images and observations for cluster separation for various $|a|,n,d$-values in practice. 

\begin{example}
    \label{itm:equal}
    In the case $n=d$, recall there are always two clusters (because $\nu=(n+d)/\gcd(n,d)=2n/n=2)$. We observe that smaller $n=d$ values correspond to larger $|a|$ values needed to achieve separation Equivalently, the $|a|$-value needed to achieve separation decreases (to the limit of $0.25$) as $n=d$ increases to $\infty$. For example, cluster separation occurs about $a=0.45$ for $n=d=10$, $a=0.35$ for $n=d=20$. See Figure~\ref{fig:two_clusters2}.

\begin{figure}[!ht]
    \centering
    \includegraphics[width=0.3\textwidth]{\detokenize{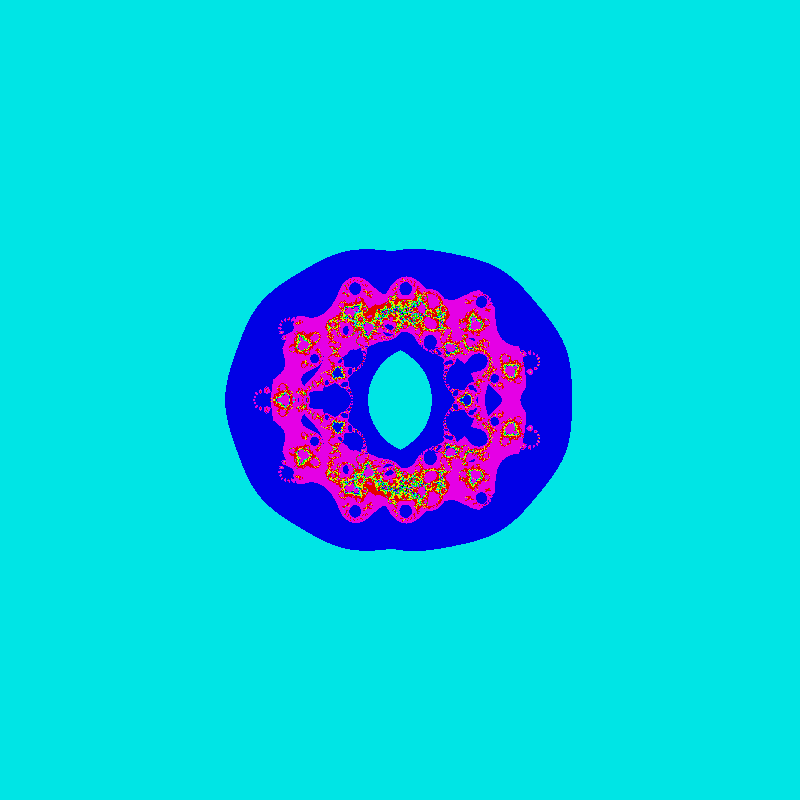}}\hfill
    \includegraphics[width=0.3\textwidth]{\detokenize{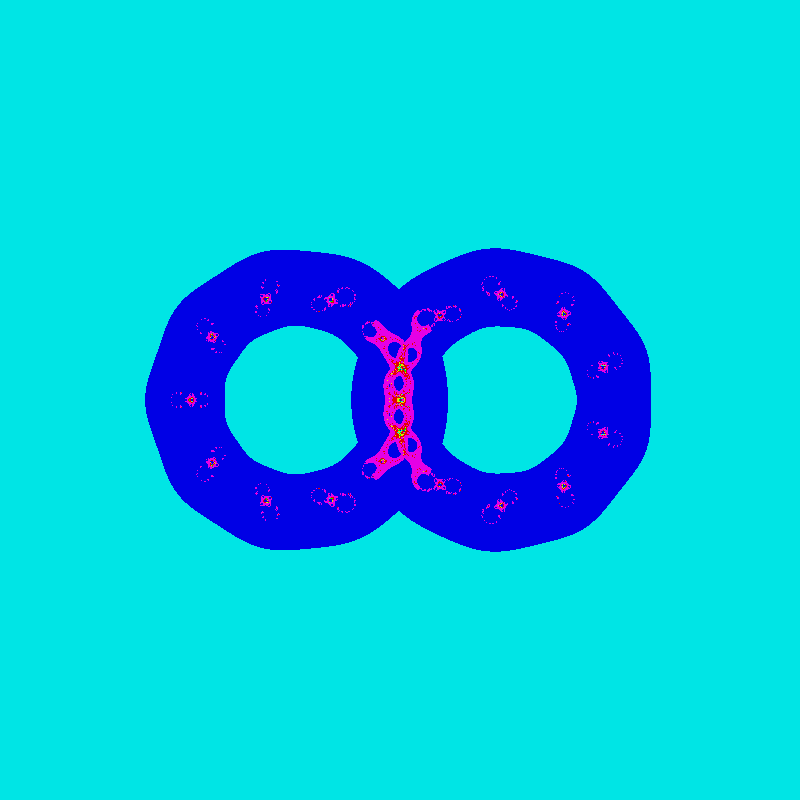}}\hfill
    \includegraphics[width=0.3\textwidth]{\detokenize{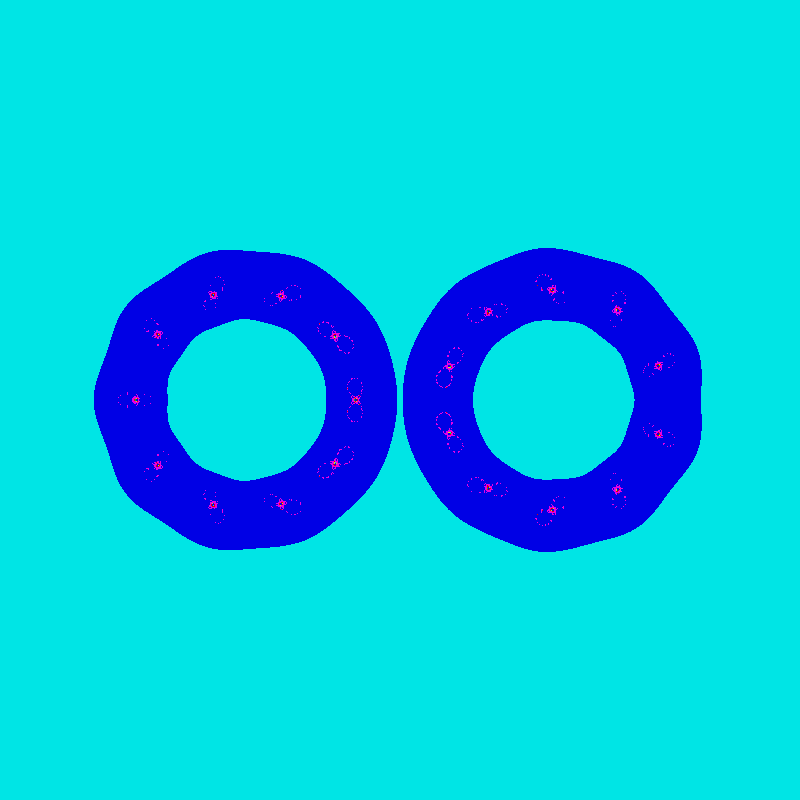}}\hfill
%
      \includegraphics[width=0.3\textwidth]{\detokenize{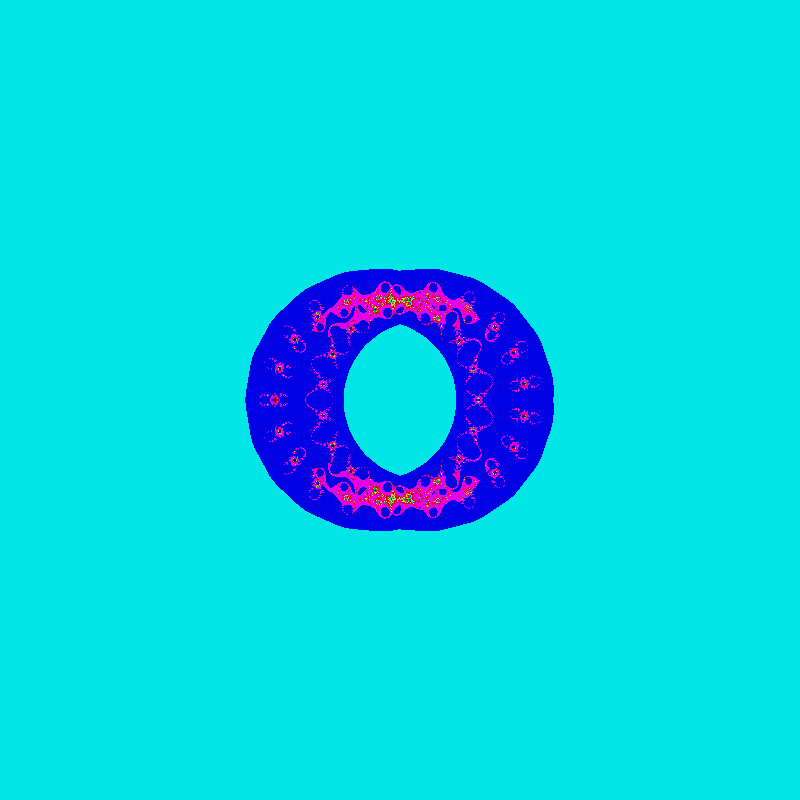}}\hfill
    \includegraphics[width=0.3\textwidth]{\detokenize{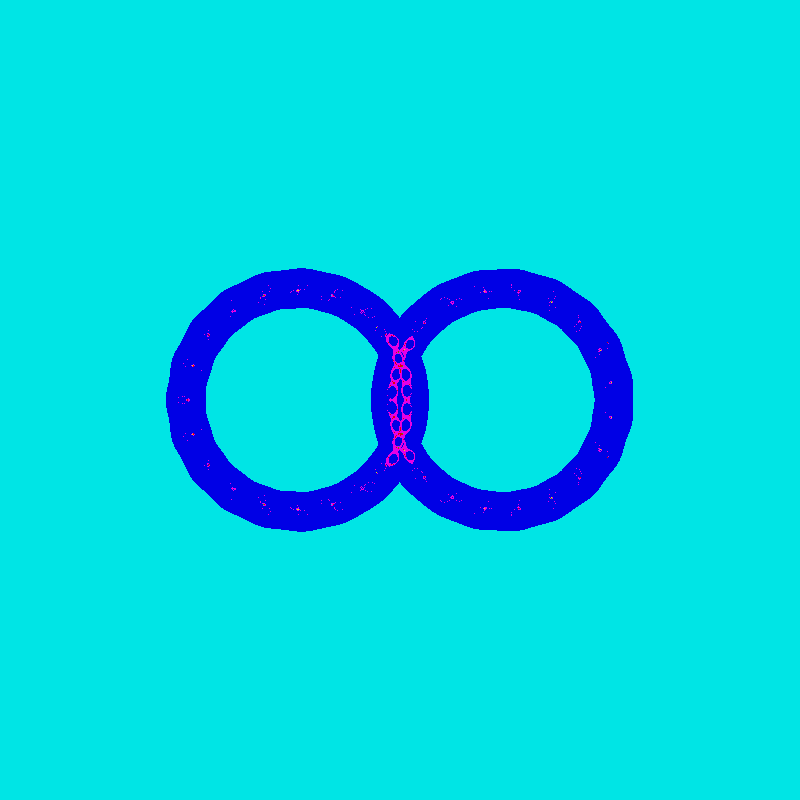}}\hfill
    \includegraphics[width=0.3\textwidth]{\detokenize{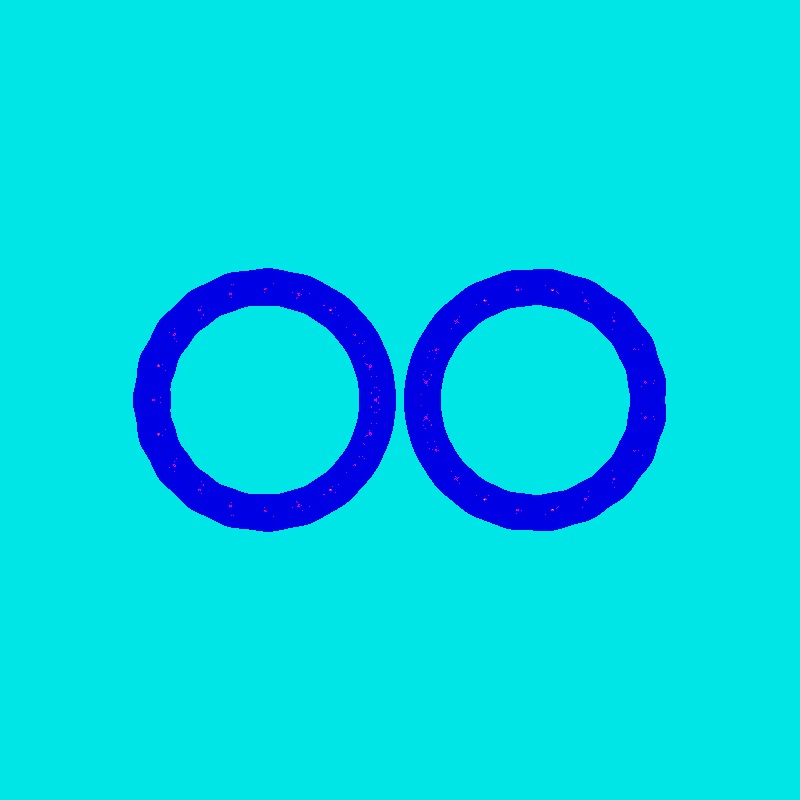}}\hfill
    \caption{Top row: $n=d=10$, and (left to right) $a=0.01,0.2,0.45$. Bottom row is $n=d=20$, and $a=0.01,0.2, 0.35$ (left to right). 
     (The $c$-window is $[-3.5,3.5]\times[-3.5,3.5]$ for all.)}
    \label{fig:two_clusters2}
\end{figure}

\end{example}

\begin{example} \label{ex:n-greater-d}
Next, we begin to examine the case of $n\neq d$. Recall we have $\nu={n+d}/{\gcd(n,d)}$ clusters or annuli, whose center points lie on a circle. In the $n=d$ case we found that increasing $n=d$ led to separation at a smaller $|a|$-value. However, when $n\neq d$, the number of clusters seems to have an effect.

\begin{enumerate}[label=(\roman*),ref=(\roman*)]

    \item \label{itm:n-greater-d-limit}
First, we show an example where a larger degree comes with also a larger $|a|$ to achieve separation (the opposite of what we observed in the $n=d$ case). Suppose we keep $d$ fixed but increase $n$ in such a way that the number of clusters, $\nu$ increases. 
Compare Figure~\ref{fig:clusters1}, $d=5, n=10,$ with separation at $ a=0.8$, to Figure~\ref{fig:clusters3}, which is also $d=5$ but larger $n=20$, however separation at larger $a=1.2.$

        \item \label{itm:n-greater-d-fixed}
    On the other hand, compare Figures~\ref{fig:clusters3} and~\ref{fig:clusters2}. Both have the same $n$, but the $d$ is larger in the second figure, Figure~\ref{fig:clusters2}, and in this case there is a smaller number of clusters in the second figure and we see separation occurs at a smaller $|a|$, as we would expect for increasing degree.

\begin{figure}[!ht]
    \centering
    \includegraphics[width=0.3\textwidth]{\detokenize{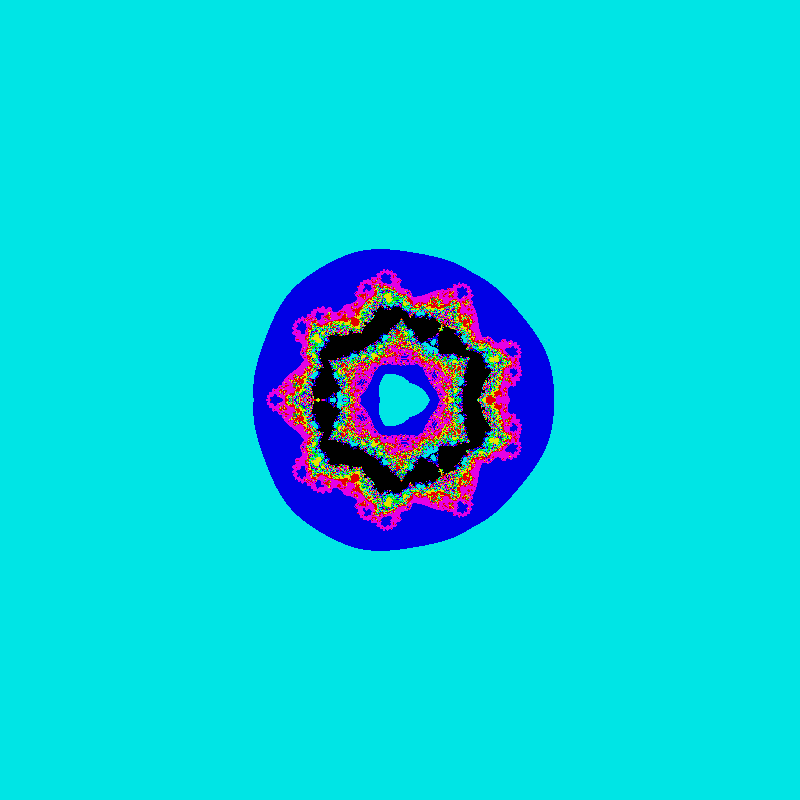}}\hfill
    \includegraphics[width=0.3\textwidth]{\detokenize{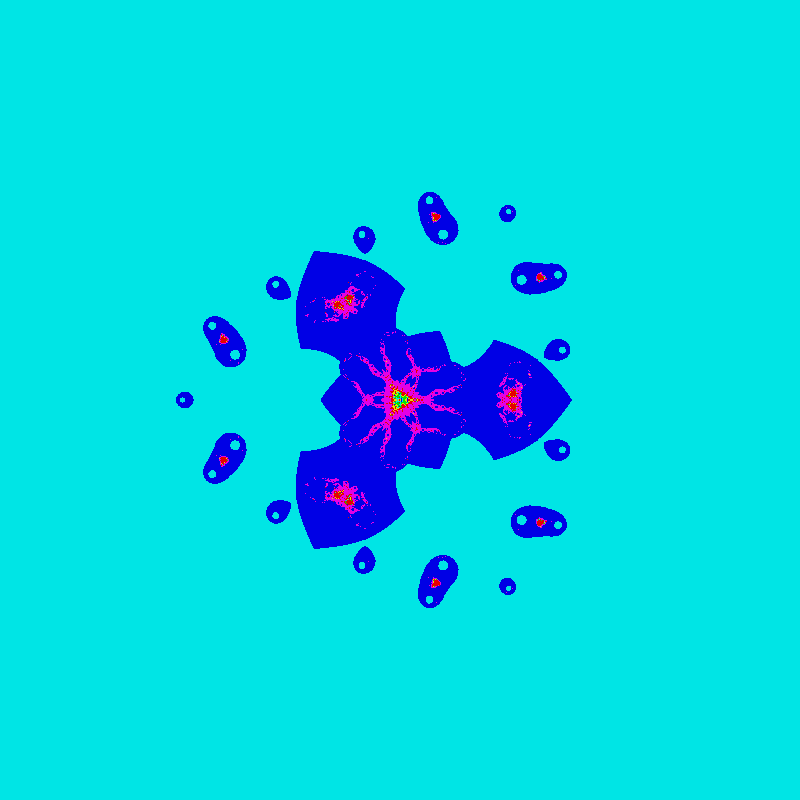}}\hfill
    \includegraphics[width=0.3\textwidth]{\detokenize{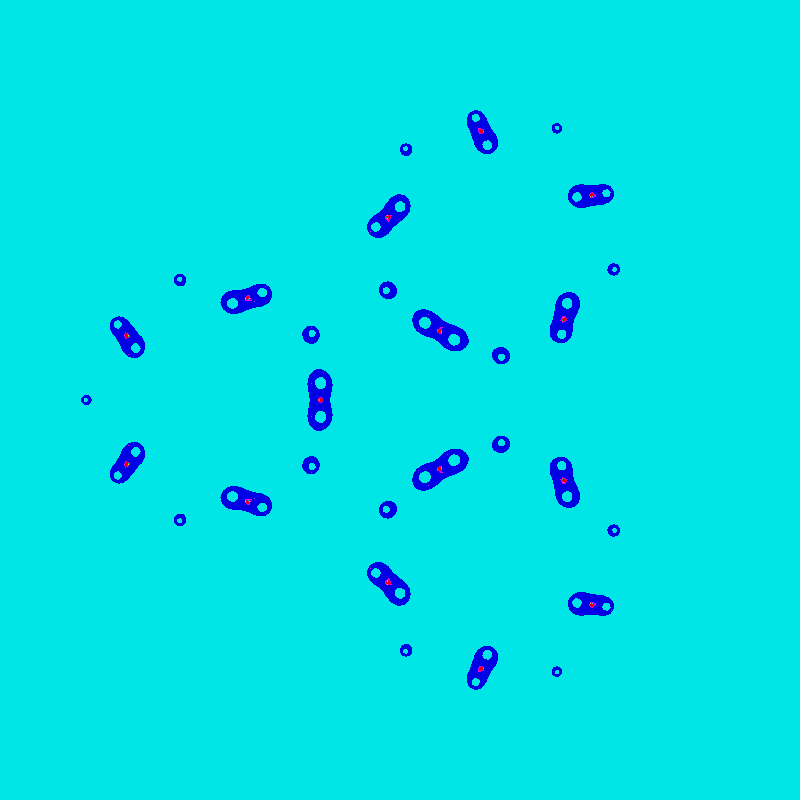}}\hfill
    \caption{Here $d=5, n=10$, $a=0.01,0.3,0.8$ ($c$-window: $[-3.5,3.5]\times[-3.5,3.5]$). The clusters separate just before $a=0.8$.}
    \label{fig:clusters1}
\end{figure}

\begin{figure}[!ht]
    \centering
    \includegraphics[width=0.3\textwidth]{\detokenize{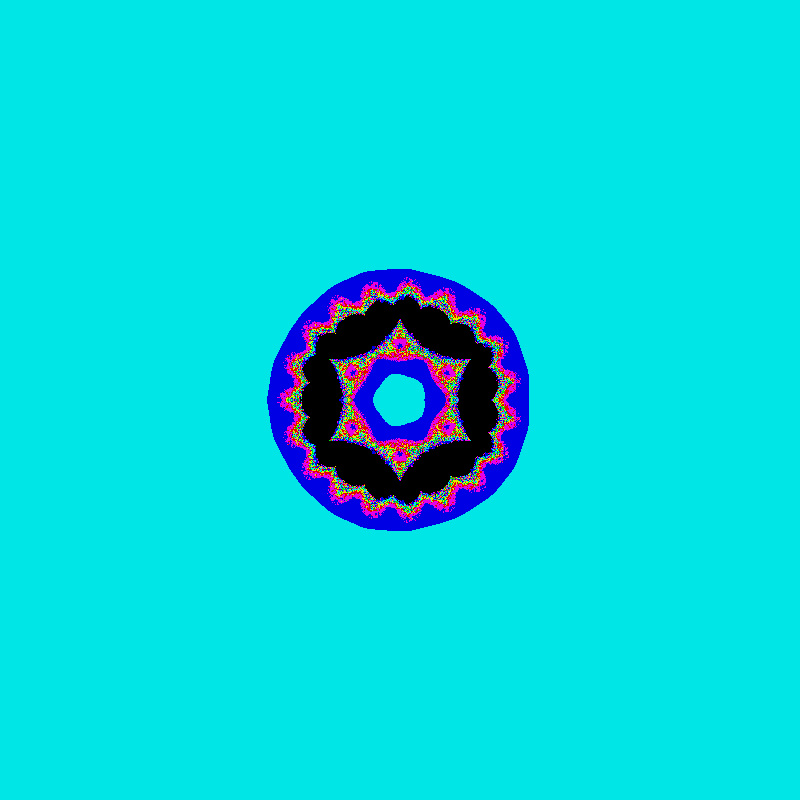}}\hfill
    \includegraphics[width=0.3\textwidth]{\detokenize{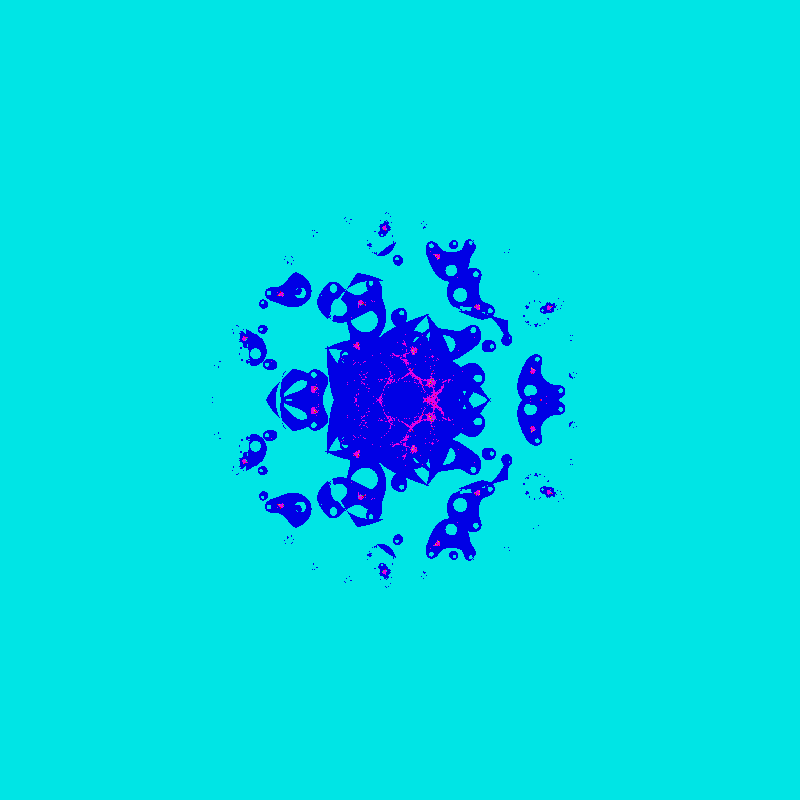}}\hfill
    \includegraphics[width=0.3\textwidth]{\detokenize{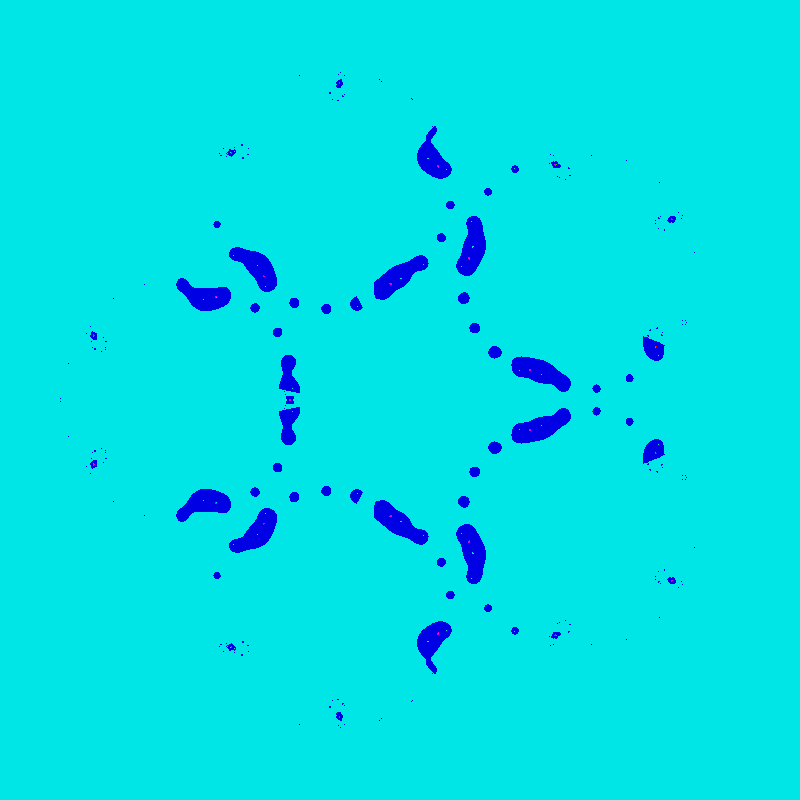}}\hfill
    \caption{Here $d=5, n=20$, and $a=0.01,0.3,1.2$ ($c$-window: $[-3.5,3.5]\times[-3.5,3.5]$). The five clusters separate just before $a=1.2$. }
    \label{fig:clusters3}
\end{figure}

\begin{figure}[!ht]
    \centering
    \includegraphics[width=0.3\textwidth]{\detokenize{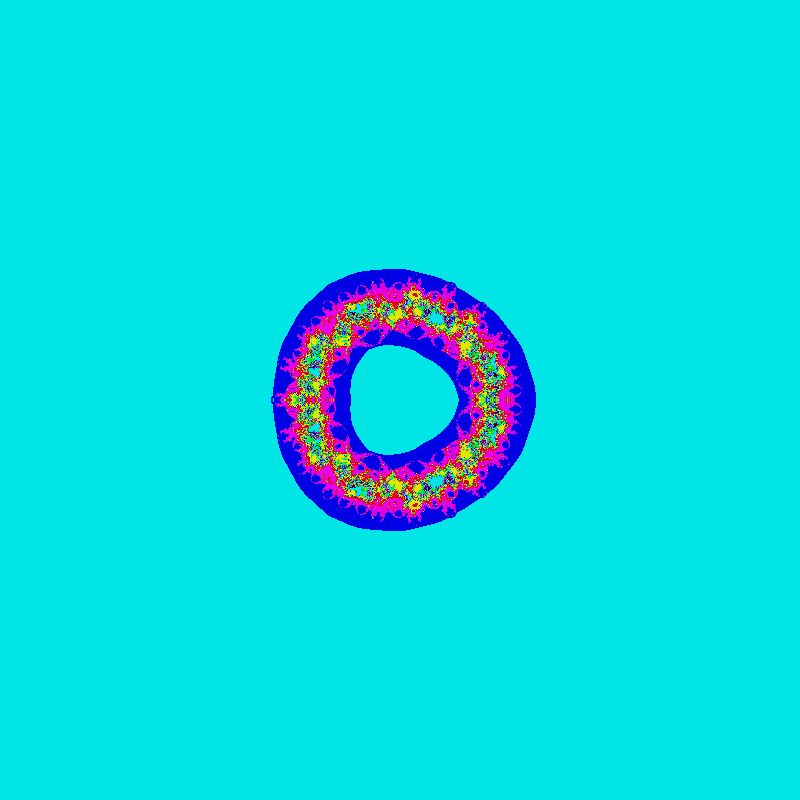}}\hfill
    \includegraphics[width=0.3\textwidth]{\detokenize{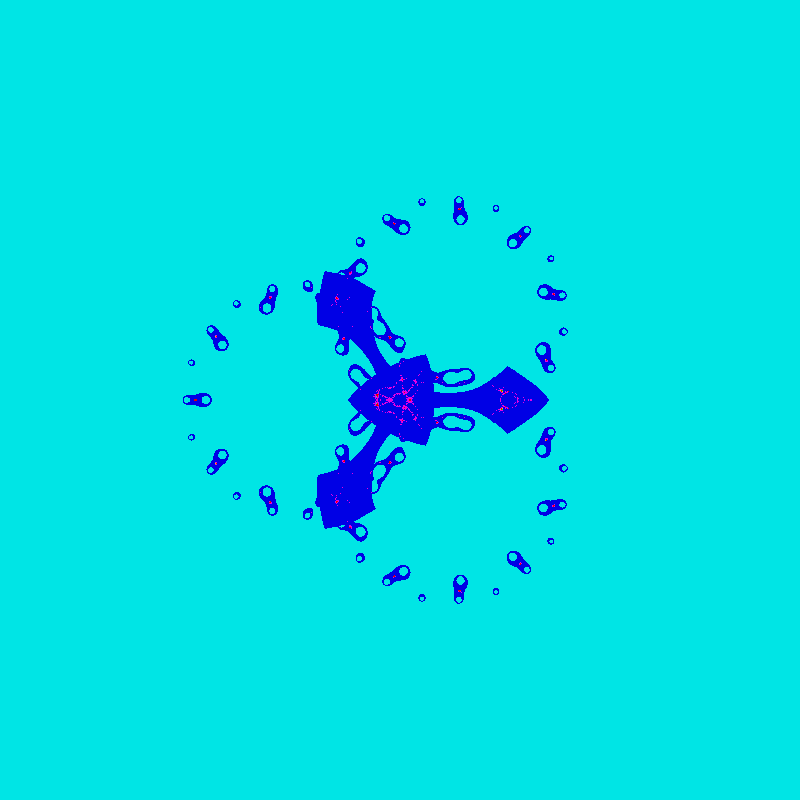}}\hfill
    \includegraphics[width=0.3\textwidth]{\detokenize{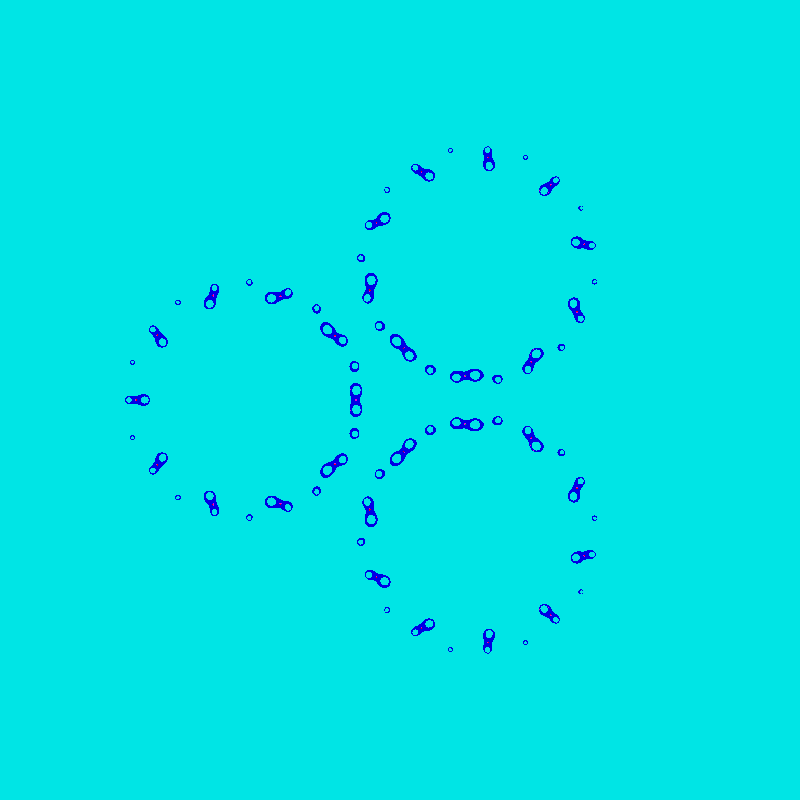}}\hfill
    \caption{Here $d=10, n=20$, and $a=0.01,0.3,0.6$ ($c$-window: $[-3.5,3.5]\times[-3.5,3.5]$). The clusters separate just before $a=0.6$.}
    \label{fig:clusters2}
\end{figure}
\end{enumerate}
\end{example}

\begin{example}
In the prior example, we had three cases in which $n$ was always larger than $d$. So, we inspected the case of $d>n$ to see if the effects were similar.  In fact, we found that changing $d$ caused more pronounced changes in the $|a|$-values needed for separation than changing $n$ did. 

Figure ~\ref{fig:diff_in_d/n} illustrates the case of the left and right pictures having $n$ and $d$ values swapped: $n=40, d=20$ on left and $n=20,d=40$ on right. Thus note the number of clusters is the same (since $(n+d)/\gcd(n,d)$ is symmetric in $n$ and $d$). The same $a=0.35$ (and same $c-$windows) are used for both pictures. But, we see the clusters separate for lower $|a|$ value in the case $d \gg n$.

\begin{figure}[!ht]
    \centering
    \includegraphics[width=0.49\textwidth]{\detokenize{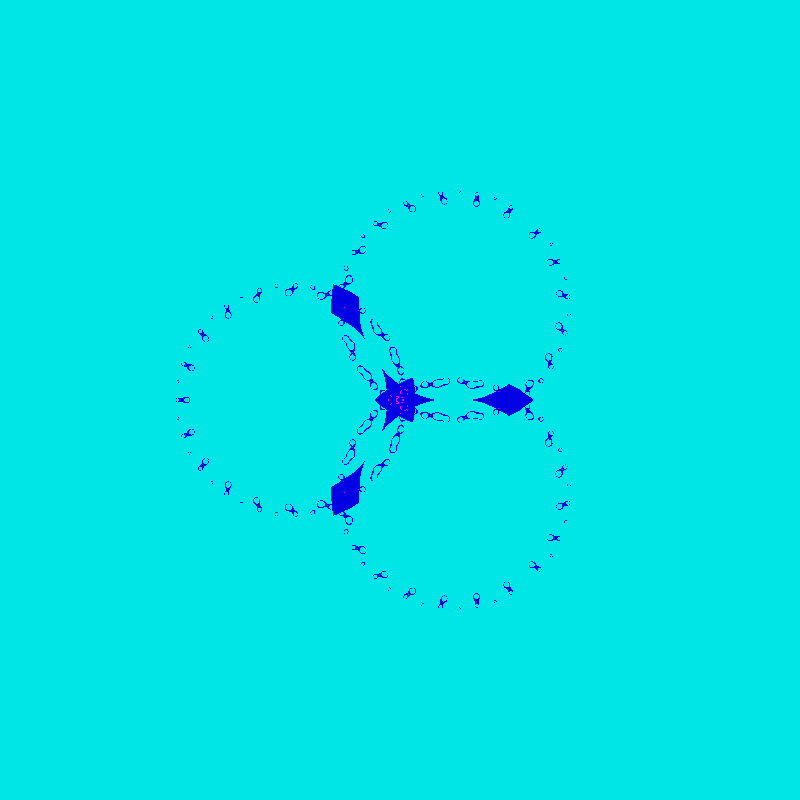}}
    \includegraphics[width=0.49\textwidth]{\detokenize{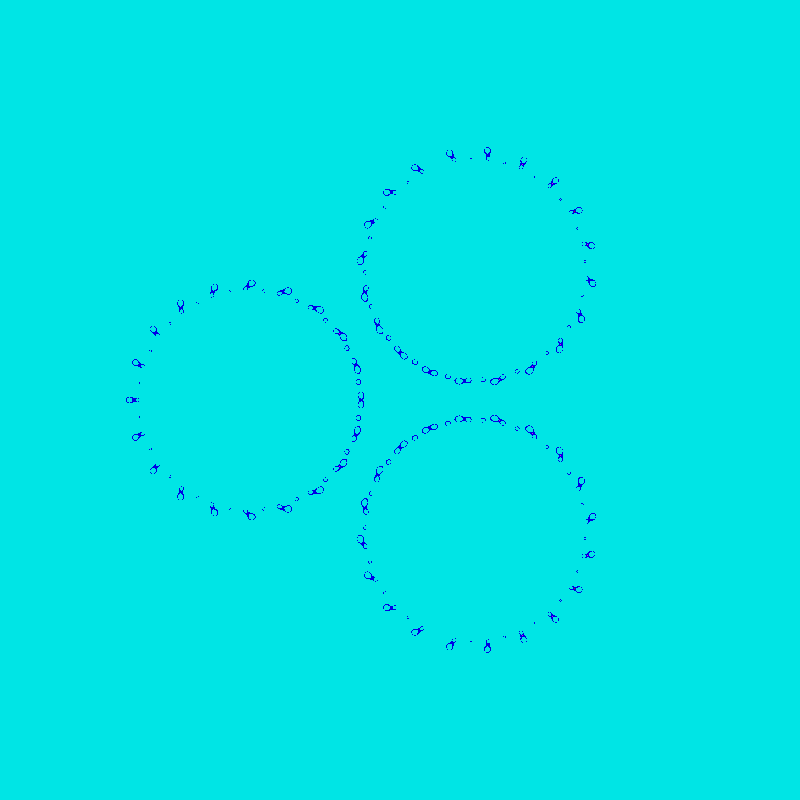}}
    \caption{$n=40$ and $d=20$ (left), $n=20$ and $d=40$ (right), $a=0.35$, ($c$-window: $[-3.5,3.5]\times[-3.5,3.5]$) (both). Note 
    the separation of clusters occurs for lower values of $|a|$ when $d \gg n$.}
    \label{fig:diff_in_d/n}
\end{figure}
    \end{example}

\subsection{Generating Parameter Space and Dynamical Plane Pictures: Interplay with our theory}

Here, we explain how to use both Dynamics Explorer and the theorems of this article to produce good dynamical and parameter images illustrating this family. 

First, for a dynamical plane, Julia set image, we have a map given by inputted $n,d,c,a$ values. The user then determines a maximum number of iterations, \ttmaxiter, and upper and lower escape radii, {\ttbig } and {\tttiny } (and the pixel width, larger means better resolution). The software iterates up to maxiter the center point of each pixel, and on each iteration tests whether the value is $>$ {\ttbig} or $<$ {\tttiny}, if so the pixel is shaded a color based on the number of iterations it took to trigger the exit condition. Pixels which never exit  up to {\ttmaxiter} are colored black. 
In the Dynamics Explorer software (as in many others) the user can choose from a few coloration options, such as rainbow or shades of one color, etc., for the color cycling. There is also a ``color density'' slider which allows the user to adjust the speed of the cycling. In most of our examples, we had the density adjusted so the colors changes at the maximum speed allowed. 

Now to connect this with the theory of this article, Proposition~\ref{prop:K-in-annulus} provides, for a given $n,d,c,a$, a theoretical value $p$ for the inner radii and  $s$ for the outer radii of an annulus guaranteed to contain a Julia set. So, you might think if you set {\tttiny} $\le p$ and {\ttbig} $\ge s$, you are essentially guaranteed to detect all of the Julia set (see Example~\ref{ex:J-annuli}).  However, for efficiency, the algorithm calculates the square norm of a pixel rather than just the norm, so what you actually need to do is set {\tttiny } and {\ttbig} in relation to $p^2$ and $s^2$. 

Choosing values of {\ttbig} and {\tttiny} which are close to $p^2,s^2$ give you the fastest algorithm, as it takes fewer iterates to exit, for escaping points. However, the user may want to increase {\ttbig} or decrease {\tttiny} as it changes which sets of pixels exit on the same iteration so it can change the colors in the image. For example, in Figure~\ref{fig:in/out-param} the darker blue star shape in the middle changes in size with the value of {\ttbig}: as you increase {\ttbig}, the annuli show up better around the edges as more pixels become the darker blue; however, the central star set of dark blue also expands and gradually swallows the clusters, leaving only a few pink speckles that are difficult to see. So, the user may adjust these values to obtain the best picture. For example, for lower $n,d$ values we tended to obtain good pictures with smaller {\ttbig} values, whereas for larger $n,d$ values a much larger {\ttbig} is required. The {\tttiny} has a smaller effect.

\begin{figure}[!ht]
    \centering
    \includegraphics[width=0.237\textwidth,trim=0 0 0 0,clip]{\detokenize{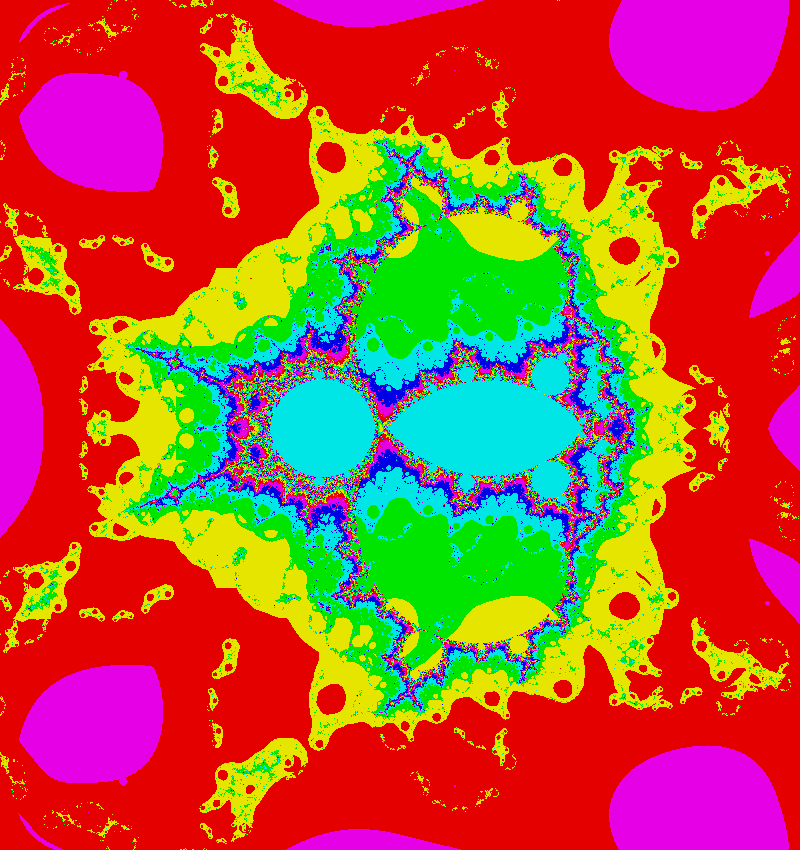}}\hspace{0.15cm}%
    \includegraphics[width=0.237\textwidth,trim=0 0 0 0,clip]{\detokenize{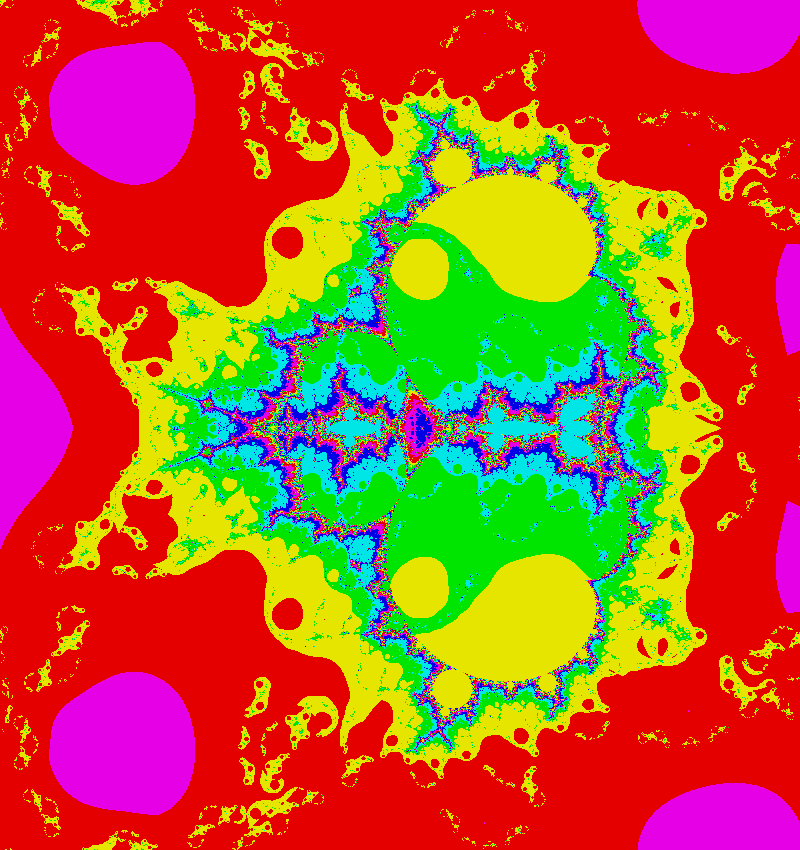}} \hspace{0.1cm}%
    \includegraphics[width=0.237\textwidth,trim=0 0 0 0,clip]{\detokenize{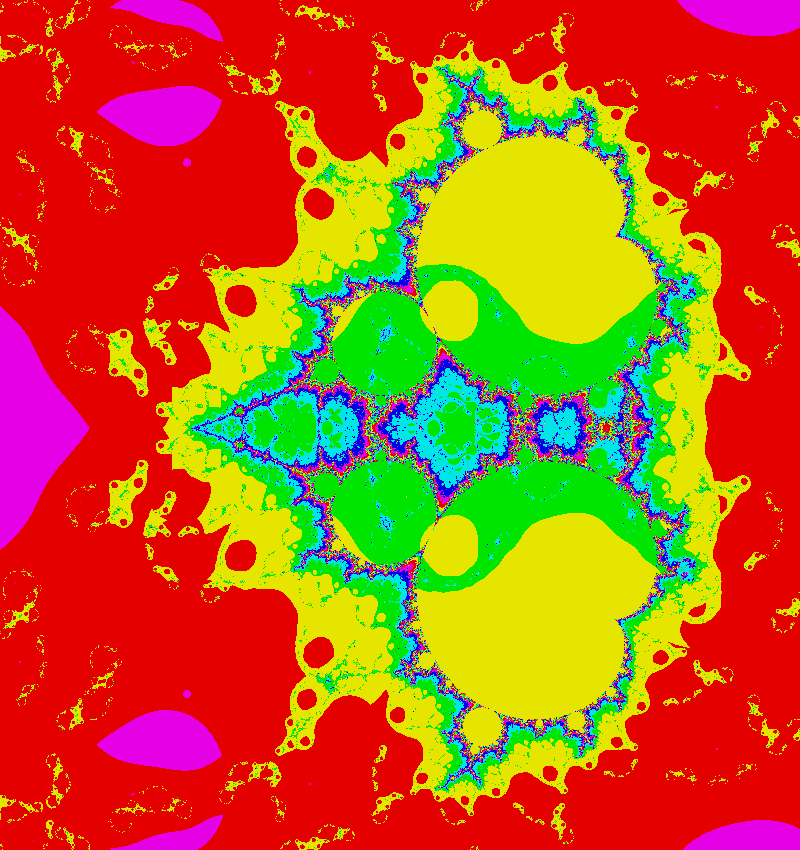}}\hspace{0.15cm}%
    \includegraphics[width=0.237\textwidth,trim=0 0 0 0,clip]{\detokenize{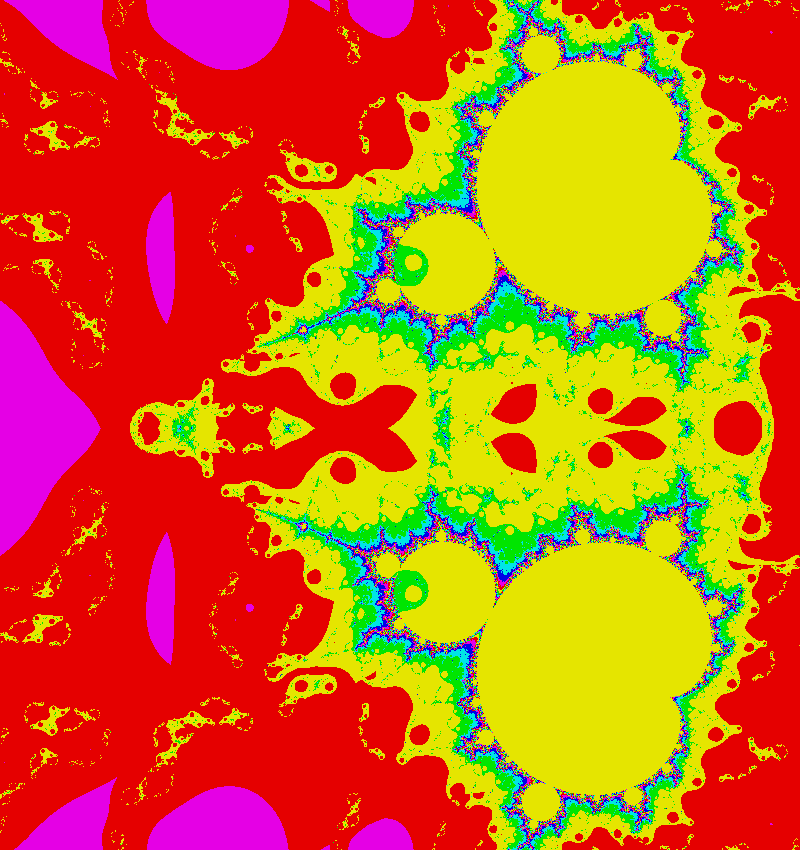}}
    \caption{$d=5,n=10$ and $a=0.26,0.265,0.27,0.28$ ($c$-window: $[0.87,1.026]\times [-0.082,0.084]$). Observe how baby Mandelbrots for different critical orbits can overlap. }
    \label{fig:babyMs-interesect}
\end{figure}

Now for parameter space pictures, remember that $\Mnda$ is the collection of all $c$-values such that at least one critical orbit is bounded. There are $\nu$ critical values. Hence, at each $c$-value, every one of the critical values' orbits is tested, and assigned black or a color based on rate of escape, if it escapes. Then, to assign just one color to the pixel, the RGB values of all $\nu$ colors are {\em averaged}. This is why, for example in Figure~\ref{fig:p_to_d} upper right, the baby Mandelbrot set is not black: only one of the three critical orbits is bounded for that set, the other two escape (each cluster is associated with one critical value and each is disjoint). 
For a more complicated scenario, for smaller $a$ when the clusters intersect, the colors for multiple critical values are rapidly changing in a small area - observe Figure~\ref{fig:babyMs-interesect} showing a series of two magnified baby Mandelbrots intersecting as $a$ varies. Observe how the color pattern for each critical orbit interferes with the other.  

Finally, for the value $T$ of Lemma~\ref{lem:MGMdisk-for-Mdna}, recall we don't provide a calculation, we prove that the boundedness locus is in some bounded disk $\DD_T$ for some $T$, it's an existence proof. But, when we draw a picture of the boundedness locus, and zoom out sufficiently to convince ourselves we have captured all of it, we can essentially estimate $T$ from the image generated. Thus, the theory and the computer generated images each contribute to our understanding of the dynamics in their own way.

\section{Future work}
\label{sec:future}

Finally, we suggest here a few avenues for future work. 
 
\begin{figure}[!ht] 
    \centering
    \includegraphics[width=0.29\textwidth]{\detokenize{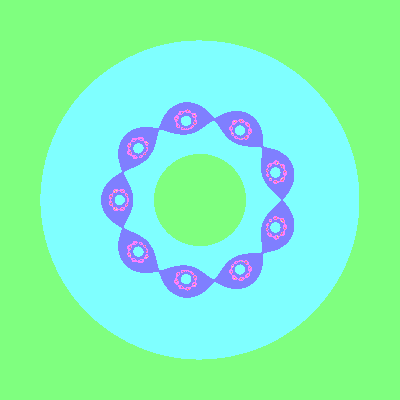}}\hfill
    \includegraphics[width=0.29\textwidth]{\detokenize{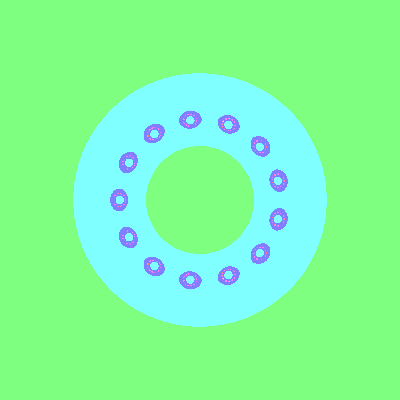}}\hfill
    \includegraphics[width=0.29\textwidth]{\detokenize{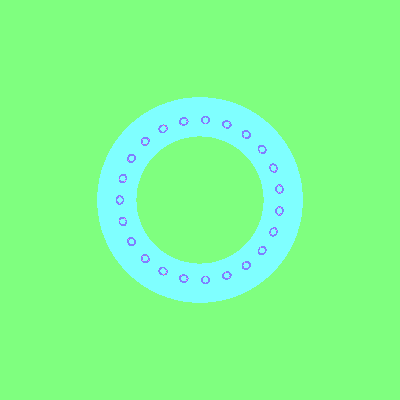}}\hfill
    \caption{$c=0.001,a=1,(n,d)=(4,5),(6,7),(11,12)$, and we draw  the window $[-2.5,2.5]\times[-2.5,2.5]$. This $n_k,d_k$ are consecutive, so relatively prime, hence $\gcd=1$ for all $k$.}
    \label{fig:n-d-rel-prime}
\end{figure}

\begin{figure}[!ht] 
    \centering
    \includegraphics[width=0.3\textwidth]{\detokenize{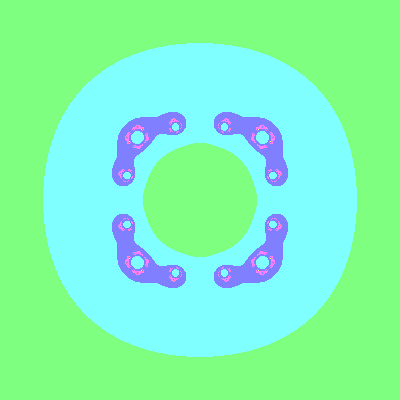}}\hfill
    \includegraphics[width=0.3\textwidth]{\detokenize{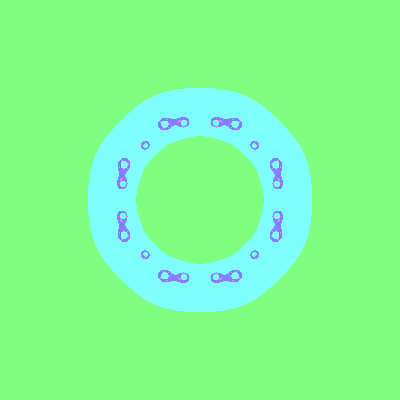}}\hfill
    \includegraphics[width=0.3\textwidth]{\detokenize{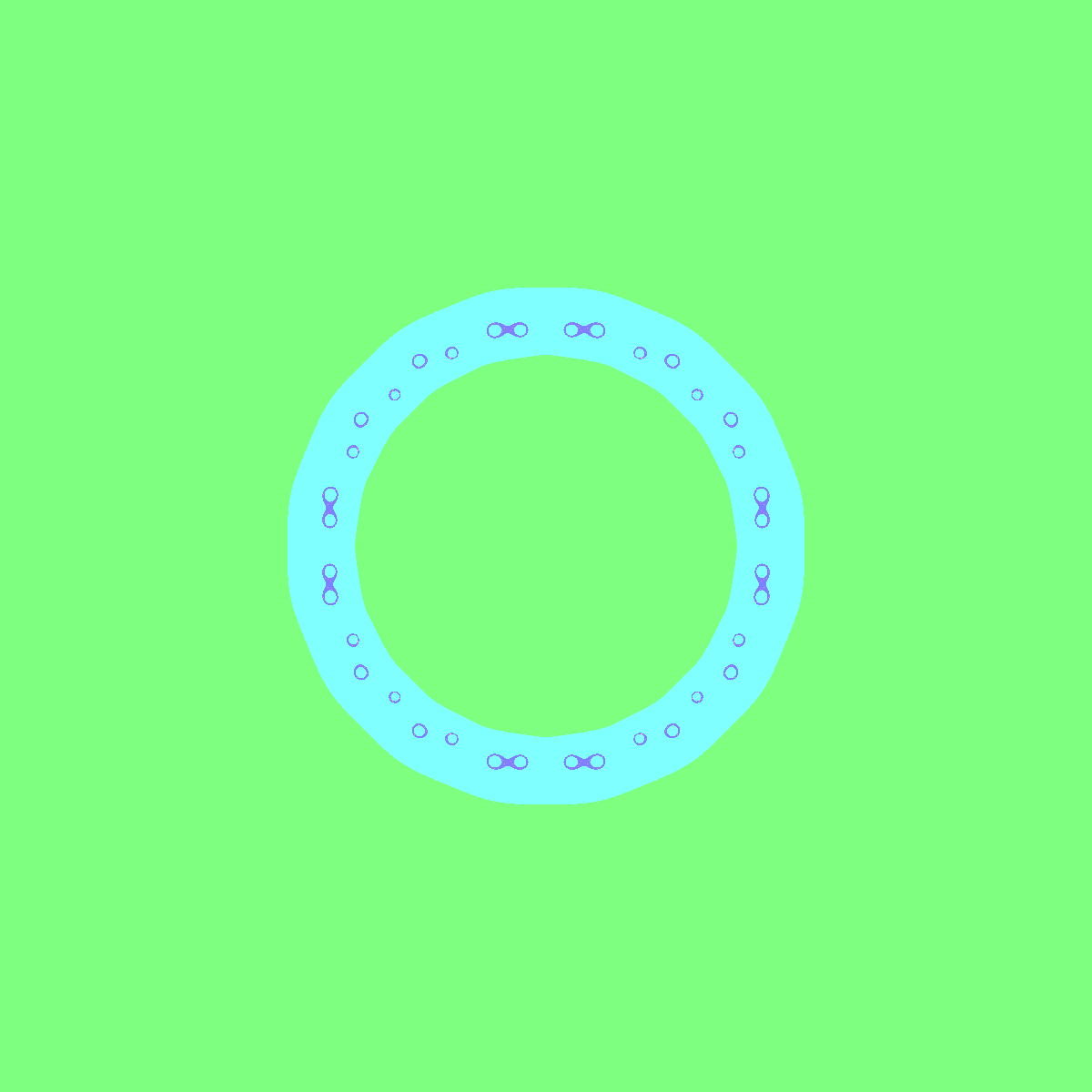}}\hfill
    \caption{$c=1,a=1,(n,d)=(4,8),(8,12),(16,20)$ and we draw the window 
    $[-2.5,2.5]\times[-2.5,2.5]$. Here $\gcd(n,d)=4$.}
    \label{fig:gcd-4-n-d-grow}
\end{figure}

 In Theorem~\ref{thm:in+outJ1}, we first show that the limit of the Julia sets always lies in $S^1$; an open question is when precisely there is equality in the limit, and whether there is a limiting shape in some cases when $\ell_k \nrightarrow \infty$. For example, if $\ell_k=\gcd(n_k,d_k)=1$ for all $k$, it appears in some cases that the Julia set tends to fill up $S^1$, see Figure~\ref{fig:n-d-rel-prime}. Figure~\ref{fig:gcd-4-n-d-grow} is a case in which the $\gcd$ is fixed (at $4$). 

 Similarly, one could ask whether in Theorem~\ref{thm:parameter_n_d_to_inf} we always need $\ell = \gcd(n,d)\to\infty$ as $n,d\to\infty$ in order to have equality in the limit? 

In \cite{BoydSchulz}, the authors study geometric limits for $z\mapsto z^t+c$, defined in terms of the principal branch of the logarithm. We have begun an examination of that setting and expect to describe a generalization of our results to the real, and potentially even complex, exponent case, in a sequel to this article.

Another avenue would be to generalize or translate these results to other slices of parameter space (such as the $a$-planes) or to study the boundedness locus in $\CC^2$ in general. In the $a$-plane slice study, the calculations are less straightforward, since $c$ is merely added to the map but $a$ is involved in the $\frac{a}{z^d}$ term.

\bibliographystyle{plain}

\begin{thebibliography}{10}

\bibitem{Beardon}
Alan~F. Beardon.
\newblock {\em Iteration of {R}ational {F}unctions}, volume 132 of {\em Graduate Texts in Mathematics}.
\newblock Springer-Verlag, New York, 1991.
\newblock Complex analytic dynamical systems.

\bibitem{DevaneyRussell}
Paul Blanchard, Robert~L. Devaney, Antonio Garijo, and Elizabeth~D. Russell.
\newblock A generalized version of the {M}c{M}ullen domain.
\newblock {\em Internat. J. Bifur. Chaos Appl. Sci. Engrg.}, 18(8):2309--2318, 2008.

\bibitem{detool}
Brian Boyd and Suzanne Boyd.
\newblock {D}ynamics {E}xplorer, at sourceforge.
\newblock [http://sourceforge.net/projects/detool].

\bibitem{BoydHoeppner1}
Suzanne Boyd and Matthew Hoeppner.
\newblock Baby {M}andelbrot sets and spines in some one-dimensional subspaces of the parameter space for generalized {M}c{M}ullen maps.
\newblock {\em Qual. Theory Dyn. Syst.}, 24(4):Paper No. 157, 39, 2025.

\bibitem{BoydMitchell}
Suzanne Boyd and Alexander~J. Mitchell.
\newblock The boundedness locus and baby {M}andelbrot sets for some generalized {M}c{M}ullen maps.
\newblock {\em Internat. J. Bifur. Chaos Appl. Sci. Engrg.}, 33(9):Paper No. 2350107, 23, 2023.

\bibitem{BoydSchulz}
Suzanne Boyd and Michael~J. Schulz.
\newblock Geometric limits of mandelbrot and julia sets under degree growth.
\newblock {\em International Journal of Bifurcations and Chaos}, 22(12), 2012.

\bibitem{Kaschner2020}
Micah Brame and Scott Kaschner.
\newblock Geometric limits of {J}ulia sets for sums of power maps and polynomials.
\newblock {\em Complex Anal. Synerg.}, 6(3):Paper No. 22, 8, 2020.

\bibitem{DevaneyHalos}
Robert~L. Devaney.
\newblock Baby {M}andelbrot sets adorned with halos in families of rational maps.
\newblock In {\em Complex dynamics}, volume 396 of {\em Contemp. Math.}, pages 37--50. Amer. Math. Soc., Providence, RI, 2006.

\bibitem{DevaneySurvey2013}
Robert~L. Devaney.
\newblock Singular perturbations of complex polynomials.
\newblock {\em Bull. Amer. Math. Soc. (N.S.)}, 50(3):391--429, 2013.

\bibitem{devgar}
Robert~L. Devaney and Antonio Garijo.
\newblock Julia sets converging to the unit disk.
\newblock {\em Proc. Amer. Math. Soc.}, 136(3):981--988, 2008.

\bibitem{Kaschner2015}
Scott~R. Kaschner, Reaper Romero, and David Simmons.
\newblock Geometric limits of {J}ulia sets of maps {$z^n+\exp(2\pi i\theta)$} as {$n\rightarrow\infty$}.
\newblock {\em Internat. J. Bifur. Chaos Appl. Sci. Engrg.}, 25(8):1530021, 8, 2015.

\bibitem{Murali-PME}
Murali Meyer.
\newblock Morsels of {M}andelbrot sets.
\newblock {\em Pi Mu Epsilon J.}, 16(3):165--171, 2025.

\end{thebibliography}

\end{document}